\documentclass[preprint,11pt]{elsarticle}
\usepackage[latin1]{inputenc}
\usepackage{a4}
\usepackage{color}
\usepackage{amsmath}
\usepackage{amsfonts}
\usepackage{amsthm}
\usepackage{amssymb}
\usepackage{bm}
\usepackage{bbm}   
\usepackage{latexsym}
\usepackage{stmaryrd}
\usepackage{epsfig}
\usepackage{url}
\usepackage{graphics}
\usepackage{graphicx}
\usepackage{fullpage}
\usepackage{subfigure}
\usepackage{appendix}
\usepackage{color}
\usepackage{xcolor}
\usepackage{multirow}
\usepackage{hhline}
\usepackage[overload]{empheq}
\usepackage[top=2.50cm,left=2.50cm,right=2.50cm,bottom=3.15cm]{geometry}

\usepackage[T3,T1]{fontenc}
\usepackage{abraces}

\newcommand{\R}{\mathbb R}

\def\P{{\mathbb P}}

\def\1{{\mathbbm 1}}

\newcommand{\ub}[2]{\displaystyle{\underbrace{#1}_{#2}}}

\newcommand{\dd}{\mathrm{d}}

\newtheorem{theoreme}{Theorem}[section]

\newtheorem{remark}[theoreme]{Remark}

\newcommand{\Min}{\displaystyle\min}
\newcommand{\Max}{\displaystyle\max}

\newcommand{\xory}{{x\backslash y}}
\newcommand{\Int}{\displaystyle\int}
\newcommand{\Oint}{\displaystyle\oint}
\newcommand{\Sum}{\displaystyle\sum}
\newcommand{\Lim}[2]{\displaystyle\lim\limits_{#1\,\rightarrow\, #2}}
\newcommand{\Frac}{\displaystyle\frac}

\newcommand{\ov}{\overline}
\newcommand{\un}{\underline}
\newcommand{\wt}{\widetilde}

\newcommand{\wh}[1]{\widehat{#1}}

\newcommand{\tra}{^{\text{t}}}

\newcommand{\qqquad}{\qquad\qquad}
\newcommand{\lh}{\hspace*{-2mm}}

\newcommand{\pds}{\centerdot}
\newcommand{\demi}{\frac{1}{2}}

\newcommand{\Dt}{\Delta t}

\newcommand{\gradx}[1]{\nabla_x{#1}}

\newcommand{\divx}{\nabla_x\pds}

\newcommand{\dt}[1]{\Frac{\partial {\,#1}}{\partial t}}

\newcommand{\ddt}[1]{\Frac{\mathrm{d} {\,#1}}{\mathrm{d} t}}

\newcommand{\vdt}{\partial_t}
\newcommand{\vdx}{\partial_x}

\newcommand{\vdy}{\partial_y}

\newcommand{\poly}[1]{\P^{\,#1}}
\newcommand{\mrm}[1]{\mathrm{#1}}
\newcommand{\mc}[1]{\mathcal{#1}}
\newcommand{\mf}[1]{\mathfrak{#1}}

\newcommand{\Mat}[1]{\mathsf{#1}}
\newcommand{\inv}[1]{\frac{1}{#1}}
\newcommand{\Inv}[1]{\Frac{1}{#1}}

\newcommand{\bk}[1]{\llbracket{#1}\rrbracket}
\newcommand{\bs}[1]{\boldsymbol{#1}}

\def\veps{\varepsilon}

\def\<{\left<}
\def\>{\right>}
\def\({\left(}
\def\){\right)}

\DeclareSymbolFont{tipa}{T3}{cmr}{m}{n}
\DeclareMathAccent{\invbreve}{\mathalpha}{tipa}{16}
\newcommand{\ibreve}[1]{\invbreve{#1\hspace*{0.2mm}}}

\makeatletter
\def\widebreve{\mathpalette\wide@breve}
\def\wide@breve#1#2{\sbox\z@{$#1#2$}%
     \mathop{\vbox{\m@th\ialign{##\crcr
\kern0.08em\brevefill#1{0.8\wd\z@}\crcr\noalign{\nointerlineskip}%
                    $\hss#1#2\hss$\crcr}}}\limits}
\def\brevefill#1#2{$\m@th\sbox\tw@{$#1($}%
  \hss\resizebox{#2}{\wd\tw@}{\rotatebox[origin=c]{90}{\upshape(}}\hss$}
\makeatletter

\definecolor{dark-gray}{gray}{0.35}
\newcommand{\Dgray}[1]{\textcolor{dark-gray}{#1}}

\newcommand{\ext}{pdf}
\graphicspath{{FIGURES/\ext/}}

\newcommand{\scheme}{local subcell monolithic DG/FV }
\newcommand{\Scheme}{Local subcell monolithic DG/FV }

\journal{Journal of Computational Physics}
\begin{document} 
\begin{frontmatter}

\title{\Scheme convex property preserving scheme on unstructured grids and entropy consideration}

\author[imag]{Fran\c{c}ois Vilar}
\ead{francois.vilar@umontpellier.fr}
\address[imag]{IMAG, Univ Montpellier, CNRS, Montpellier, France}

\begin{abstract}
  This article aims at presenting a new local subcell monolithic Discontinuous-Galerkin/Finite-Volume (DG/FV) convex property preserving scheme solving system of conservation laws on 2D unstructured grids. This is known that DG method needs some sort of nonlinear limiting to avoid spurious oscillations or nonlinear instabilities which may lead to the crash of the code. The main idea motivating the present work is to improve the robustness of DG schemes, while preserving as much as possible its high accuracy and very precise subcell resolution. To do so, a convex blending of high-order DG and first-order FV scheme will be locally performed, at the subcell scale, where it is needed. To this end, by means of the theory developed in \cite{vilar_aplsc_1D,vilar_aplsc_2D}, we first recall that it is possible to rewrite DG scheme as a subcell FV scheme on a subgrid provided with some specific numerical fluxes referred to as DG reconstructed fluxes. Then, the monolithic DG/FV method will be defined as following: each face of each subcell will be assigned with two fluxes, a 1st-order FV one and a high-order reconstructed one, that in the end will be blended in a convex way. The goal is then to determine, through analysis, optimal blending coefficients to achieve the desire properties. Numerical results on various type problems will be presented to assess the very good performance of the design method.\\
  A particular emphasis will be put on entropy consideration. By means of this subcell monolithic framework, we will attempt to address the following questions: is this possible through this monolithic framework to ensure any entropy stability? what do we mean by entropy stability? What is the cost of such constraints? Is this absolutely needed while aiming for high-order accuracy?
\end{abstract}

\end{frontmatter}

\newpage
\section{Introduction}
\label{intro}

This paper is concerned with solving system of conservation laws, and it is well known that hyperbolic partial differential equations frequently lead to discontinuous weak solutions within a finite time frame, posing significant challenges for numerical simulations. These challenges revolve around handling discontinuities, ensuring accuracy, and maintaining solutions within an admissibility set, such as guaranteeing positive density and internal energy in gas dynamics. Addressing these constraints simultaneously is particularly difficult because they often conflict with one another, demanding sophisticated approaches in the design of numerical methods for hyperbolic problems. A large number of numerical scheme have been developed these past fifty years to achieve such goal, and one which has particularly stood out is the Discontinuous Galerkin (DG) method. This scheme, initially introduced by Reed and Hill for neutron transport \cite{Reed}, has become one of the most widely used numerical schemes, particularly in computational fluid dynamics. Significant advancements by Cockburn and Shu in a series of seminal papers, see for instance \cite{Cockburn_lcs5} and the references within, have propelled DG methods forward. These methods theoretically allow for achieving any arbitrary order of accuracy while maintaining a compact stencil, and they exhibit desirable properties such as $L_2$ stability and $hp$-adaptivity. The DG scheme is renowned for its high accuracy and precise subcell resolution, even demonstrating superconvergence in some cases. However, robustness is a critical concern alongside accuracy. High-order DG schemes are known to produce spurious oscillations in the presence of discontinuities and possibly non-physical solutions (e.g., negative density or pressure in gas dynamics), which may lead to nonlinear instability or code crashes. Therefore, stabilization techniques are essential. This fundamental issue has been extensively tackled in the past. There is thus a vast literature on that topic, among which \cite{bisw,burb01,kriv07,MYang,MLP_unstruct,Kuzmin09,Li_vertex_based_lim}. For a lot more detailed description of the state of the art limiters, we refer to \cite{vilar_aplsc_1D} and the references within.\\

These past fifteen years, great progresses have been made in the direction of stabilizing and improving the robustness of high-order DG. And to do so, two main properties were under study: convex set preserving also referred to as Invariant Domain Preserving (IDP) and high-order entropy stability. In the former, the goal is to ensure that the numerical solution remains, at all time, in a convex admissible set. This property particularly permits to guarantee global maximum principle for Scalar Conservation Laws (SCL) and for instance positivity of the density and pressure in the Euler system case. The articles on this subject have flourished in recent years. Although not exhaustive, it worth mentioning \cite{zshu1,zshu5} and \cite{IDP_renac} where some polynomial limiters have been developed to ensured such property, as well as \cite{bound_preserving_shu} where a new framework, refereed to as geometric quasilinearization, has been introduced to also develop bound-preserving numerical methods. Another wide family of schemes also concerned with this convex property preserving issue is the one gathering Flux-Corrected Transport (FCT) techniques, Algebraic Flux Corrections (AFC) and convex limiting schemes, see for instance \cite{Boris_FCT,Zalesak_FCT,IDP_guermond_2016,Kuzmin_FCT_limiting_2017,IDP_guermond_2018,IDP_guermond_2019,Kuzmin_monolithic_2020,monolithic_hajduk_2021,Kuzmin_FCT_limiting_2022}. Most of these aforementioned methods share a similar philosophy, meaning blending high-order and low-order fluxes, operators or schemes to ensure the preservation of convex properties or more generally to be IDP.\\

Now, regarding high-order entropy stable schemes, a new class of numerical methods has recently extensively grown in popularity, see for instance \cite{SBP_carpenter,DG_SBP_gassner,entropy_carpenter,SBP_gassner,SBP_shu,DG_SBP_chan,DG_SBP_renac}. Those schemes, generally referred to as entropy conservative/stable DG Spectral Element Method (DGSEM), were initially developed in the context of finite difference Summation-By-Parts (SBP) operators and Simultaneous Approximation Term (SAT) by T. Fisher and M. Carpenter in their seminal paper \cite{SBP_carpenter}. In the context of DG, they rely on the use of particular quadrature rules inducing a mass lumping type diagonalization of the mass matrix and specific collocation of the flux in order to exhibit discrete SBP properties. Then, a substitution of the flux collocation values by a combination of entropy conservative numerical fluxes ensure an entropy conservation or stability while preserving the high-order accuracy. While some SBP operators and entropy stable DGSEM scheme have been successfully extended to simplex meshes, see for instance \cite{unstructured_SBP_16,SBP_shu}, this family of numerical methods are generally restrained to one-dimension in space or tensor-product multi-dimensional grids.\\

A third direction which has been particularly embraced this past decade and has shown some of the most promising developments is subcell techniques. Here, the idea is to subdivide the bad cells, and to adopt a special procedure with the hope of curing the negative aspects of the original scheme. Some examples of this strategy can be found in \cite{peraire_2012,munz_subFV_2014,dumbser_subFV_lim_tri}. For example, in \cite{peraire_2012}, the authors use a convex combination between high-order DG schemes and first-order Finite-Volumes (FV) on a subgrid, allowing them to retain the very high accurate resolution of DG in smooth areas and ensuring the scheme's robustness in the presence of shocks. Similarly, in \cite{munz_subFV_2014,dumbser_subFV_lim_tri}, after having detected the troubled zones, cells are then subdivided into subcells and a robust first-order FV scheme, or alternatively other robust scheme (second-order TVD FV scheme, WENO scheme, \dots), is performed on the subgrid in troubled cells. Let us emphasize that these subcell techniques offer several advantages. They preserve the high accuracy of DG schemes in smooth regions by applying corrections only where necessary. This local modification approach ensures that the majority of the grid remains unaffected by the stabilization process, maintaining computational efficiency and accuracy. Let us however emphasize that if a correction is used at the subgrid level, all the subcells contained in a bad cell are generally impacted. Since one of the main advantage of high-order scheme is to be able to use coarse grids while still being very precise, one can see that there is a waste of information here, as well as unnecessary computational effort made. This is particularly the case in the vicinity of discontinuities since the polynomials are globally modified. This problem was addressed in the one-dimensional case in \cite{vilar_aplsc_1D} and in the two-dimensional unstructured case in \cite{vilar_aplsc_2D}. This new technique relies on the reformulation of DG schemes as a FV-like scheme defined on a subgrid, through the definition of particular fluxes referred to as reconstructed fluxes. Then, after computing a DG candidate solution and check if this solution is admissible, one returns if needed to the previous time step and correct locally, at the subcell scale, the numerical solution. In the subcells where the solution was detected as bad, one substitutes the DG reconstructed flux on the subcell boundaries by a robust first-order numerical flux. And for subcells detected as admissible, one keeps the high-order reconstructed flux which allows to retain conservativity as well as  the very high accurate resolution of DG schemes. Consequently, only the solution inside troubled subcells and their first neighbors will have to be recomputed. Elsewhere, the solution remains unchanged. This correction procedure is then extremely local, and has proved in different contexts its high capability to ensure a stable and robust behavior while maintaining the very high accuracy of DG schemes, see \cite{vilar_aplsc_1D,vilar_aplsc_NSW_1D,vilar_aplsc_2D,vilar_aplsc_NSW_1D_fixed}.\\

Finally, in recent years the interest in combining these three family of schemes and techniques, namely convex limiting methods, high-order entropy stable schemes and subcell techniques has grown tremendously. Have therefore emerged new methods, as \cite{Kuzmin_subcell_flux_limiting,Pazner_2021,subcell_shock_rueda,subcell_limiting_rueda,Kuzmin_Gassner_monolithic_2024,Chan_knapsack}, which combine, at the subcell scale, high-order and low-order schemes to ensure different properties on the numerical solution, while trying to preserve accuracy. For instance in \cite{Kuzmin_Gassner_monolithic_2024}, the authors develop a subcell monolithic scheme blending high-order DGSEM based on Gauss-Lobatto quadrature points and a first-order FV scheme. Different conditions on the blending coefficients, corresponding to different properties as positivity or local maximum principle for the concern of spurious oscillations, are presented. Similarly, in \cite{Chan_knapsack}, a monolithic Gauss-Lobatto DGSEM and FV scheme is also presented but along with a particular blending procedure, based of the resolution of a continuous Knapsack optimization problem, enabling the preservation of the high-order accuracy of the scheme while ensuring a cell entropy inequality. Let us underline that the numerical scheme presented in the present article falls also in this category. Indeed, the aim of this paper is to introduced a new monolithic scheme in which DG and first-order FV methods will be blended, locally at the subcell scale, to ensure any convex property as well as different entropy stabilities, while trying to preserve as much as possible the high accuracy of DG schemes. Let us emphasize that the monolithic schemes presented in \cite{Kuzmin_Gassner_monolithic_2024,Chan_knapsack}, because being based on Gauss-Lobatto solution representation and flux collocation, are limited to one-dimension in space, or by tensor-product extension to multi-dimensions on Cartesian grids. And because the theory developed in \cite{vilar_aplsc_2D}, namely reformulating DG scheme into a FV-like one, is very general in the sens that it can be extended to any dimension and any type of grids and cell subdivision, we aim here at presenting a monolithic scheme applicable to generic polytopal meshes. Let us emphasize however that, if the whole theory and scheme are indeed developed on any type of grid, numerical results are only shown on triangular meshes. The practical implementation of those schemes on generic polygonal grids is still an ongoing project.\\

Now, to present this \scheme scheme, the remainder of this paper is organized as follows: we recall in Section~\ref{sect_DG_as_FV} how unstructured grid DG scheme can be reinterpreted as a subgrid FV-like scheme, through the definition of particular fluxes that we referred to as reconstructed fluxes. While this part mainly relies on Section~2 of our previous article, \cite{vilar_aplsc_2D}, this theoretical analysis will be taken further here as the case where the number of subcells does not fit the dimension of the functional space will be addressed. Then, the \scheme scheme will be introduced in Section~\ref{sect_monolithic}. Practically, each face of each subcell will be assigned two fluxes, one reconstructed flux giving the equivalency with a high-order DG scheme and one first-order FV numerical flux. These two fluxes will then be blended in a convex manner through a blending coefficient between zero and one. The goal is now to determine, through analysis, the optimal coefficients to reach the desired properties while trying to maintain as much as possible the high accuracy of the scheme. Following this, we present in Section~\ref{sect_entropy} different definitions for the blending coefficients to reach different types of entropy stability, meaning from a discrete subcell entropy inequality for any entropy to a semi-discrete cell entropy inequality for a given entropy. Only the latter one will proved to allow high-order accuracy preservation. Numerical results and a preliminary conclusion on those entropy stabilities will be given at the end of this section. Finally, section~\ref{sect_GLMP} is devoted to the definition of the blending coefficients ensuring different maximum principles, and by this make the monolithic scheme IDP. These theoretical parts as well as the different numerical results will be presented for both SCL and the Euler compressible gas dynamics system.

\newpage
\section{DG scheme reformulation}
\label{sect_DG_as_FV}

This section is devoted to the demonstration of the equivalency between DG schemes and a FV-like method on a subgrid provided the definition of particular fluxes. This theoretical part mainly relies on Section~2 of our previous article, \cite{vilar_aplsc_2D}. Consequently, only the essential ingredients of such reformulation will be recalled at this time. Let us yet note that this theoretical analysis will be taken further here than in \cite{vilar_aplsc_2D}, as the case where the number of subcells does not fit the dimension of the functional space will be addressed.
To remain as simple as possible, two-dimensional SCL will be considered in this section. The system extension is perfectly straightforward. Let then $u=u(\bs{x},t)$, for $\bs{x} \in \omega\subset \mathbb{R}^2$ and $t \in [0,T]$, be the solution of the following problem

\begin{subequations}
\label{lcs_eq2D}
\begin{alignat}{2}[left = \empheqlbrace\,]
&\vdt{\,u}(\bs{x},t)+\divx{\bs{F}\(u(\bs{x},t)\)}=0, \qqquad& (\bs{x},t)\in\,\omega\times[0,T], \label{lcs1}\\
&u(\bs{x},0)=u_0(\bs{x}), &\bs{x}\in\,\omega, \label{lcs2}
\end{alignat}\\[-4mm]
\end{subequations}

where $u_0$ is the initial data and $\bs{F}(u)$ the flux function. For the subsequent discretization, let us introduce the following notation. $\{\omega_c\}_c$ would be a generic partition of the domain $\omega$ into non-overlapping cells, with $|\omega_c|$ being the size of $\omega_c$. We also partition the time domain in intermediate times $(t^n)_n$ with $\Delta t^n=t^{n+1}-t^n$ the $n^{th}$ time step. In the DG frame, the numerical solution is considered piecewise polynomial over the domain, and hence developed on each cell onto $\P^{\,k} (\omega_c)$, the set of polynomials of degree up to $k$ defined on cell $\omega_c$. This space approximation theoretically leads to a $(k+1)^{th}$ space order accurate scheme. Let $u_h^c=\sum_{m=1}^{N_k}u_m^c(t)\,\sigma_m^c(\bs{x})$ be the restriction of $u_h$, the piecewise polynomial approximation of the solution $u$, over the cell $\omega_c$, where the $u_m^c$ are the $N_k$ successive components of $u_h^c$ over the polynomial basis $\{\sigma_m^c\}_m$. We recall that in the two-dimensional case $N_k=\frac{(k+1) (k+2)} {2}$. The coefficients $u_m^c$ are the solution moments to be computed. To this end, by means of the weak formulation of equation \eqref{lcs1} on $\omega_c$, restricting the solution functional space and the test function space to $\P^{\,k}(\omega_c)$ and then substituting the solution $u$ by its approximated polynomial counterpart $u_h^c$, one gets

\begin{align}
  \label{DG_2D_1}
  \Int_{\omega_c} \dt{u_h^c}\,\psi\,\dd V=\Int_{\omega_c} \bs{F}(u_h^c)\pds\gradx{\psi}\,\dd V-\Int_{\partial \omega_c}\hspace*{-2mm}\psi\;\mc{F}_n\,\dd S, \qquad \forall \,\psi\in\P^{\,k}(\omega_c).
\end{align}

The DG numerical solution $u_h^c$ is then the unique polynomial function defined in $\,\P^{\,k}(\omega_c)$ satisfying equation \eqref{DG_2D_1} for all function $\psi\in\P^{\,k}(\omega_c)$. In \eqref{DG_2D_1}, the numerical flux function $\mc{F}_n$, in addition to ensure the scheme conservation, is the cornerstone of any FV or DG scheme regarding fundamental considerations as stability, positivity and entropy among others. In this context, this numerical flux is defined as a function of the two states on the left and right of each interface, $\mc{F}_n=\mc{F}\(u^-,u^+,\bs{n}\)$, with $u^-=\Lim{\epsilon}{0^+}u_h^c(\bs{x}_i-\epsilon\,\bs{n},\,t)$ and $u^+=\Lim{\epsilon}{0^+}u_h^v(\bs{x}_i+\epsilon\,\bs{n},\,t)$, where $\omega_v$ is a face neighboring cell of $\omega_c$, while $\bs{x}_i$ and $\bs{n}$ respectively stand for a point and the unit outward normal of the separating interface. From now on, we refer to the set containing the face neighboring cells of $\omega_c$ as $\mc{V}_c$. The numerical flux function is generally obtained through the resolution of an exact or approximated Riemann problem. In this context of SCL, we make use of the following well-known general numerical flux definition
\begin{align}
  \label{num_flux}
  \mc{F}(u^-,u^+,\bs{n})=\Frac{\big(\bs{F}(u^-)+\bs{F}(u^+)\big)}{2}\pds\bs{n}-\Frac{\gamma(u^-,u^+,\bs{n})}{2}\,(u^+-u^-).
\end{align}
Under condition $\gamma(u^-,u^+,\bs{n})\geq \max_{w\,\in\, I(u^-,u^+)}\big(|\bs{F}'(w)\pds\bs{n}|\big)$, where $I(a,b)=[\min(a,b),\max(a,b)]$, the numerical flux \eqref{num_flux} is nothing but an E-Flux, \cite{Osher_Eflux,Tadmor_Escheme}. A FV scheme relying on such a numerical flux will be positivity-preserving and ensure a discrete entropy inequality for any entropy. Let us emphasize that for systems, the whole theory and scheme development that will follow can be easily extended to classical numerical fluxes, as HLL and HLL-C.\\

Now, taking in \eqref{DG_2D_1} the test function $\psi$ among the polynomial basis functions leads to the following linear system allowing the calculation of the solution moments $u_m^c$
\begin{align}
  \label{DG_2D_2}
  \Sum_{m=1}^{N_k}\ddt{u_m^c}\,\Int_{\omega_c} \hspace*{-2mm}\sigma^c_m\,\sigma^c_p\,\dd V=\Int_{\omega_c} \bs{F}(u_h^c)\pds\gradx{\sigma^c_p}\,\dd V-\Int_{\partial \omega_c}\hspace*{-2mm}\sigma^c_p\;\mc{F}_n\,\dd S, \qquad \forall \,p\in\,\llbracket 1, N_k\rrbracket.
\end{align}
Terms $\int_{\omega_c} \bs{F}(u_h^c)\pds\gradx{\sigma^c_p}\,\dd V$ and $\int_{\partial \omega_c}\sigma^c_p\;\mc{F}_n\,\dd S$ are respectively referred to as volume and surface integrals.

\begin{remark}
  \label{quadrature}
  Let us emphasize that these volume and surface integrals have to be computed in practice. Considering complex SCL with non-polynomial flux or non-linear systems as the Euler compressible gas dynamics one with non-convex equation of state for instance, exact integration may be difficult nay impossible. Generally, people either use quadrature rules, as originally introduced in \cite{Cockburn_lcs2}, or a collocation of the flux, as it is done in nodal DG \cite{Hesthaven_nodal_DG} or in DGSEM \cite{DG_SBP_gassner}. As demonstrated in \cite{Cockburn_lcs4}, quadrature rules exact for polynomial of degree respectively  $2\,k$ for volume integrals and $2\,k+1$ for surface ones have to be used to reach the excepted accuracy. In this article, such approach is chosen. In the remainder, $\oint$ will refer to quadrature approximated integration, while $\int$ holds for exact integration. Obviously, for polynomials of degree up to $2\,k$ and $2\,k+1$ respectively for volume and surface integrals, $\int$ and $\oint$ are indeed equivalent.  
\end{remark}

\begin{remark}
  \label{quadrature_entropy}
  In \cite{jiang}, G.-S. Jiang and C.-W. Shu proved that DG schemes solving SCL do ensure a cell entropy inequality, for the square entropy $\eta(u)=\demi\,u^2$. Similarly in \cite{hou}, this square stability analysis has been extended by S. Hou and X.-D. Liu to symmetric systems. However, those demonstrations rely on exact calculation of integrals, and thus does not applied if one uses quadrature rules. They are also limited to the square entropy.
\end{remark}

\begin{figure}[!ht]
  \begin{center}
    \subfigure[Quad/tri: $N_s=N_k$]{\includegraphics[height=3.8cm]{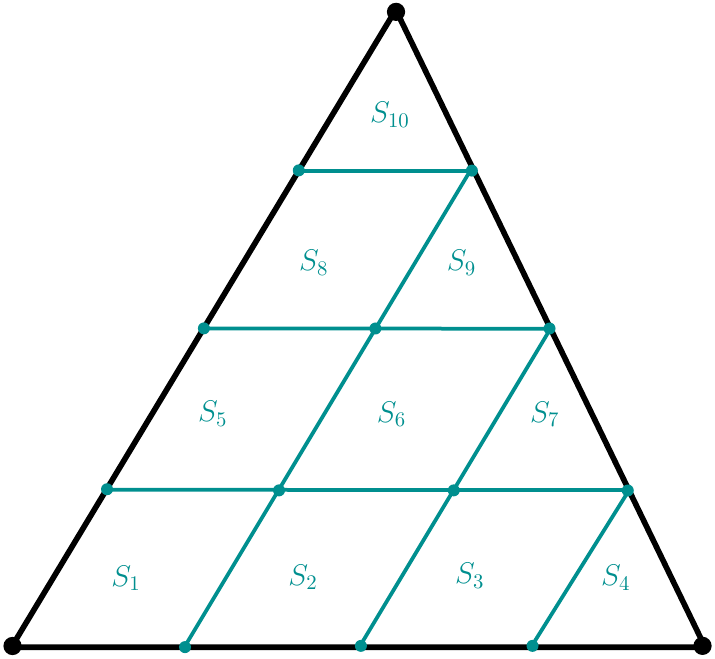}\label{fig_tri1}}\hspace*{8mm}
    \subfigure[Voronoi-type: $N_s=N_k$]{\includegraphics[height=3.8cm]{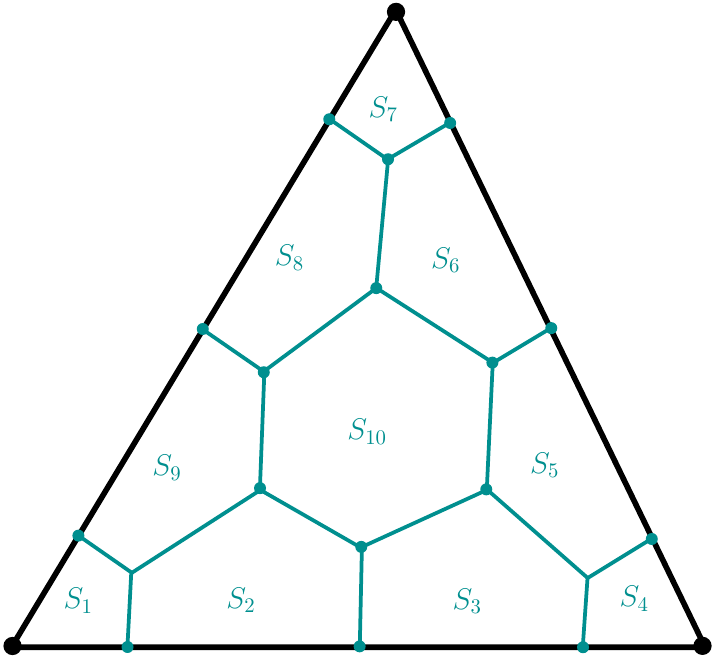}\label{fig_tri2}}\hspace*{8mm}
    \subfigure[Tri: $N_s>N_k$]{\includegraphics[height=3.8cm]{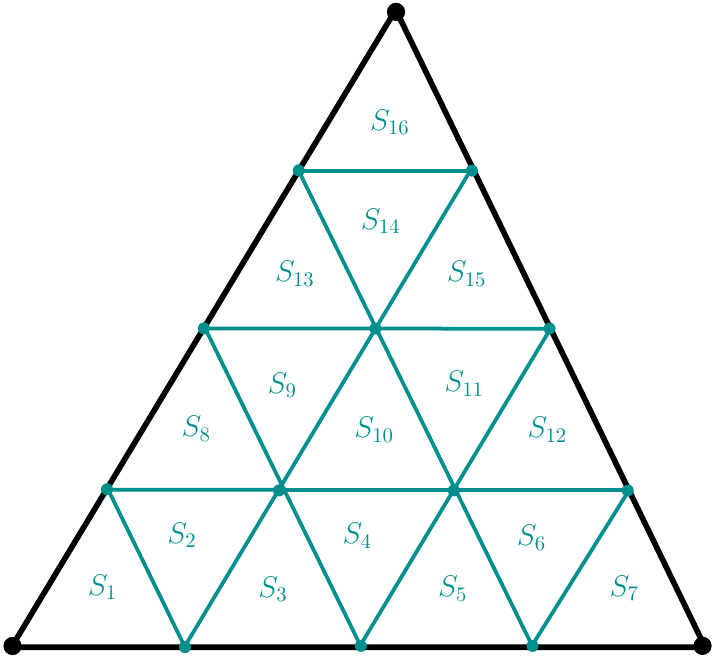}\label{fig_tri3}}
    \caption{Examples of easily generalizable subdivisions for a triangular cell and a $\P^{\,3}$ DG scheme ($N_k=10$)}
  \label{fig_subdiv_tri}
  \end{center}
\end{figure}

 In \eqref{DG_2D_2}, we identify $\int_{\omega_c} \sigma^c_m\,\sigma^c_p\,\dd V=\(M_c\)_{mp}$ as the generic coefficient of the symmetric mass matrix $M_c\in\mc{M}_{N_k}$. The scheme \eqref{DG_2D_2} can then be reformulated in a compact matrix-vector form as
\begin{align}
  \label{DG_2D_3}
  M_c\,\ddt{U_c}=\Phi_c,
\end{align}
with $(U_c)_m=u_m^c$ the solution vector filled with the polynomial moments, and where the so-called DG residuals $\Phi_c\in\R^{N_k}$ is defined as $(\Phi_c)_m=\oint_{\omega_c} \bs{F}(u_h^c)\pds\gradx{\sigma^c_m}\,\dd V-\oint_{\partial \omega_c}\hspace*{-2mm}\sigma^c_m\;\mc{F}_n\,\dd S$. Now, aiming at reformulating DG scheme \eqref{DG_2D_3} as a subgrid FV-like scheme, let us subdivide the mesh cells into subcells, similarly to what we did in \cite{vilar_aplsc_1D,vilar_aplsc_2D}. Let us emphasize that here we can relax the constraint requiring the number of subcells to match the dimension of the functional space. Hence, the subdivision can be chosen very freely. If $N_k$ stands for the number of degrees of freedom in a given cell, let define $N_s$ as the number of subcells in the cell. The only constraint we impose here is to have enough subcells not to be under-resolved, hence we impose $N_s\geq N_k$. In Figure~\ref{fig_subdiv_tri}, we present some easily generalizable subdivisions for triangle cells. Let us mention that in the first two, Figures~\ref{fig_tri1} and \ref{fig_tri2}, $N_s=N_k$. This is no more the case in the last one, Figure~\ref{fig_tri3}, as the number of subcells exceed by a lot the dimension of $\P^{\,3}$.\\

\begin{remark}
  \label{rem_3D}
  Let us emphasize that, while only triangular grids are considered for numerical applications, see Sections~\ref{sect_results_entropy} and \ref{sect_results_GLMP}, the following demonstration as well as the \scheme scheme presented in the remainder are not limited to this case. Any grid made of generic polygonal cells can be considered. Furthermore, apart from the coding aspects, the present analysis and monolithic scheme can also be directly extended to 3D geometries.
\end{remark}

  That being said, let us consider a cell $\omega_c$ and its subdivision into $N_s$ subcells $S_m^c$, for $m\in\bk{1,N_s}$. Then, we define the numerical solution subcell mean values, also referred to as submean values, as $\ov{u}_m^{\,c}=\inv{|S_m^c|}\int_{S_m^c}u_h^c\,\dd V$. To express DG scheme as a subgrid FV-like method, we want find the so-called reconstructed fluxes $\wh{F_{mp}}$ such that
\begin{align}
  \label{subcell_FV_DG}
  \ddt{\ov{u}_m^c}=-\Inv{|S_m^c|}\,\Sum_{S_{p}^v\,\in\mc{V}_m^{\,c}} l_{mp}\, \wh{F_{mp}}.
\end{align}
In equation \eqref{subcell_FV_DG}, $\mc{V}_m^{\,c}$ denotes the set of face neighboring subcells of $S_m^c$, while $l_{mp}$ stands for the length of the interface $f_{mp}$ between subcells $S_m^c$ and $S_p^v$. Let us highlight that $S_{p}^v\,\in \,\mc{V}_m^{\,c}$ can either be inside cell $\omega_c$ or in one of its neighbors $\omega_v\in\mc{V}_c$. As in \cite{vilar_aplsc_2D}, we impose on the boundary of cell \;$\omega_c$, so for $S_p^v\not\subset \omega_c$, that the reconstructed flux is nothing but the DG numerical flux, $i.e.\,$ $l_{mp}\,\wh{F_{mp}}=\oint_{f_{mp}} \mc{F}\(u_h^c,\,u_h^v,\,\bs{n}_{mp}\)\;\dd S$. As details of the proof have been given in \cite{vilar_aplsc_2D}, let us simply recall the final formula to compute the subcell interior faces reconstructed fluxes

\begin{align}
  \label{recons_flux_definition}
  \wh{F_c}=-A_c\tra\,\mc{L}_c^{-1}\(D_c\,P_c\,M_c^{-1}\,\Phi_c+B_c\).
\end{align}

\begin{figure}[!ht]
  \begin{center}
    \subfigure[Quad/tri subdivision]{\includegraphics[width=5.2cm]{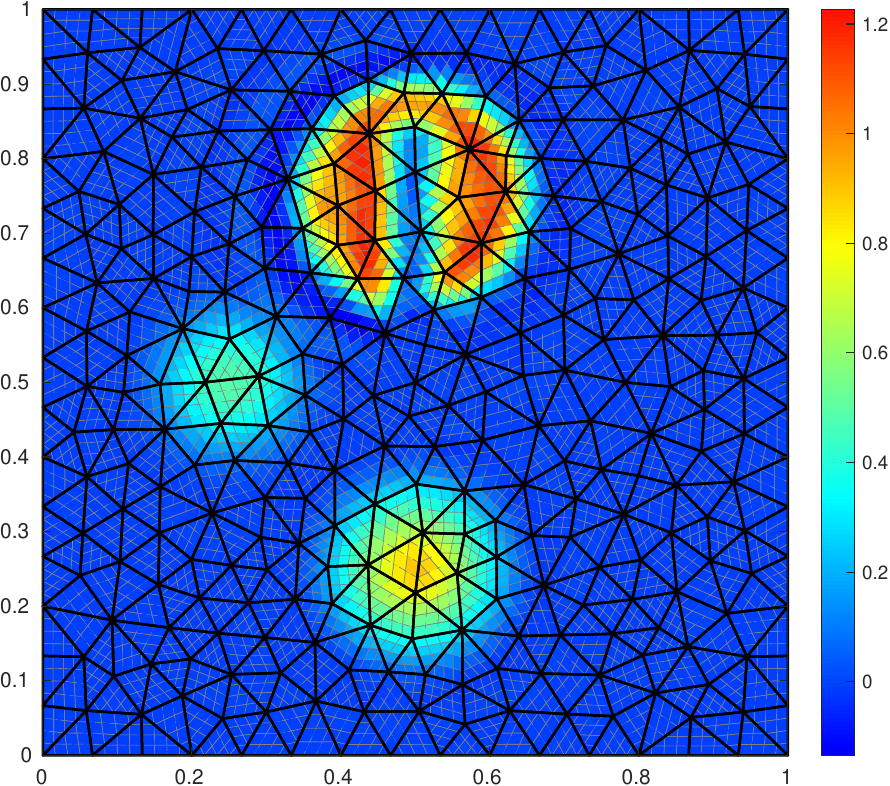}}\hspace*{3mm}
    \subfigure[Voronoi-type subdivision]{\includegraphics[width=5.2cm]{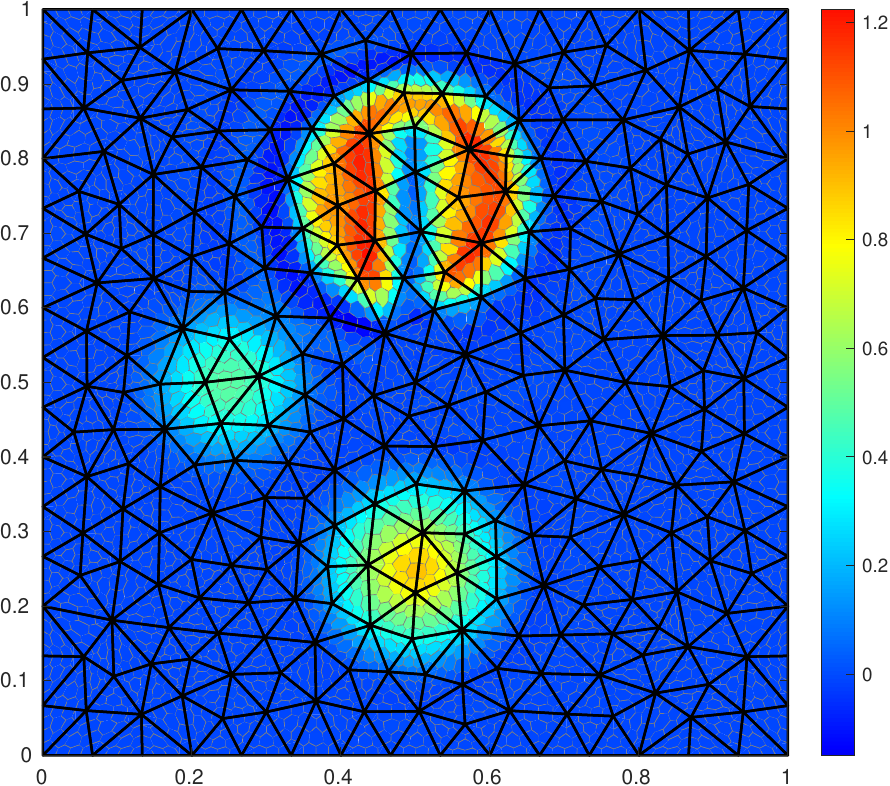}}\hspace*{3mm}
    \subfigure[Triangular subdivision]{\includegraphics[width=5.2cm]{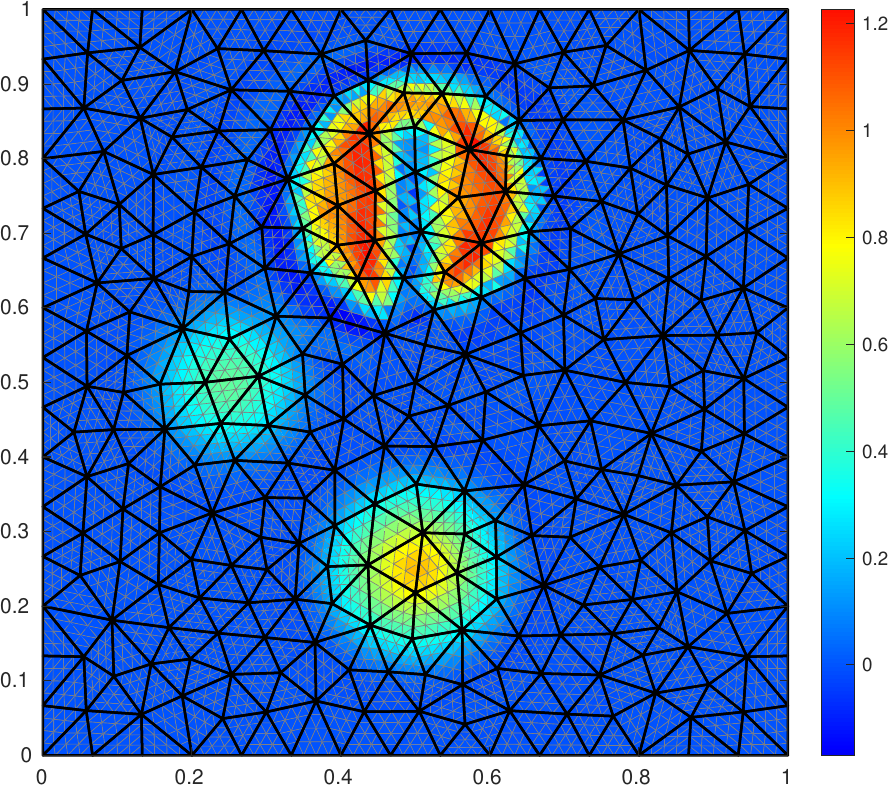}}
    \caption{$\P^{\,3}$reconstructed flux FV schemes on 576 cells: subcells mean values}
    \label{Fig_DG_reconstructed}
  \end{center}
\end{figure}

In \eqref{recons_flux_definition}, if $N_f^c$ denotes the number of subcells' faces inside $\omega_c$, meaning not belonging to $\partial \omega_c$, the vector $\wh{F_c}\in\R^{N_f^c}$ then contains all the interior faces reconstructed fluxes weighted by the face length, $i.e.\,$ $l_{mp}\,\wh{F_{mp}}$. Matrix $A_c\in\mc{M}_{N_s\times N_f^c}$ stands for the adjacency matrix, $\mc{L}_c^{-1}\in\mc{M}_{N_s}$ the generalized inverse of the graph Laplacian matrix of the subdivision, $D_c=\text{diag}(|S_1^c|,\dots,|S_{N_k}^c|)\in\mc{M}_{N_s}$ the subcells volume matrix, $P_c\in\mc{M}_{N_s\times N_k}$ the projection matrix such that $(P_c)_{mp}=\inv{|S_m^c|}\,\int_{S_m^c} \sigma^c_p\,\dd V$ and $(B_c)_m=\oint_{\partial S_m^c\cap\partial \omega_c}\mc{F}_n\;\dd S$ the cell boundary contribution. Definition of all these matrices can be found in \cite{vilar_aplsc_2D}. Let us just recall that matrices $A_c$ and $\mc{L}_c^{-1}$ only depends on the cell subdivision connectivity, while $D_c$, $P_c$ and $M_c$ depends on the chosen basis function and their values on the subcells. All those matrices can be computed initially, once and for all. The only time dependent quantities in \eqref {recons_flux_definition} are $\Phi_c$, the DG residual which is computed and available in any DG code, and term $B_c$ which is required to close the linear system to solve and ensure that \eqref{recons_flux_definition} is indeed the unique solution. Let us mention that different cell subdivisions will lead to different reconstructed flux values, but will still be equivalent to the same unique DG numerical solution, see Figure~\ref{Fig_DG_reconstructed}.

\begin{remark}
  \label{remark_least_square}
  To go from the polynomial representation of the solution to its submean values, we make use of the projection matrix $P_c$ as $\ov{U}_c=P_c\,U_c$, where $\ov{U}_c\in\R^{N_s}$ is the vector containing all the subcell mean values in cell $\omega_c$. Now, working with the piecewise constant representation of the numerical solution on the subcells through equation \eqref{subcell_FV_DG}, one still needs the polynomial representation of the solution in the computation of the DG residual. Then, to go from the submean values $\ov{u}_m^c$ to the polynomial moments $u_m^c$, we make use of the following least square procedure
  \begin{align}
    \label{least_square}
    U_c=(P_c\tra\,P_c)^{-1}\, P_c\tra\;\ov{U}_c.
  \end{align}
  In the light of \eqref{least_square}, we consider a subdivision to be admissible if matrix $P_c\tra\,P_c$ is indeed invertible, which has been the case for all subdivisions we have studied. Let us emphasize that in the case where $N_s=N_k$, relation \eqref{least_square} reduces to $U_c=P_c^{-1}\,\,\ov{U}_c$.
\end{remark}

To make sure that, regardless the type of cell subdivision, the FV-like scheme \eqref{subcell_FV_DG} provided with the reconstructed fluxes definition \eqref{recons_flux_definition} does indeed produce the DG numerical solution defined in equation \eqref{DG_2D_1}, let us run the classical solid body rotation test case taken from \cite{LeVeque2}. To this end, we then consider \eqref{lcs1} with a divergence-free velocity field corresponding to a rigid rotation, defined by $\bs{F}(u,\bs{x})=(\demi-y,\, x-\demi)\tra\,u$. We apply this solid body rotation to an initial datum which includes both a plotted disk, a cone and a smooth hump. We run this test as a FV-like scheme \eqref{subcell_FV_DG} associated with definition \eqref{recons_flux_definition}, where the DG residual has been computed as for a $\P^{\,3}$ DG scheme. The three types of subdivision displayed in Figure~\ref{fig_subdiv_tri} have been considered, see Figure~\ref{Fig_DG_reconstructed}.

\begin{figure}[!ht]
  \begin{center}
    \includegraphics[width=9cm]{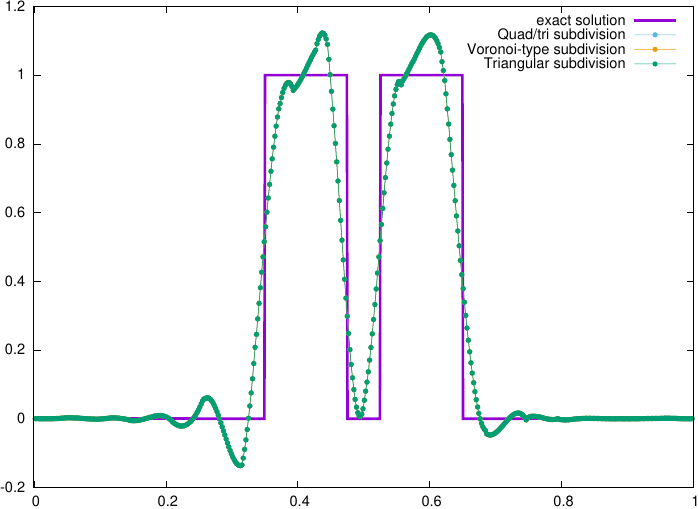}
    \caption{$\poly{3}$ reconstructed flux FV-type scheme on 576 cells: polynomial solution values on line $y=0.75$}
  \label{Fig_DG_reconstructed_profile}
  \end{center}
\end{figure}

In the light of Figures~\ref{Fig_DG_reconstructed} and \ref{Fig_DG_reconstructed_profile}, the three computations, involving three different types of cell subdivision, do produce the same numerical solution, which is nothing but the one that a $\P^{\,3}$ DG code would have produced.
Figures~\ref{Fig_DG_reconstructed} and \ref{Fig_DG_reconstructed_profile}, although demonstrating how accurate DG schemes are, also highlight the need of further limitation or correction to ensure an admissible behavior. Indeed, while the unique entropic weak solution is supposed to remained bounded by the minimum and maximum of the initial datum $u_0$ (the so-called maximum principle), the numerical solutions in Figure~\ref{Fig_DG_reconstructed_profile} clearly violate this principle. In the context of systems, as the Euler compressible gas dynamics one for instance, this maximum principle translates into a positivity principle, where some quantities have to remain positive, as the density and internal energy in the aforementioned Euler case. The non-preservation of the positivity of the numerical solution is absolutely critical, as it can lead to the crash of the simulation code. In both Figures~\ref{Fig_DG_reconstructed} and \ref{Fig_DG_reconstructed_profile}, one can furthermore clearly see the well-known Gibbs phenomenon, for which the approximation of a discontinuity through a high-order scheme will produce spurious oscillations. Those phenomena, along with the capture of non-entropic weak solutions, require some stabilization or correction techniques. This paper aims at presenting a \scheme scheme, where DG scheme will be blended at the subcell scale with a first-order FV scheme, in order to combine the best of the two worlds, namely accuracy and robustness.

\section{Local subcell monolithic DG/FV scheme}
\label{sect_monolithic}

\subsection{Blended fluxes and intermediate Riemann states}
\label{subsect_blended_fluxes}

The previous reformulations of DG scheme into subcell FV-like scheme through the definition of reconstructed fluxes enable us to construct our \scheme scheme. To this end, each face $f_{mp}$ of each subcell $S_m^c$ will be assigned two fluxes, one reconstructed flux $\wh{F_{mp}}$ giving the equivalency with a high-order DG scheme and one first-order FV numerical flux $\mc{F}_{mp}^\text{\,\tiny FV}=\mc{F}\(\ov{u}_m^c,\ov{u}_p^v,\bs{n}_{mp}\)$, where $\bs{n}_{mp}$ is outward unit normal of face $f_{mp}$. Then, these two fluxes will be blended in a convex manner through a blending coefficient $\theta_{mp}\in[0,1]$ as in following
\begin{align}
  \label{blended_flux}
  \wt{F_{mp}}=\mc{F}_{mp}^\text{\,\tiny FV}+\theta_{mp}\,\underbrace{\(\wh{F_{mp}}-\mc{F}_{mp}^\text{\,\tiny FV}\)}_{\Delta F_{mp}^{\phantom{DG}}}.
\end{align}
A blending coefficient set to 0 will lead to a first-order FV numerical flux, while a coefficient set at 1 will induce a high-order DG reconstructed flux. The \scheme then writes as follows
\begin{align}
  \label{subcell_monolithic}
  \ddt{\ov{u}_m^c}=-\Inv{|S_m^c|}\,\Sum_{S_{p}^v\,\in\mc{V}_m^{\,c}} l_{mp}\, \wt{F_{mp}}.
\end{align}
The goal is now to determine, through analysis, the optimal coefficients to reach the desired properties while trying to maintain as much as possible the high accuracy of the scheme. To do so, we rewrite the monolithic scheme \eqref{subcell_monolithic} as a Godunov-like scheme. But first, as only the semi-discrete version of the analysis and the monolithic scheme have been presented, we make use of SSP Runge-Kutta (RK) time integration method \cite{Osher} to achieve high-accuracy in time. In the light of the fact that these multistage time integration methods write as convex combinations of first-order forward Euler scheme, the monolithic DG/FV scheme will be presented for the simple case of this latter time numerical scheme, for sake of simplicity. The semi-discrete scheme \eqref{subcell_monolithic} provided with first-order forward Euler time integration writes
\begin{align}
  \label{discrete_subcell_monolithic}
  \ov{u}_m^{c,n+1}=\ov{u}_m^{c,n}-\Frac{\Dt}{|S_m^c|}\,\Sum_{S_{p}^v\,\in\mc{V}_m^{\,c}} l_{mp}\, \wt{F_{mp}},
\end{align}
where all the quantities involved in the definition of the blended flux $ \wt{F_{mp}}$ are taken at time level $n$ (at the previous Runge-Kutta stage in a RK time integration). Defining $\gamma_{mp}=\gamma\(\ov{u}_m^{c,n},\ov{u}_p^{v,n},\bs{n}_{mp}\)$ and recalling that $\sum_{S_{p}^v\,\in\mc{V}_m^{\,c}}l_{mp}\,\bs{n}_{mp}=\bs{0}$, let us now rewrite $\ov{u}_m^{c,n+1}$ as a convex combination of quantities defined at the previous time step
\begin{align*}
  \ov{u}_m^{c,n+1}&=\ov{u}_m^{c,n}-\Frac{\Dt}{|S_m^c|}\,\Sum_{S_{p}^v\,\in\mc{V}_m^{\,c}} l_{mp}\, \wt{F_{mp}} \; \Dgray{\pm\,\Frac{\Dt}{|S_m^c|}\Sum_{S_{p}^v\,\in\mc{V}_m^{\,c}}l_{mp}\,\gamma_{mp}\,\ov{u}_m^{c,n}+\Frac{\Dt}{|S_m^c|}\bs{F}(\ov{u}_m^{c,n})\pds\Sum_{S_{p}^v\,\in\mc{V}_m^{\,c}}l_{mp}\,\bs{n}_{mp}},\\[3mm]
  &= \Big(1-\Frac{\Dt}{|S_m^c|}\,\Sum_{S_{p}^v\,\in\mc{V}_m^{\,c}}l_{mp}\,\gamma_{mp}\Big)\,\ov{u}_m^{c,n}+\Frac{\Dt}{|S_m^c|}\,\Sum_{S_{p}^v\,\in\mc{V}_m^{\,c}}l_{mp}\,\gamma_{mp}\,\Big(\ov{u}_m^{c,n}-\Frac{\wt{F_{mp}}-\bs{F}(\ov{u}_m^{c,n})\pds\bs{n}_{mp}}{\gamma_{mp}}\Big).
\end{align*}
Then, defining the left blended Riemann intermediate state $\wt{u_{mp}}^-=\ov{u}_m^{c,n}-\frac{\wt{F_{mp}}-\bs{F}(\ov{u}_m^{c,n})\,\pds\,\bs{n}_{mp}}{\gamma_{mp}}$, the previous expression can finally be recast into the following convex form
\begin{align}
  \label{convex_combo}
  \ov{u}_m^{c,n+1}= \Big(1-\Frac{\Dt}{|S_m^c|}\,\Sum_{S_{p}^v\,\in\mc{V}_m^{\,c}}l_{mp}\,\gamma_{mp}\Big)\,\ov{u}_m^{c,n}+\Frac{\Dt}{|S_m^c|}\,\Sum_{S_{p}^v\,\in\mc{V}_m^{\,c}}l_{mp}\,\gamma_{mp}\,\wt{u_{mp}}^-.
\end{align}
Consequently, $\ov{u}_m^{c,n+1}$ does indeed write as a convex combination of previous time step quantities under the standard CFL condition in this subcell context
\begin{align}
  \label{CFL}
  \Dt\leq\Frac{|S_m^c|}{\Sum_{S_{p}^v\,\in\mc{V}_m^{\,c}}l_{mp}\,\gamma_{mp}}.
\end{align}
By mean of the numerical flux definition \eqref{num_flux}, the left blended Riemann intermediate state $\wt{u_{mp}}^-$ can be rewritten into the following form
\begin{align}
  \label{blended_intermediate}
  \wt{u_{mp}}^-=\Frac{\ov{u}_m^{c,n}+\ov{u}_p^{v,n}}{2}-\Frac{\big(\bs{F}(\ov{u}_p^{v,n})-\bs{F}(\ov{u}_m^{c,n})\big)\pds\bs{n}_{mp}}{2\,\gamma_{mp}}-\theta_{mp}\Frac{\Delta F_{mp}}{\gamma_{mp}}=u_{mp}^{\ast,\text{\,\tiny FV}}-\theta_{mp}\Frac{\Delta F_{mp}}{\gamma_{mp}},
\end{align}
where $u_{mp}^{\ast,\text{\,\tiny FV}}$ is nothing but the first-order FV Riemann intermediate state. It is straightforward to prove that, since $\gamma_{mp}\geq\max_{w\,\in\, I(\ov{u}_m^{c,n},\ov{u}_p^{v,n})}\big(|\bs{F}'(w)\pds\bs{n}_{mp}|\big)$, we have then $u_{mp}^{\ast,\text{\,\tiny FV}}\in I(\ov{u}_m^{c,n},\ov{u}_p^{v,n})$, see Appendix~\ref{sect_FV_DMP}.\\

For sake of clarity, let us specify conservativity relation. Obviously, we have that $\bs{n}_{pm}=-\bs{n}_{mp}$ as well as $\mc{F}_{pm}^\text{\,\tiny FV}=-\mc{F}_{mp}^\text{\,\tiny FV}$, $\wh{F_{pm}}=-\wh{F_{mp}}$ and $\wt{F_{pm}}=-\wt{F_{mp}}$, while $\theta_{pm}=\theta_{mp}$ and $u_{pm}^{\ast,\text{\,\tiny FV}}=u_{mp}^{\ast,\text{\,\tiny FV}}$. In the light of these relations, it is clear the right blended Riemann intermediate state $\wt{u_{mp}}^+:=\wt{u_{pm}}^-$ hence writes $\wt{u_{mp}}^+=u_{mp}^{\ast,\text{\,\tiny FV}}+\theta_{mp}\frac{\Delta F_{mp}}{\gamma_{mp}}$. It then appears that, where in first-order FV scheme we have only one Riemann intermediate state, here we have two, $\wt{u_{mp}}^\pm=u_{mp}^{\ast,\text{\,\tiny FV}}\pm\theta_{mp}\frac{\Delta F_{mp}}{\gamma_{mp}}$, which both rely on the admissible first-order one $u_{mp}^{\ast,\text{\,\tiny FV}}$. Consequently, introducing $G$, a convex admissible set where the solution has to remain in, if the numerical initial solution does lie in $G$, then it is always possible to find blending coefficients $\theta_{mp}$ to ensure that $\ov{u}_m^{c,n}$ remains in $G$ during the whole calculation.

\begin{remark}
  \label{initialization}
  To ensure the approximated solution to be in $G$ at the initial time, the initialization has to be carried out by computing the subcell mean values and then use the projection matrix $P_c$ to recover the cell polynomial representation, and not by a $L_2$ projection or an interpolation as it is generally done in DG schemes.
\end{remark}

As long as the first-order FV scheme used should achieve the desired properties, one can find blending coefficients for the high-order local subcell monolithic scheme to do as well. The different conditions on the blending coefficients will be formulated as inequalities. Thus, to combine several properties, one just have to take the minimum of the corresponding conditions.\\

Let us enlighten that it has been previously observed, \cite{vilar_aplsc_1D,vilar_aplsc_2D}, that a stiff transition from first-order to a fully high-order scheme would produce more oscillatory solutions. It will be the case in 2D if a subcell will be assigned with first-order FV fluxes on some faces (corresponding to $\theta_{mp}=0$) as well as fully high-order reconstructed fluxes (corresponding to $\theta_{mp}=1$) on some other faces. In order to avoid such strong variation in fluxes accuracy, we will make use blending coefficients smoothers.

\subsection{Blending coefficients smoothening}
\label{smoother}

In our previous work \cite{vilar_aplsc_2D}, \textit{a posteriori} blending of high-order reconstructed fluxes with first-order FV fluxes have been performed with arbitrary blending coefficients. And, in the context of non-linear problems, to avoid to yield too strong transition from high to low orders, a wider blending stencil with increasing coefficients (increasing according to the sequence $\{0, \inv{4}, \demi, \frac{3}{4}, 1\}$) was used. Here, in this \textit{a priori} monolithic framework, we make use of a simple procedure to avoid very stiff order transition. Let us mention that we only make use of this coefficient smoother for solving non-linear problems, as it is no needed in the linear case. Furthermore, two different versions will be used, as a more constraining one will help in the context of non-convex fluxes SCL.

\subsubsection{Coefficients smoothening n°1}
\label{smoother1}

Each subcell $S_m^c$ will be given a blending coefficient $\theta_m^{\,c}$ defined as the average of its faces blending coefficients, as in follows
\begin{align}
  \label{subcell_blending_coeff}
  \theta_m^{\,c}=\Inv{\#\mc{V}_m^{\,c}}\,\Sum_{S_{p}^v\,\in\,\mc{V}_m^{\,c}} \theta_{mp}.
\end{align}

Then, each subcell's face blending coefficient might be potentially reduced to the average of subcells' blending coefficients of every subcells sharing a node with face $f_{mp}$, as
\begin{align*}
  \wt{\theta}_{mp}=\min\(\theta_{mp},\,\Inv{\#\mc{V}_{mp}}\,\Sum_{S_q^v\,\in\,\mc{V}_{mp}} \theta_q^{\,v}\).
\end{align*}
Here, $\mc{V}_{mp}$ is the set containing all the subcells that share at least one node with face $f_{mp}$.

\subsubsection{Coefficients smoothening n°2}
\label{smoother2}

This second smoother is nothing but the first one where the different averaged values are substituting by minimum values. Consequently, each subcell $S_m^c$ will be given a blending coefficient $\theta_m^{\,c}$ defined as the minimum of its faces blending coefficients, as in follows
\begin{align*}
  \theta_m^{\,c}=\Min_{S_{p}^v\,\in\,\mc{V}_m^{\,c}} \theta_{mp}.
\end{align*}

Then, each subcell's face blending coefficient might be potentially reduced by taking
\begin{align*}
  \wt{\theta}_{mp}=\min\(\theta_{mp},\,\Min_{S_q^v\,\in\,\mc{V}_{mp}}\theta_q^{\,v}\).
\end{align*}

Those two smoothening techniques have proved to improve the numerical results while preserving all the different properties presented in the next sections. While in Section~\ref{sect_GLMP} we focus on imposing positivity and local maximum principles to control spurious oscillations, which will prove to produce the best results, in the next section we first address the different questions regarding entropy stability.

\section{Entropy stabilities}
\label{sect_entropy}

This section is devoted to entropy stability. By means of this \scheme framework, we will attempt to address the following questions: Is it possible to find $\theta_{mp}$ the blending coefficients ensuring entropy stability? What do we mean by entropy stability? What is the cost of such constraints? Is this absolutely needed while aiming for high-order accuracy? To this end, we first introduce the definition of blending coefficients ensuring different type of entropy stabilities, while discussing the cost of such properties. Numerical results are then presented to confirm the developed theory, and to help us answer the stated queries.\\

For sake of simplicity, entropy stability will be addressed here in the simple case of SCL. Nonetheless, the extension to systems is perfectly straightforward. In the remainder, let then $\eta(u)$ be a strictly convex entropy, while $\bs{\phi}(u)$ be the associated entropy flux. $v(u)=\eta'(u)$ refers to the entropy variable, $\bs{\psi}(u)=v(u)\,\bs{F}(u)-\bs{\phi}(u)$ and $\bs{\Psi}(v)=\bs{\psi}\big(u(v)\big)$ to the entropy potential flux. Thanks to the entropy convexity, the mapping between $u$ and $v$ is indeed a diffeomorphism.

\subsection{Discrete subcell entropy stability for any entropy $\Longrightarrow$ First-order}
\label{subsect_entropy_1}

First, let us enforce, at the discrete level and subcell scale, an entropy inequality for any given entropy. To this end, we seek a blending coefficient $\theta_{mp}$ to make the blended flux $\wt{F_{mp}}$ an E-flux. Such flux does yield the desired property, as recalled in Appendix~\ref{sect_FV_entropy}. Let introduce $\wt{\gamma_{mp}}$ such that
\begin{align*}
  \wt{F_{mp}}=\Frac{\bs{F}\big(\ov{u}_m^{c,n}\big)+\bs{F}\big(\ov{u}_p^{v,n}\big)}{2}\pds\bs{n}_{mp}-\ub{\big(\gamma_{mp}-2\,\theta_{mp}\,\Delta F_{mp}\big)}{\wt{\gamma_{mp}}}\,\Frac{\ov{u}_p^{v,n}-\ov{u}_m^{c,n}}{2}.
\end{align*}
Then, to ensure a discrete entropy inequality, at the subcell level, for any entropy, it is sufficient to take $\theta_{mp}$ such that $\wt{\gamma_{mp}}\geq\max_{w\,\in\, I(\ov{u}_m^{c,n},\ov{u}_p^{v,n})}\big(|\bs{F}'(w)\pds\bs{n}_{mp}|\big)$. This conditions leads to the following condition: if $\Delta F_{mp}\,.\,\(\ov{u}_p^{v,n}-\ov{u}_m^{c,n}\)>0$ then
\begin{align}
  \label{theta_entropy_1}
  \theta_{mp}\leq \min\(1,\,\Frac{\big(\gamma_{mp}-\gamma_{\text{max}}\big)\,\big(\ov{u}_p^{v,n}-\ov{u}_m^{c,n}\big)}{2\,\Delta F_{mp}}\),
\end{align}
where $\gamma_{\text{max}}:=\Max_{w\,\in\, \mrm{I}\,(\ov{u}_m^{c,n},\ov{u}_p^{v,n})}\Big(|\bs{F}'(w)\pds\bs{n}_{mp}|\Big)$.

\begin{remark}
  \label{rem_entropy_1}
  While this particular choice does ensure the desired entropy property, it does also, as expected, lead to a \textbf{first-order} accurate scheme, as displayed in Table~\ref{table_order_adv_1D_entropy} and Figure~\ref{fig_advect_entropy_compar}. As stated in \cite{Tadmor_entropy_acta}, an E-scheme is indeed first-order accurate, and in light of condition \eqref{theta_entropy_1}, one can see for instance that $\theta_{mp}$ will be partially set to zero in the simple case of linear advection or if one uses global Lax-Friedrichs numerical flux.
\end{remark}

\subsection{Semi-discrete subcell entropy stability for one given entropy $\Longrightarrow$ Second-order}
\label{subsect_entropy_2}

Because condition \eqref{theta_entropy_1} would lead to a first-order scheme, we may relaxed our expectations for sake of accuracy. Instead of a discrete entropy stability, for any entropy, one may aim for a semi-discrete entropy inequality, for a given entropy - entropy flux pair $(\eta,\,\bs{\phi})$. To this end, we make use of the following Tadmor two-point entropy conservation/dissipation condition, see \cite{Tadmor_entropy,Tadmor_entropy_acta}
\begin{align*}
  \wt{F_{mp}}\,\Big(v\big(\ov{u}_p^{v,n}\big)-v\big(\ov{u}_m^{c,n}\big)\Big) \leq \Big(\bs{\psi}\big(\ov{u}_p^{v,n}\big)-\bs{\psi}\big(\ov{u}_m^{c,n}\big)\Big)\pds\bs{n}_{mp}.
\end{align*}
As the first-order FV flux does ensure such inequality, see Appendix~\ref{sect_FV_Tadmor}, a sufficient condition for the blended flux to do as well is: if $\Delta F_{mp}\,.\,\Big(v\big(\ov{u}_p^{v,n}\big)-v\big(\ov{u}_m^{c,n}\big)\Big)>0$ then
\begin{align}
  \label{theta_entropy_2}
  \theta_{mp}\leq \min\(1,\,\Frac{\(\Frac{\bs{\psi}\big(\ov{u}_p^{v,n}\big)-\bs{\psi}\big(\ov{u}_m^{c,n}\big)}{v\big(\ov{u}_p^{v,n}\big)-v\big(\ov{u}_m^{c,n}\big)}\)\pds\bs{n}_{mp}-\mc{F}_{mp}^{\,\text{FV}}}{\Delta F_{mp}^{\phantom{\int^{B^X}}}}\)
\end{align}

\begin{remark}
  \label{rem_entropy_2}
  As one would expect, \cite{Tadmor_entropy_acta}, condition \eqref{theta_entropy_2} does decrease the accuracy to \textbf{second-order}, see Figure~\ref{fig_advect_entropy_compar}. Practically, it has been observed that this condition even further reduces the accuracy, as illustrated in Table~\ref{table_order_adv_1D_entropy}.
\end{remark}

\subsection{Semi-discrete cell entropy stability for one given entropy $\Longrightarrow$ High-order}
\label{subsect_entropy_3}

To get a second-order accurate scheme, we had no choice but to relax a bit the aimed entropy stability, by getting a semi-discrete entropy stable scheme, at the subcell level, for one given pair $\(\eta,\,\bs{\phi}\)$. In order to preserve the order of accuracy of this monolithic DG/FV scheme, we go from a subcell semi-discrete entropy stability to a cell based one, again for a given entropy. To do so, let us first introduce $\;\{\varphi_m^c\}_m$, a particular set of $\P^k$ basis functions. Let us emphasize that in the case where $N_s>N_k$, those functions form a spanning sets. Those functions, previously introduced in \cite{vilar_aplsc_1D,vilar_aplsc_2D} and that we refer to as sub-resolution basis functions. Those functions, can be seen as the $L_2$ projection of the subcell indicator functions $\1_{S_m^c}(\bs{x})$ onto $\P^k(\omega_c)$. They are defined such that $\forall\,\psi\in\P^k(\omega_c)$ and $\forall\,m=1,\ldots,N_s$
\begin{align}
  \label{subresolution}
  \Int_{\omega_c}\varphi_m\,\psi\,\dd V=\Int_{S_m^c}\psi\,\dd V.
\end{align}
We note $\un{v}_m^c$ the corresponding moments such that $v_h^c=\sum_{m=1}^{N_s} \un{v}_m^c\,\varphi_m^c$.
%
%
Now, let us express the time variation of a given entropy $\eta$ over cell $\omega_c$
\begin{align*}
\Delta \eta_c:=\ddt{}\Oint_{\omega_c}\eta\big(u_h^c\big)\;\dd V=\Oint_{\omega_c} v\big(u_h^c\big)\;\vdt{}u_h^c\;\dd V.
\end{align*}
Referring by $v_h^c$ the $L_2$ projection of the entropy variable $v\big(u_h^c\big)$ onto $\P^k(\omega_c)$, and by means of \eqref{subresolution}, the entropy variation can be put into the following simple expression 
\begin{align*}
\Delta \eta_c=\Int_{\omega_c} v_h^c\;\vdt{}u_h^c\;\dd V=\Sum_{m=1}^{N_s} \un{v}_m^c\,\Int_{\omega_c} \varphi_m^c\;\vdt{}u_h^c\;\dd V=\Sum_{m=1}^{N_s} \un{v}_m^c\,\Int_{S_m^c} \vdt{}u_h^c\;\dd V.
\end{align*}
In the light of \eqref{least_square} and remark~\ref{remark_least_square}, the entropy time variation leads to the following compact form
\begin{align}
  \label{def_entropy_diag}
\Delta \eta_c=\Sum_{m=1}^{N_s} |S_m^c|\,\un{v}_m^c\,\ddt{\ov{u}_m^c}.
\end{align}
In the case where $N_s>N_k$, this last equality has to be understood in a least square sens. The use of the semi-discrete scheme \eqref{subcell_monolithic} provides us with the following expression
\begin{align}
  \label{def_entropy_3}
\Delta \eta_c=-\Sum_{m=1}^{N_s} \un{v}_m^c\Sum_{S_{p}^v\,\in\mc{V}_m^{\,c}} l_{mp}\, \wt{F_{mp}}.
\end{align}
Now, let $\mf{f}_c$ be the set containing the subcell's faces of any subcell in $\omega_c$, while $\breve{\mf{f}}_c$ would be the set of subcell's faces inside $\omega_c$, meaning not belonging to $\partial \omega_c$. By means of previous definition, $\#\,\breve{\mf{f}}_c=N_f^c$, and if $f_{mp}\,\in\,\mf{f}_c\setminus\,\breve{\mf{f}}_c:=\ibreve{\mf{f}}_c$ that means $f_{mp}\subset\partial \omega_c$. Manipulating the two sums in \eqref{def_entropy_3} and recalling that the reconstructed fluxes are set to the DG numerical flux on the cell boundary, we are able to recast the cell entropy time evolution as
%
\begin{align*}
  \Delta \eta_c\;\;&=\Sum_{f_{mp}\,\in\,\breve{\mf{f}}_c} l_{mp}\,\(\un{v}_p^c-\un{v}_m^c\)\, \wt{F_{mp}}\;-\lh\Sum_{f_{mp}\,\in\,\ibreve{\mf{f}}_c} \(1-\theta_{mp}\)\,l_{mp}\,\un{v}_m^c\, \mc{F}_{mp}^\text{\,\tiny FV}\\
  &\hspace{4.9cm}-\lh\Sum_{f_{mp}\,\in\,\ibreve{\mf{f}}_c} \theta_{mp}\;\un{v}_m^c\, \oint_{f_{mp}}\lh \mc{F}\(u_h^c,\,u_h^v,\,\bs{n}_{mp}\)\,\dd S.
\end{align*}
Adding and retrieving \;$\sum_{f_{mp}\,\in\,\ibreve{\mf{f}}_c} \theta_{mp}\; \oint_{f_{mp}}\lh v\big(u_h^c\big)\,\mc{F}\(u_h^c,\,u_h^v,\,\bs{n}_{mp}\)\,\dd S$\; to the previous relation, the entropy variation can be separated into two terms, $i.e.$\; $\Delta \eta_c=\text{A}+\text{B}$, where
\begin{align}
  \label{entropy_A}
  \text{A}=\Sum_{f_{mp}\,\in\,\breve{\mf{f}}_c} l_{mp}\,\(\un{v}_p^c-\un{v}_m^c\)\, \wt{F_{mp}}+\lh\Sum_{f_{mp}\,\in\,\ibreve{\mf{f}}_c} \theta_{mp}\; \oint_{f_{mp}}\lh \Big(v\big(u_h^c\big)-\un{v}_m^c\Big)\, \mc{F}\(u_h^c,\,u_h^v,\,\bs{n}_{mp}\)\,\dd S
\end{align}
and
\begin{align}
  \label{entropy_B}
  \text{B}=-\lh\Sum_{f_{mp}\,\in\,\ibreve{\mf{f}}_c} l_{mp}\,\(\(1-\theta_{mp}\)\,\un{v}_m^c\, \mc{F}_{mp}^\text{\,\tiny FV}+\frac{\theta_{mp}}{l_{mp}}\oint_{f_{mp}}\lh v\big(u_h^c\big)\,\mc{F}\(u_h^c,\,u_h^v,\,\bs{n}_{mp}\)\,\dd S\).
\end{align}
A sufficient to ensure a correct cell entropy inequality would then be to yield that
\begin{align}
  \label{entropy_3_CS1}
  \text{A}\leq\lh\Sum_{f_{mp}\,\in\,\ibreve{\mf{f}}_c} l_{mp}\,\(\(1-\theta_{mp}\)\,\bs{\Psi}\big(\un{v}_m^c\big)+\frac{\theta_{mp}}{l_{mp}}\oint_{f_{mp}}\lh \bs{\psi}\big(u_h^c\big)\,\dd S\)\pds\bs{n}_{mp}.
\end{align}
Indeed, this sufficient condition would ensure the following inequality
\begin{align}
  \label{entropy_3_variation_1}
  \Delta \eta_c\leq-\lh\lh\Sum_{f_{mp}\,\in\,\ibreve{\mf{f}}_c} \lh l_{mp}\(\(1-\theta_{mp}\)\,\Big(\un{v}_m^c\, \mc{F}_{mp}^\text{\,\tiny FV}-\bs{\Psi}\big(\un{v}_m^c\big)\pds\bs{n}_{mp}\Big)+\frac{\theta_{mp}}{l_{mp}}\oint_{f_{mp}}\lh \Big(v\big(u_h^c\big)\,\mc{F}_n-\bs{\psi}\big(u_h^c\big)\pds\bs{n}_{mp}\Big)\,\dd S\).
\end{align}
The right hand side contains two contributions, a low-order one and a high-order one. Regarding the latter, if the numerical flux used in DG for the calculation of the reconstructed fluxes ensures the two-point Tadmor condition
\begin{align}
  \label{entropy_FV_Tadmor}
  \mc{F}\(u_L,u_R,\bs{n}\)\,\big(v(u_R)-v(u_L)\big)\leq\big(\bs{\psi}(u_R)-\bs{\psi}(u_L)\big)\pds\bs{n},
\end{align}
then the high-order part induces
\begin{align*}
  -\inv{l_{mp}}\oint_{f_{mp}}\lh \Big(v\big(u_h^c\big)\,\mc{F}_n-\bs{\psi}\big(u_h^c\big)\pds\bs{n}_{mp}\Big)\,\dd S\leq-\inv{l_{mp}}\oint_{f_{mp}}\lh \phi^\ast\big(u_h^c,\,u_h^v,\,\bs{n}_{mp}\big)\,\dd S:=-\wh{\phi_{mp}},
\end{align*}
where the consistent numerical entropy flux $\phi^\ast$ is defined as
\begin{align}
  \label{num_entr_flux_Tadmor}
  \phi^\ast(u_L,u_R,\bs{n})=\Frac{\big(v(u_L)+v(u_R)\big)}{2}\,\mc{F}(u_L,u_R,\bs{n})-\Frac{\big(\bs{\psi}(u_L)+\bs{\psi}(u_R)\big)}{2}\pds\bs{n}.
\end{align}
Similarly, by means of an \textit{appropriate} FV flux $\mc{F}_{mp}^\text{\,\tiny FV}$, the low-order contribution in \eqref{entropy_3_variation_1} produces 
\begin{align}
  \label{entropy_3_rel_FV}
  -\Big(\un{v}_m^c\, \mc{F}_{mp}^\text{\,\tiny FV}-\bs{\Psi}\big(\un{v}_m^c\big)\pds\bs{n}_{mp}\Big)\leq-\phi^\ast\Big(u\big(\un{v}_m^c\big),\,u\big(\un{v}_p^v\big),\,\bs{n}_{mp}\Big):=-\phi_{mp}^\text{\,\tiny FV}.
\end{align}
Combining the low and high contributions, the sufficient condition \eqref{entropy_3_CS1} ensures that
\begin{align}
  \label{entropy_3_variation_2}
  \ddt{}\Oint_{\omega_c}\eta\big(u_h^c\big)\;\dd V\leq-\lh\Sum_{f_{mp}\,\in\,\ibreve{\mf{f}}_c} \lh l_{mp}\,\(\(1-\theta_{mp}\)\,\phi_{mp}^\text{\,\tiny FV}+\theta_{mp}\,\wh{\phi_{mp}}\):=-\lh\Sum_{f_{mp}\,\in\,\ibreve{\mf{f}}_c} \lh l_{mp}\,\wt{\phi_{mp}},
\end{align}
which guarantees the entropy stability over cell $\omega_c$, for a given entropy $\eta$.

\begin{remark}
  \label{rem_modif_FV}
  We previously said that by means of \textit{appropriate} FV fluxes $\mc{F}_{mp}^\text{\,\tiny FV}$, relation \eqref{entropy_3_rel_FV} is ensured. To this end, the first-order FV fluxes, previously defined as $\mc{F}_{mp}^\text{\,\tiny FV}=\mc{F}\(\ov{u}_m^c,\ov{u}_p^v,\bs{n}_{mp}\)$, have to be modified for this entropy consideration as follows
  \begin{align}
  \label{modif_FV}
    \mc{F}_{mp}^\text{\,\tiny FV}=\mc{F}\Big(u\big(\un{v}_m^c\big),\,u\big(\un{v}_p^v\big),\,\bs{n}_{mp}\Big).
  \end{align}
  While this definition, along with condition \eqref{entropy_3_CS1}, guarantees entropy stability, global maximum and positivity preserving principles, introduced in the next Section~\ref{sect_GLMP}, may not be assured anymore.
\end{remark}

Lastly, to ensure the sufficient condition \eqref{entropy_3_rel_FV}, let us show that is possible to reformulate it as a continuous Knapsack problem, similarly to \cite{Chan_knapsack}. In this aforementioned paper, Y. Lin and J. Chan have introduced a way to ensure a semi-discrete cell entropy inequality for monolithic SEMDG schemes, by ultimately solving a continuous Knapsack optimization problem. This technique has been applied to systems of conservation laws in one-dimensional or tensor-product multi-dimensional settings. Investigations are currently carried out in the extension to simplex elements, but based on connecting quadrature points to build the SBP low-operator and not on any subcell representation. The main difference between the monolithic framework used in \cite{Kuzmin_Gassner_monolithic_2024,Chan_knapsack} and the one presented here resides in the fact that the former are for now mostly limited to tensor-product Cartesian grids, as the summation-by-parts property of the scheme derives from specific quadrature points based solution approximation and flux collocation. The framework presented here is by construction multi-dimensional and can be theoretically applied to any type of grids with great flexibility in the choice of cell sub-partition. Let us emphasize nonetheless that in \cite{Kuzmin_Gassner_monolithic_2024,Chan_knapsack}, thanks to the Gauss-Lobatto representation, the solution point-values can be simultaneously considered as subcell mean values. This characteristic, while mainly limiting the scheme to one-dimensional or tensor-product multi-dimensional geometries, makes the entropy analysis way simpler, as only one set of data is involved. Here, two sets of data are required in the entropy development, see definition \eqref{def_entropy_diag}, namely the solution subcell mean values $\ov{u}_m^c$ and the entropy solution sub-resolution moments $\un{v}_m^c$. This is the reason why, as said in Remark~\ref{rem_modif_FV}, the first-order FV numerical fluxes require to be modified if one aims at this cell entropy stability.\\

Let us finally formulate the sufficient condition \eqref{entropy_3_rel_FV} as the following continuous Knapsack problem
\begin{align}
  \label{Knapsack_pb}
  \text{\textbf{C}}_c\pds\bs{\Theta}_c\leq \text{D}_c,
\end{align}
where the vector $\bs{\Theta}_c=\(\theta_1,\,\ldots,\,\theta_{\#\ibreve{\mf{f}}_c},\,\theta_{\#\ibreve{\mf{f}}_c+1},\,\ldots,\,\theta_{\#\mf{f}_c}\)\tra$ contains all the boundary subcells' faces and interior subcells' faces, where $\text{D}_c$ the right hand side is defined as
\begin{align}
  \label{Knapsack_D}
  \text{D}_c\;\;=\Sum_{f_{mp}\,\in\,\ibreve{\mf{f}}_c} l_{mp}\,\bs{\Psi}\big(\un{v}_m^c\big)\pds\bs{n}_{mp}\;\;-\Sum_{f_{mp}\,\in\,\breve{\mf{f}}_c} l_{mp}\,\(\un{v}_p^c-\un{v}_m^c\)\, \mc{F}_{mp}^\text{\,\tiny FV},
\end{align}
while vector $\text{\textbf{C}}_c$ writes as, with $f_i=f_{mp}$

\begin{align*}
  \text{C}_i^c\;\;=\left\{\begin{array}{ll}
  \Oint_{f_{mp}}\lh\Big( \big(v(u_h^c)-\un{v}_m^c\big)\,\mc{F}_n-\big(\bs{\psi}(u_h^c)-\bs{\Psi}(\un{v}_m^c)\big)\pds\bs{n}_{mp}\Big)\,\dd S, \qquad &\forall i=1,\ldots,\#\ibreve{\mf{f}}_c,\\[7mm]
  l_{mp}\,\(\un{v}_p^c-\un{v}_m^c\)\, \Delta F_{mp},  &\forall i=\#\ibreve{\mf{f}}_c+1,\ldots,\#\mf{f}_c.
  \end{array}\right.
\end{align*}\vspace*{2mm}

Let us first state that \eqref{Knapsack_pb} is indeed solvable as $\text{D}_c\geq 0$, hence in the worse case $\bs{\Theta}_c$ can be set to zero. The positivity of $\text{D}_c$ is easily verifiable as, through \eqref{modif_FV} and \eqref{num_entr_flux_Tadmor}, it directly follows that $-\(\un{v}_p^c-\un{v}_m^c\)\, \mc{F}_{mp}^\text{\,\tiny FV}\leq -\big(\bs{\Psi}(\un{v}_m^c)-\bs{\Psi}(\un{v}_p^v)\big)\pds\bs{n}_{mp}$, and $\text{D}_c$ can be recast into
\begin{align*}
  \text{D}_c\geq \Sum_{m=1}^{N_s} \(\Sum_{S_{p}^v\,\in\mc{V}_m^{\,c}}l_{mp}\,\bs{n}_{mp}\) \pds \bs{\Psi}(\un{v}_m^c)=0.
\end{align*}
Following the steps of \cite{Chan_knapsack}, we efficiently solved \eqref{Knapsack_pb} through a Greedy algorithm (see \cite{Chan_knapsack} for a detailed algorithm) by finding all the $\theta_{mp}$ under the constraint $0\leq\theta_{mp}\leq\theta_{mp}^e\leq 1$, where the $\theta_{mp}^e$ can be any additional constraint on the blending coefficient, while maximizing $\;\sum_{f_{mp}\in\,\mf{f}_c} \theta_{mp}\,$.

\begin{remark}
  \label{rem_entropy_3}
  While conditions \eqref{theta_entropy_1} and \eqref{theta_entropy_2} would respectively reduced the accuracy to first and second order, condition \eqref{Knapsack_pb} does allow the preservation of the high-order accuracy of the scheme, see for instance Table~\ref{table_order_adv_1D_entropy} and Figure~\ref{fig_advect_entropy_compar}. To make sure of it, let us raise that substituting $u_h$ by a smooth function $u$ leads to \,{\normalfont $\text{\textbf{C}}_c\pds\bs{1}-\text{D}_c=|\omega_c|\,\mc{O}(h_c^{k+1})$}, where $h_c$ is the diameter of cell $\omega_c$. Indeed, after some simple manipulations, it is possible to write that {\normalfont
  \begin{align*}
    \text{\textbf{C}}_c\pds\bs{1}-\text{D}_c=E_{\partial\omega_c}\big(\bs{\phi}(u)\pds\bs{n}\big)+E_{\omega_c}\big(\bs{F}(u)\pds\gradx{} v_h^c\big)-E_{\partial\omega_c}\big(v_h^c\,\bs{F}(u)\pds\bs{n}\big)+\Int_{\omega_c}\hspace*{-1mm}\big(v(u)-v_h^c\big)\,\divx{\bs{F}(u)}\,\dd V,
  \end{align*}}
  where $E_{\Omega}(f)=\oint_\Omega f \,\dd \Omega-\int_\Omega f \,\dd \Omega$. Then, using similar arguments as in \cite{Cockburn_lcs4}, and since we make use of quadrature rules respectively exact for polynomials up to degree $2 k$ over $\omega_c$ and $2 k+1$ over $\partial \omega_c$, we obtain the desired result. Then, by means of the greedy algorithm, see \cite{Chan_knapsack}, it is possible to state that $\(\theta_{mp}-1\)\,\Delta F_{mp}=\frac{|\omega_c|}{l_{mp}}\,\mc{O}(h_c^{k+1})=\mc{O}(h_c^{k+2})$. In the end, as the blended flux can be re-expressed as $\wt{F_{mp}}=\wh{F_{mp}}+\(\theta_{mp}-1\)\,\Delta F_{mp}=\wh{F_{mp}}+\mc{O}(h_c^{k+2})$, the monolithic DG/FV scheme reduces to pure DG, up to $\mc{O}(h_c^{k+2})$, in this smooth solution context.
\end{remark}

\subsection{Numerical results: entropy stabilities}
\label{sect_results_entropy}

To end this section concerned with entropy stabilities and to confirm the previous stated results regarding the different entropy stability and their respective cost in accuracy, we run some numerical tests on some classical problems. Let us emphasize that a lot more problems and test cases will be considered in the section devoted to maximum principles, Section~\ref{sect_results_GLMP}, as the following one is simply devoted to the questions regarding entropy announced in the introduction of this part. In the following, if not stated otherwise, the subcell mean values will be displayed.

\subsubsection{1D linear advection case}

\begin{figure}[!ht]
  \begin{center}
    \includegraphics[height=6.cm]{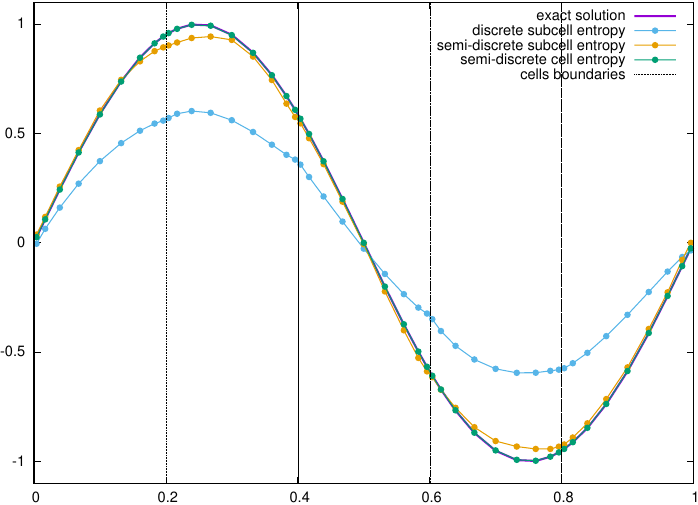}
    \caption{$\poly{8}$-DG/FV solutions on 5 cells with $\eta(u)=\demi u^2$}
    \label{fig_adv_sinus_9th_entropy}
  \end{center}
\end{figure}

First, let us consider the very simple case of 1D linear advection $\vdt u +a\vdx u=0$, where the velocity is set to $a=1$. We start from a smooth initial condition $u_0(x)=\sin(2\pi x)$ and assume periodic boundary conditions. We assess the scheme accuracy after one full period. In Figure~\ref{fig_adv_sinus_9th_entropy}, numerical solutions obtained by means of the $\P^{\,8}$ monolithic DG/FV scheme based on the three different blending coefficients, conditions~\eqref{theta_entropy_1}, \eqref{theta_entropy_2} and \eqref{Knapsack_pb}, ensuring the three different types of entropy stability, are plotted.

As expected, one can see how the first condition does reduce the accuracy to first-order, the second to second-order while the third one is the only one able to preserve the high accuracy, as with only 5 cells the numerical solution is extremely close to the exact one. 
\begin{table}[!ht]
  \begin{center}
    \begin{tabular}{|c||c|c||c|c||c|c|}
      \hline & \multicolumn{2}{c||}{Entropy stability n°1} & \multicolumn{2}{c||}{Entropy stability n°2} & \multicolumn{2}{c|}{Entropy stability n°3}\\ \hline\hline $h$ & $E_{L_2}$ & $q_{L_2}$ & $E_{L_2}$ & $q_{L_2}$ & $E_{L_2}$ & $q_{L_2}$\\
      \hline $1/{1}$ & 6.97E-1 & 0.21 & 5.76E-1 & 1.35 & 1.32E-2 & 4.64\\
      \hline $1/{2}$ & 6.01E-1 & 0.48 & 2.26E-1 & 1.43 & 5.27E-4 & 6.36\\
      \hline $1/{4}$ & 4.29E-1 & 1.54 & 8.40E-2 & 1.41 & 6.41E-6 & 5.98\\
      \hline $1/{8}$ & 2.64E-1 & 0.92 & 3.16E-2 & 1.22 & 1.02E-7 & 5.89\\
      \hline $1/{16}$ & 1.48E-1 & - & 1.36E-2 & - & 1.72E-9 & -\\
      \hline 
    \end{tabular}
  \end{center}\vspace*{-3mm}
    \caption{Convergence rates for the linear advection case for $\poly{5}$-DG/FV monolithic scheme}
  \label{table_order_adv_1D_entropy}
\end{table}
The rates of convergence gathered in Table~\ref{table_order_adv_1D_entropy} further confirm this result. Let us emphasize that, as previously proved, in this smooth solution context the \scheme scheme with the third entropy stability condition reduces to a pure DG scheme, up to machine precision.\\

We follow with the classical case of the linear advection of a composite signal, introduced in \cite{jiang_weno}. This signal is composed by the succession of a Gaussian, rectangular, triangular and parabolic signals.

\begin{figure}[!ht]
  \begin{center}
    \includegraphics[height=7.5cm]{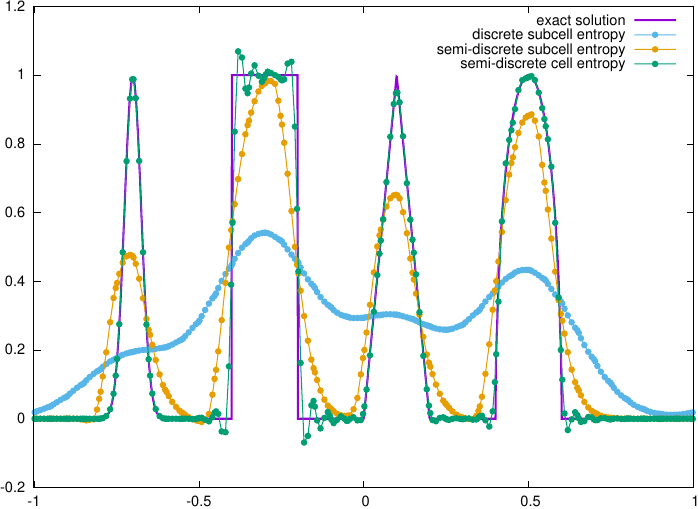}
    \caption{$\poly{5}$-DG/FV solutions on 40 cells and $\eta(u)=\demi u^2$}
    \label{fig_advect_entropy_compar}
  \end{center}
\end{figure}

\begin{figure}[!ht]
  \begin{center}
    \subfigure[$k_e=-0.00001$]{\includegraphics[height=6.cm]{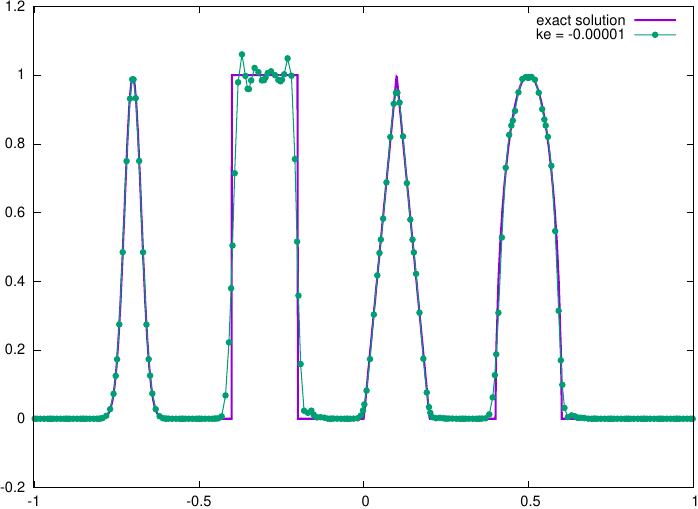}\label{fig_advection_smooth_map}}\hspace*{5.mm}
    \subfigure[$k_e=1.00001$]{\includegraphics[height=6.cm]{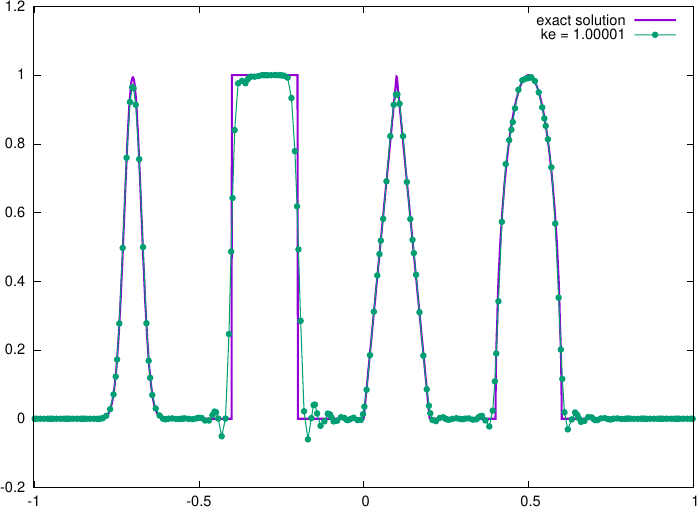}\label{fig_advection_smooth_profile}}
    \caption{$\poly{5}$-DG/FV submean values on 40 cells: $\epsilon=0.25$ and $\eta(u)=|u-k_e|^{1+\epsilon}\backslash(1+\epsilon)$}
    \label{fig_advect_entropy_kruzkov}
  \end{center}
\end{figure}

Figure~\ref{fig_advect_entropy_compar}, in which monolithic $\poly{5}$-DG/FV solutions ensuring the three types of entropy stability are displayed, confirms furthermore our previous conclusion on accuracy and entropy. However, let us note that the numerical solution ensuring the semi-discrete cell entropy inequality, for $\eta(u)=\demi\,u^2$, is very close to what a pure DG scheme would produce and thus exhibits the same pathologies. To be able to see a more significant impact of this entropy inequality enforcement, let us use different entropies $\eta(u)=|u-k_e|^{1+\epsilon}\backslash(1+\epsilon)$. Those can be seen as a smoothed version of the Kruzkov's entropies. The coefficient $\epsilon>0$ is set here at $\epsilon=0.25$, while different values of $k_e$ will be used. In Figure~\ref{fig_advect_entropy_kruzkov}, we make use of two different ones, respectively $k_e=-0.00001$ and $k_e=1.00001$.
One can clearly see how the numerical solution has now been impacted by the semi-discrete cell entropy stability, and how this condition put the emphasis around $u=k_e$ making in the first case the numerical solution almost positive while in the second case roughly less than one.

\subsubsection{1D Buckley non-convex case}

\begin{figure}[!ht]
  \begin{center}
    \subfigure[Test case n°1]{\includegraphics[height=6.cm]{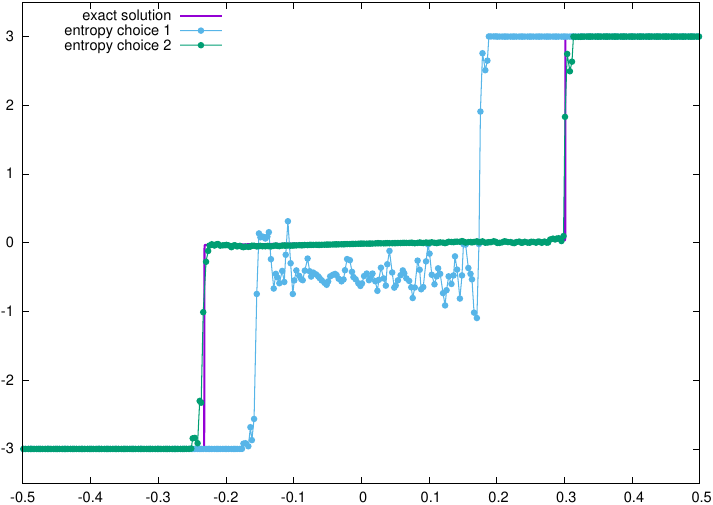}\label{fig_buckley_shu_entr}}\hspace*{5.mm}
    \subfigure[Test case n°2]{\includegraphics[height=6.cm]{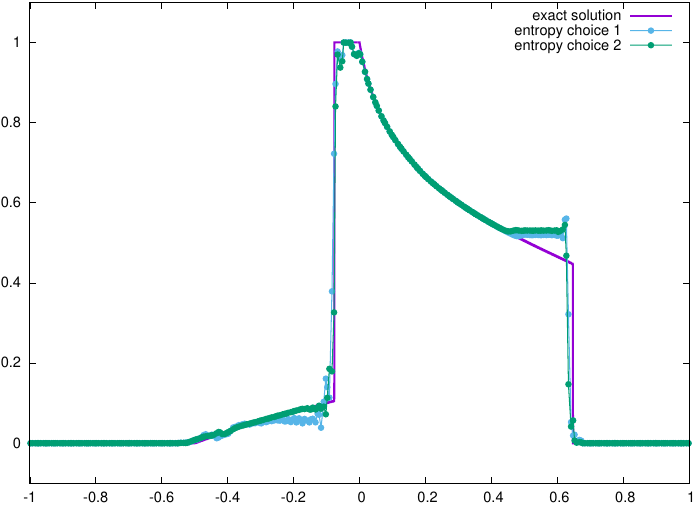}\label{fig_buckley_entr}}
    \caption{$\poly{3}$-DG/FV submean values on 80 cells: $\eta_1(u)=\demi u^2$ and $\eta_2(u)=\int \text{atan}(20 u)\, \dd u$}
    \label{fig_buckley_entr_2}
  \end{center}
\end{figure}

Now, we address the challenging 1D Buckley problem. The Buckley equation is defined as $\vdt{u}+\vdx{F(u)}=0$, where the non-convex flux function writes $F(u)=\frac{4\,u^2}{4\,u^2+(1-u)^2}$. As said in Remark~\ref{quadrature}, since the flux function is now a complex rational function, it is not practical to analytically integrate the volume integrals. And due to that, entropy stability proved in \cite{jiang} does not hold anymore. Furthermore, approximated integration or collocation of the flux may also produce some aliasing effects, see \cite{vilar_aplsc_1D} for some examples. Here, we want to observe the benefit of entropy stability and check if a semi-discrete cell entropy inequality for one given entropy is practically enough to capture the correct unique entropic weak solution. To this end, two different test cases will be addressed. The first one has been introduced by T. Chen and C.-W. Shu in \cite{SBP_shu}. This Riemann problem, consisting in an initial discontinuity located at $x=0$ and taking -3 and 3 as left and right values, admits an entropic solution containing two shock waves connected by a flat rarefaction that is close to 0. We run this test case on a domain $[-0.5,\,0.5]$ and end at time $t=1$. As said in the aforementioned paper, in this case the choice of the entropy function is critical. To corroborate this statement, we make use of the same two entropy functions as in \cite{SBP_shu}, meaning the energy one $\eta(u)=\demi\,u^2$ and $\eta(u)=\int \arctan(20\,u)\,\dd u$ a mollified version of the Kruzkov's entropy $|u|$. To solely observe the effect of entropy stability, and not other choices of blending coefficients ensuring other properties, see Section~\ref{sect_GLMP}, we additionally use here the DG maximum principle preserving limitation of X. Zhang and C.-W. Shu \cite{zshu1} to make sure the computation goes through. This additional limitation will only be used here as in the next section maximum principle and positivity will be ensured directly by appropriate choice of blending coefficients. In Figure~\ref{fig_buckley_shu_entr}, $\poly{3}$-DG/FV monolithic scheme is used on 80 cells, ensuring the high accuracy preserving semi-discrete cell entropy stability, with the two different entropies discussed before.
In Figure~\ref{fig_buckley_shu_entr}, it is clear how the scheme using the energy entropy has failed to capture the entropic weak solution, while using another entropy, putting the stress around 0 because being an approximation of $|u|$, has solved this issue. This is perfectly consistent to the results obtained in \cite{SBP_shu}. Now, making use of the same two entropies, we run a second test case in which we start from the initial solution $u^0(x)=1$ if $x\in\,[-\demi,0]$ and $u^0(x)=0$ elsewhere. The results obtained again by means of $\poly{3}$-DG/FV monolithic scheme is used on 80 cells are shown in Figure~\ref{fig_buckley_entr}. One can clearly see how none of these two choices of entropy has enabled the scheme to capture the entropic solution. A proper entropy has to be designed in some retro-engineering process to fit not only the PDE considered but also the test case to hope to capture the correct solution. And this may be even not feasible if the solution presents very complex structures. Let us emphasize that similar results would have been obtained ensuring a semi-discrete subcell entropy inequality, meaning by means of condition \eqref{theta_entropy_2}, as this entropy stability is ensured only for one chosen entropy. Only the first condition \eqref{theta_entropy_1} enforcing discrete entropy stability for any entropy will succeed in capturing the unique weak solution in both cases. But, as recalled, it will reduce the accuracy to first-order.

\subsubsection{2D KPP non-convex case}

We now turn our attention to the 2D KPP problem proposed by Kurganov, Petrova, Popov (KPP) in \cite{Popov}. For this particular problem, the flux function is given by $\bs{F}(u)=\(\sin(u),\,\cos(u)\)\tra$. Considering the computational domain $[-2,2]\times[-2.5,1.5]$, the initial condition reads as follows

\begin{align*}
  u^0(x)=\left\{\begin{array}{ll}
  \frac{7\,\pi}{2}\qquad &\text{if }\; x<\demi,\\[2mm]
  \frac{\pi}{4}\qquad &\text{if }\; x>\demi.
  \end{array}\right.
\end{align*}

\begin{figure}[!ht]
  \begin{center}
    \subfigure[$1^\text{th}$-order FV on\, 209184 cells]{\includegraphics[height=7.3cm]{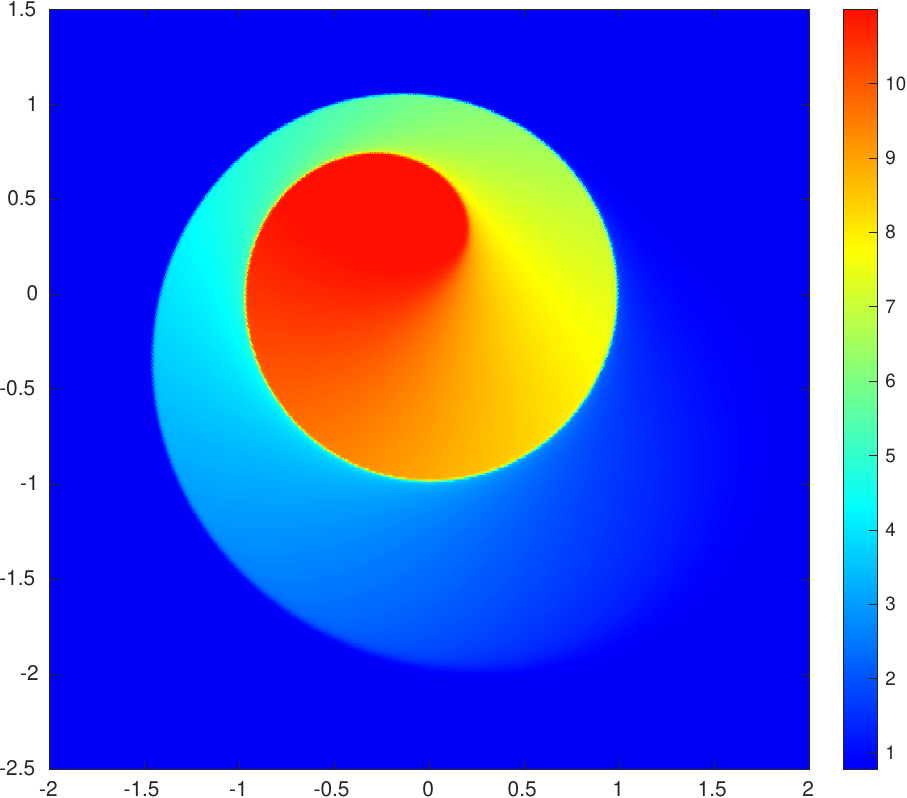}\label{fig_KPP_order1_entr}}\hspace*{5.mm}
    \subfigure[$\P^{\,3}$-DG/FV on 1054 cells]{\includegraphics[height=7.3cm]{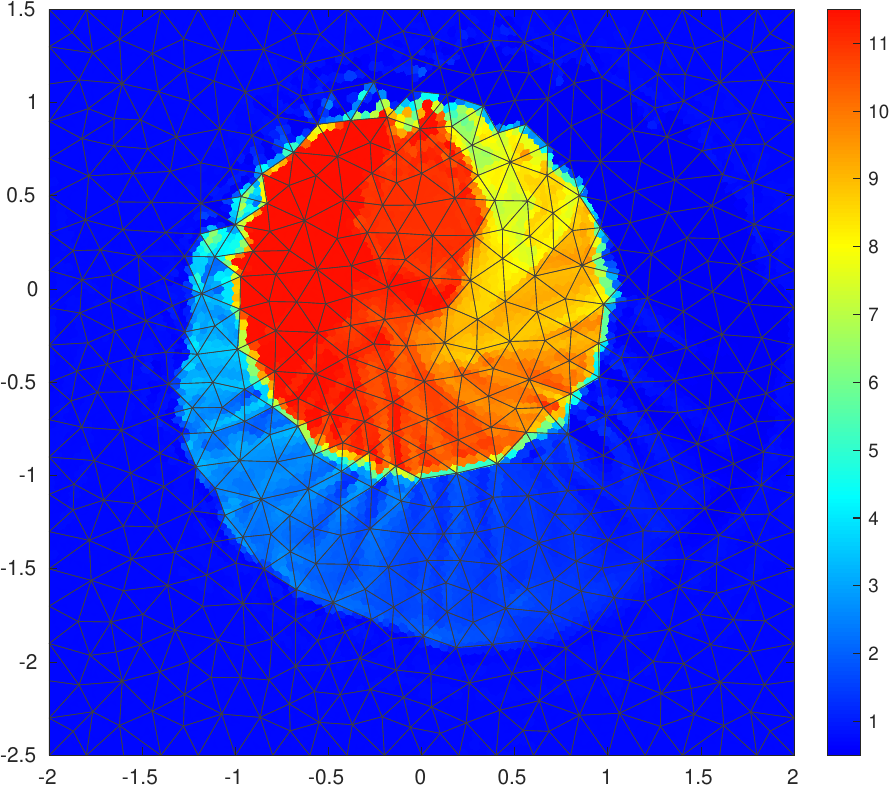}\label{fig_KPP_order4_entr}}
    \caption{$\P^{\,3}$-DG/FV entropic scheme: non-entropic solution}
    \label{fig_KPP_entr}
  \end{center}
\end{figure}

This test is very challenging to many high-order schemes as the solution has a two-dimensional composite wave structure. Thus, generally, to be able to capture such rotation composite structure, very fine grids are used. Here, we compare a referential solution, obtained through a first-order FV scheme on a very fine grid made of 209184 triangular cells, with the one obtained with the $\poly{3}$-DG/FV monolithic scheme with the accuracy preserving semi-discrete cell entropy stability for $\eta(u)=\demi u^2$, on a coarse mesh made of 1054 cells, see Figure~\ref{fig_KPP_entr}. Because that it has been previously observed in \cite{vilar_aplsc_2D} that in this particular case, voronoi-type cell subdivision \ref{fig_tri2} produces slightly better results, this cell sub-partition is then used here.
One can notice how this entropy stability condition does not guarantee the capture of the correct weak entropy solution. And refining the mesh would not cure this issue.

\subsubsection{1D modified Sod shock tube problem}
\label{test_Sod_mdf_entropy}

\begin{figure}[!ht]
  \begin{center}
    \subfigure[Density]{\includegraphics[height=5.7cm]{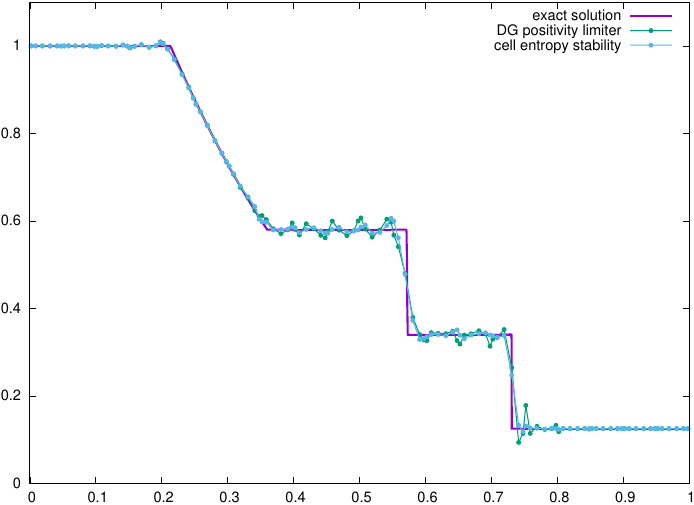}\label{fig_modif_sod_entr_den}}\hspace*{5.mm}
    \subfigure[Velocity]{\includegraphics[height=5.7cm]{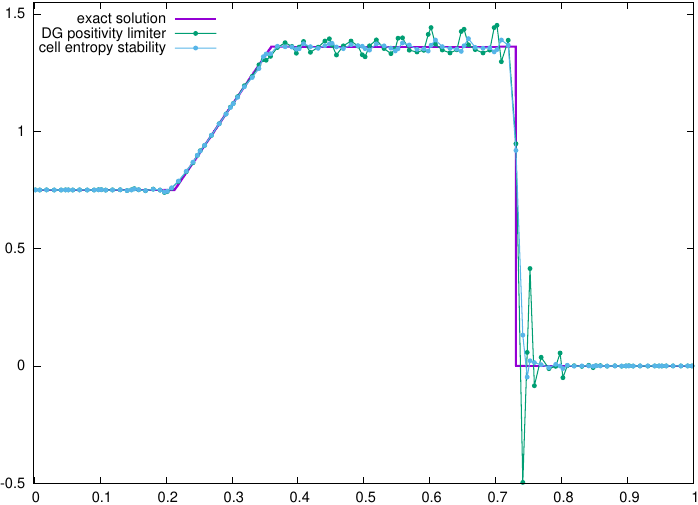}\label{fig_modif_sod_entr_vel}}
    \caption{$\poly{5}$ pure DG with positivity limiter and DG/FV monolithic scheme with cell entropy stability on 20 cells}
    \label{fig_modif_sod_entr}
  \end{center}
\end{figure}

Finally, to end this section regarding entropy stability, we now consider the 1D version of the Euler compressible gas system of equations \eqref{Euler_system}. A classical test case where numerical schemes may capture non-physical weak solution is the modified Sod shock tube problem, \cite{toro_book}. This problem is a modified version of the popular Sod's test \cite{Sod}; the solution consists of a right shock wave, a right traveling contact wave and a left sonic rarefaction wave. This test is very useful in assessing the entropy satisfaction property of numerical methods, as some of them may present a shock at the sonic point in the rarefaction. In Figure~\ref{fig_modif_sod_entr} both solutions obtained by means of $\poly{5}$ pure DG scheme with the X. Zhang and C.-W. Shu positivity-preserving limiter \cite{zshu1} and our $\poly{5}$-DG/FV monolithic scheme with the accuracy preserving semi-discrete cell entropy stability. The chosen entropy in the Euler case is $\eta(\Mat{U})=-\rho\,\log\(p/\rho^\gamma\)$. In both cases, Rusanov type of numerical flux has been used.

Firstly, it is important to note that while pure DG without positivity limiter crashes and thus fails to produce a solution, monolithic scheme with cell entropy stability does run without any need of positivity limiter. Secondly, one can see how the cell entropy stability reduces the amplitude of spurious oscillations. However, let us emphasize that both schemes successfully capture the correct entropic solution, as none of them present the non-physical shock in the rarefaction wave. Similar results can be obtained with the use of global Lax-Friedrichs, HLL and HLL-C numerical fluxes.

\subsection{Conclusion on entropy stabilities?}
\label{sect_conclusion_entropy}

Let us recall the questions we were asking ourselves at the beginning of this section. First, is it possible to find $\theta_{mp}$, the blending coefficients, ensuring entropy stability? The answer is obviously yes. But what do we mean by entropy stability? We have shown how to ensure three types of entropy stabilities, meaning a fully-discrete subcell entropy inequality for any entropy, a semi-discrete subcell entropy one for a given entropy and a semi-discrete cell entropy one for a given entropy. The follow-up question is then what are the costs of such constraints? Well, if one wants to ensure entropy stability for any entropy, the scheme has no choice but to reduce to a first-order one. Relaxing this objective and aiming for an entropy stability for only one entropy, we may achieve second or even arbitrary high-order of accuracy. However, to indeed preserve the high-order accuracy of DG scheme, the cell entropy inequality requires to modify the definition of the subcell's faces FV fluxes in the monolithic DG/FV scheme, see Remark~\ref{rem_modif_FV}, which may invalidate other properties as positivity for instance. Finally, regarding the last question, is entropy stability absolutely needed? In the design of first-order scheme or any robust scheme forming the safe basis of a MOOD-type \cite{mood1,dumbser_subFV_lim_tri} or monolithic scheme, it totally is. A FV scheme with the Roe solver without entropy fix would fail for instance to capture the correct solution in the modified Sod test case \ref{test_Sod_mdf_entropy}. That being said, if the goal is to go to very high order of accuracy, one has no choice but to relax its expectations and aim for an entropy stability for only one entropy. And then, it may always be possible to design a test case, considering a particular scalar conservation law with complex fluxes, that will trick high-order entropy conservative/stable schemes and make them to fail to capture the unique entropic solution. From our experience, while high-order entropy stability does slightly introduce numerical diffusion and consequently reduces spurious oscillations, it is generally not enough to capture the correct entropic solution in complicated cases. An additional shock capturing technique is added, which brings further diffusion and hence does the trick. Now, considering systems of equations, as the Euler compressible gas one, generally no high-order entropy stability is needed as a pure DG scheme with entropic stable numerical fluxes would be enough to capture the entropic solution, see Figure~\ref{fig_modif_sod_entr}. It may not be the case considering complex equations of state, with non-concave pressure for instance, but in classical settings we are not aware of any test case where an high-order DG scheme, based on entropic numerical fluxes as Rusanov or HLL, does capture a non-entropic solution.\\ 

While entropy stability may not be absolutely required, preserving global bounds of some variables generally is, as negative density in the compressible gas case would lead to a crash of the code. The next section aims at defining the blending coefficients to guarantee the numerical solution to remain in a convex admissible set, as well as imposing a local maximum principle to reduce the apparition of spurious oscillations. Those conditions will prove to introduce enough numerical diffusion for the scheme to capture the correct entropic solution as well as to produce the best results.

\section{Global and local maximum principles}
\label{sect_GLMP}

This section is devoted to the definition of blending coefficients to ensure different maximum principles. It is important to emphasize that we do distinguish physical maximum principles, that the unique entropic solution should ensure, from numerical maximum principles used in the following to avoid, as much as possible, the apparition of non-physical oscillations.

\subsection{Physical maximum principles}
\label{subsect_physical_maximum}

In here, we present the minimum requirements on the numerical solution to ensure the simulation code to be robust.

\subsubsection{Scalar conservation laws and global maximum principle}
\label{para_physical_maximum_SCL}

Considering SCL, equation \eqref{lcs_eq2D}, if the initial datum yields $u_0\in[\alpha,\,\beta]$, then the unique entropic weak solution $u$ ensures $u(., t)\in[\alpha,\,\beta]$ for any time $t$. To guarantee that the numerical solution submean values do ensure such property, it is sufficient, by convexity of relation \eqref{convex_combo}, that the blended Riemann intermediate state $\wt{u_{mp}}^\pm$ also remain in $[\alpha,\,\beta]$. Since $u_{mp}^{\ast,\text{\,\tiny FV}}$ does, a sufficient condition is to take $\theta_{mp}$ such that 
\begin{align}
  \label{GMP_SCL}
  \theta_{mp}\leq \min\Big(1,\,\Big|\Frac{\gamma_{mp}}{\Delta F_{mp}}\Big|\,\min\big(\beta-u_{mp}^{\ast,\text{\,\tiny FV}},\,u_{mp}^{\ast,\text{\,\tiny FV}}-\alpha\big)\Big).
\end{align}

\subsubsection{Euler system and positivity preservation}
\label{subsubsect_physical_maximum_Euler}

Let us consider the Euler compressible gas dynamics system
\begin{subequations}
\label{Euler_system}
\begin{alignat}{2}[left = \empheqlbrace\,]
&\vdt{\,\Mat{U}}(\bs{x},t)+\divx{\bs{\Mat{F}}\big(\Mat{U}(\bs{x},t)\big)}=0, \qqquad& (\bs{x},t)\in\,\omega\times[0,T], \label{Euler1}\\
&\Mat{U}(\bs{x},0)=\Mat{U}_0(\bs{x}), &\bs{x}\in\,\omega, \label{Euler2}
\end{alignat}
\end{subequations}
with $\Mat{U}=\(\rho,\,\bs{q},\,E\)\tra$ and $\bs{\Mat{F}}\big(\Mat{U}\big):=\(\bs{F}^\rho,\,\bs{F}^{\bs{q}},\,\bs{F}^E\)\tra=\(\bs{q},\,\bs{v}\otimes\bs{q}+p\,I_d,\,(E+p)\,\bs{v}\)\tra$. The conserved variables $\rho$, $\bs{q}=\rho\,\bs{v}$ and $E$ then respectively stand for the density, momentum and total energy, while $\bs{v}$ characterizes the fluid velocity. The thermodynamic closure is given by the equation of state $p=p(\rho,\, \veps)$ where $\veps=E-\tfrac{1}{2}\rho\Vert \bs{v}\Vert^2$ denotes the internal energy. In this paper, we make use of a gamma gas law, $i.e.$\; $p=(\gamma-1)\,\veps$, where $\gamma$ is the polytropic index of the gas. Although the whole theory presented here has been introduced in the simple case of scalar conservation laws, the extension to the system case is perfectly straightforward.\\

Now, defining the admissible convex set $G=\big\{\Mat{U}=\(\rho,\,\bs{q},\,E\)\tra, \quad \rho> 0,\,p>0\big\}$, we want to ensure that any subcell mean value $\ov{u}_m^c$ remains in $G$ at all time. Following similar steps than in \cite{Kuzmin_monolithic_2020,Kuzmin_Gassner_monolithic_2024}, we first ensure the positivity of the density and then of the internal energy. Firstly, we introduce a first temporary blending coefficient $\theta_{mp}^{(1)}$ such that
\begin{align}
  \label{positivity_density}
  \theta_{mp}^{(1)}\leq \min\Big(1,\,\Big|\Frac{\gamma_{mp}}{\Delta F_{mp}^\rho}\Big|\,\rho_{mp}^{\ast,\text{\,\tiny FV}}\big)\Big).
\end{align}
Then, defining quantities $A_{mp}$, $B_{mp}$ and $M_{mp}$ as in the following
\begin{align}
  \label{positivity_A_B_M}
  \left\{\begin{array}{l}
   A_{mp}=\Frac{1}{(\gamma_{mp})^2}\,\(\demi\,\big\Vert\Delta \bs{F}_{mp}^{\bs{q}}\big\Vert^2-\theta_{mp}^{(1)}\,\Delta F_{mp}^\rho\,\Delta F_{mp}^E\),\\[5mm]
   B_{mp}=\Frac{1}{\gamma_{mp}}\,\(\bs{q}_{mp}^{\ast,\text{\,\tiny FV}}\pds\Delta \bs{F}_{mp}^{\bs{q}}-\rho_{mp}^{\ast,\text{\,\tiny FV}}\,\Delta F_{mp}^E-\theta_{mp}^{(1)}\,E_{mp}^{\ast,\text{\,\tiny FV}}\,\Delta F_{mp}^\rho\),\\[5mm]
   M_{mp}=\rho_{mp}^{\ast,\text{\,\tiny FV}}\,E_{mp}^{\ast,\text{\,\tiny FV}}-\demi\,\big\Vert \bs{q}_{mp}^{\ast,\text{\,\tiny FV}}\big\Vert^2,
  \end{array}\right.
\end{align}

where $M_{mp}>0$ in the case of a positivity-preserving FV scheme, we introduce a second temporary blending coefficient as
\begin{align}
  \label{positivity_energy}
  \theta_{mp}^{(2)}\leq \min\(1,\,\Frac{M_{mp}}{\big|B_{mp}\big|+\max\big(0, A_{mp}\big)}\).
\end{align}
Finally, the blending coefficients will be defined as the product of these two, $i.e.$\; $\theta_{mp}=\theta_{mp}^{(1)}\,\theta_{mp}^{(2)}$.

\begin{remark}
  \label{remark_positivity_HLL-C}
  The previous formula have been given for a numerical flux of form \eqref{num_flux}, which includes global Lax-Friedrichs and Rusanov fluxes, producing a single FV intermediate state $\Mat{U}_{mp}^{\ast,\text{\,\tiny FV}}$. Everything can be easily extended to other types of numerical fluxes, as HLL or HLL-C fluxes for instance. In those latter cases, for sake of simplicity, we can set $\gamma_{mp}$, present in the previous formula, as $\gamma_{mp}=\Max(| S_{mp}^L|, | S_{mp}^R|)$, where $S_{mp}^L$ and $S_{mp}^R$ are the smallest and highest velocities of the two acoustic waves, see \cite{batten} for a proper definition of those velocities to ensure a positivity-preserving behavior. Doing so, the first-order FV numerical flux will then produce two FV Riemann intermediate states $\Mat{U}^{\ast,\,\pm}_{mp}$, left and right, being defined as a convex combination of the HLL or HLL-C intermediate states and the initial left and right states $\ov{\Mat{U}}_m^c$ and $\ov{\Mat{U}}_p^v$.
\end{remark}

\begin{remark}
  \label{remark_NaN}
  It is essential to emphasize that those maximum or positivity principles impose the subcell mean values to remain in $G$, the admissible set. However, nothing is said on the values of the solution polynomial reconstruction, required in the computation of the DG residuals to define the reconstructed fluxes. It is thus perfectly possible that the polynomial solution in a cell yields a non-admissible value, at a cell interface for instance, which will lead to non-admissible reconstructed fluxes (even possibly \textit{NaN} values). Obviously, this situation is automatically treated in this monolithic framework as, if the reconstructed flux $\wh{F_{mp}}$ presents a pathological value, as \textit{NaN} for instance, the blending coefficient will be set to zero. In the end, this means that if the uncorrected DG solution is nowhere to be saved inside the cell and the DG code would have then crashed, the monolithic scheme will then reduced to a first-order FV scheme applied on each subcell contained in the pathological cell, see for instance the results obtained for the Mach 20 hypersonic flow test case \ref{test_Mach_20}.
\end{remark}

\subsection{Local maximum principles}
\label{subsect_local_maximum}

Now, to reduce as much as possible the apparition of spurious oscillations in the approximation of discontinuous solution, we choose to impose a local maximum principle, at the subcell level.

\subsubsection{Scalar conservation laws and local maximum principle}
\label{subsubsect_local_maximum_SCL}

For SCL, thanks to their hyperbolic nature, the solution at a point should remained bounded by the minimum and maximum values of the solution at a previous time, taken in a large enough domain including the point under consideration. To ensure a low oscillatory behavior of the numerical solution, we will mimic such property at the discrete level. To do so, we will impose the submean value $\ov{u}_m^{c,\,n+1}$ on subcell $S_m^c$ to be bounded by the submean values at the previous time step (or the previous RK step in the general case) in a given subcells set, as
\begin{align}
  \label{NAD}
\alpha_m^c:=\Min_{S_q^w\,\in \,\mc{N}(S_m^c)}\big(\ov{u}_q^{\,w,\,n}\big) \;\leq\;  \ov{u}_m^{\,c,n+1} \;\leq\; \Max_{S_q^w\,\in \,\mc{N}(S_m^c)}\big(\ov{u}_q^{\,w,\,n}\big):=\beta_m^c,
\end{align}
where $\mc{N}(S_m^c)$ is some set of $S_m^c$ neighboring subcells, including subcell $S_m^c$, yet to be defined. The wider the set $\,\mc{N}(S_m^c)$ is, the softer this local maximum principle will be. Reversely, a smaller set would lead to a larger first-order FV contribution in flux blending. In this work, similarly to the one used in the non-linear numerical results section of \cite{vilar_aplsc_2D}, $\mc{N}(S_m^c)$ will be constituted by subcell $S_m^c$, as well as all its face and node neighboring subcells $S_q^v$, either they belong to the same cell or not. By introducing $\mc{P}_m^c$ the set of vertices of subcell $S_m^c$, as well as $\mc{V}_p$ the set of subcells that share $\bs{x}_p$ as a vertex, $i.e. \;\; S_p^v\in\mc{V}_q\; \Longrightarrow\; \bs{x}_q\in \mc{P}_p^v$, this definition of $\mc{N}(S_m^c)$ can be rewritten as $ \mc{N}(S_m^c)=\bigcup_{ \bs{x}_p\,\in\, \mc{P}_m^c}\,\mc{V}_p$.\\

Now, to guarantee condition \eqref{NAD} for any subcell, in the light of relation \eqref{convex_combo} it is sufficient to ensure that $\wt{u_{mp}}^-\in[\alpha_m^c,\,\beta_m^c]$ along with $\wt{u_{mp}}^+\in[\alpha_p^v,\,\beta_p^v]$. Since $u_{mp}^{\ast,\text{\,\tiny FV}}$ insures both conditions, it is sufficient to take $\theta_{mp}$ such that 
\begin{align}
  \label{LMP_SCL}
  \theta_{mp}\leq \min\(1,\,\Big|\Frac{\gamma_{mp}}{\Delta F_{mp}}\Big|\,
  \left\{\begin{array}{ll}\min\big(\beta_p^v-u_{mp}^{\ast,\text{\,\tiny FV}},\,u_{mp}^{\ast,\text{\,\tiny FV}}-\alpha_m^c\big) \quad&\text{if } \Delta F_{mp}>0\\[3mm]\min\big(\beta_m^c-u_{mp}^{\ast,\text{\,\tiny FV}},\,u_{mp}^{\ast,\text{\,\tiny FV}}-\alpha_p^v\big) \qquad&\text{if } \Delta F_{mp}<0
  \end{array}\right.\).
\end{align}

\begin{remark}
  \label{rem_smooth_relax}
  Let us enlighten that the local maximum principle \eqref{NAD} relies on subcell mean values. And because this constraint is not concerned with the whole polynomial set of values, it is very well-known that one has to relax it to preserve scheme accuracy in the presence of smooth extrema.
\end{remark}

\paragraph{Smooth extrema relaxation}
In order to preserve high-accuracy in the vicinity of smooth extrema, we make use of a subcell level version of the smooth detector we introduced in our previous article, \cite{vilar_aplsc_2D}. This latter is closed related to the smoothness indicator for finite elements, \cite{Kuzmin_smooth}. The basic idea of this detector is the following: the numerical solution is supposed to exhibit a smooth extrema if at least the linearized version of the numerical solution spatial derivatives present a monotonous profile. To this end, let us introduce the following subcell linear reconstructions

\begin{subequations}
\label{grad_uh}
\begin{alignat}{2}[left = \empheqlbrace\,]
&v_x^{\,m}(\bs{x})=\ov{\vdx \,u_h^{c}}^{\,m}+\ov{\gradx\, (\vdx\, u_h^{c})}^{\,m}\pds (\bs{x}-\bs{x}_m^c), \label{dx_uh}\\[3mm]
&v_y^{\,m}(\bs{x})=\ov{\vdy \,u_h^{c}}^{\,m}+\ov{\gradx\, (\vdy\, u_h^{c})}^{\,m}\pds (\bs{x}-\bs{x}_m^c). \label{dy_uh}
\end{alignat}\\[-4mm]
\end{subequations}

In \eqref{grad_uh}, $\bs{x}_m^c$ denotes the centroid of subcell $S_m^c$, while $\ov{\partial_\xory \,u_h^{c}}^{\,m}$ and $\ov{\gradx\, (\partial_\xory\, u_h^{c})}^{\,m}$ are nothing but the averaged values on $S_m^c$ of the successive partial derivatives of $u_h^c$. In practice, this smoothness indicator works as a vertex-based limiter on $v_\xory^{\,m}(\bs{x})$. Due to their linearity, functions $v_\xory^{\,m}(\bs{x})$ attain their extrema at the vertices $\bs{x}_q\in\mc{P}_m^c$. Then, we consider that the exact weak solution underlying the numerical solution $u_h$ presents a smooth profile in subcell $S_m^c$ if, for any vertex $\bs{x}_q\in\mc{P}_m^c$, the linearized spatial derivative functions ensure the following constraint
\begin{align}
  \label{smooth_constraint}
  v_{x,\,q}^{\min}\leq v_x^m(\bs{x_q})\leq v_{x,\,q}^{\max} \qquad \text{and} \qquad v_{y,\,q}^{\min}\leq v_y^m(\bs{x_q})\leq v_{y,\,q}^{\max},
\end{align}
where $v_{\xory,\,q}^{\min}$ and $v_{\xory,\,q}^{\max}$ are respectively defined as
\begin{align}
  \label{min_max}
  v_{\xory,\,q}^{\min}=\Min_{v\in \mc{V}_q} v_\xory^{\,m}(\bs{x}_q)\qquad \text{and} \qquad v_{\xory,\,q}^{\max}=\Max_{v\in \mc{V}_q} v_\xory^{\,m}(\bs{x}_q).
\end{align}
Practically, if for any vertex $\bs{x}_q\in\mc{P}_m^c$, conditions \eqref{smooth_constraint} are ensured, we then consider that the numerical solution presents a smooth profile on subcell $S_m^c$. Finally, if both the solution is considered smooth in both subcells $S_m^c$ and $S_p^v$, the blended coefficient constraint through the local maximum condition \eqref{NAD} is relaxed. This procedure allows in practice the preservation of smooth extrema along with the order of accuracy for smooth problems, see Section~\ref{sect_results_GLMP}. Let us emphasize that this smooth extrema is not needed for second-order approximation. Furthermore, considering third-order \scheme scheme, the smoothness detection has to be performed at the cell level, instead of the subcell level, as in this case the second derivatives $\gradx\, (\partial_\xory\, u_h^{c,\,n})$ are constant over the cell.

\subsubsection{Euler system and local maximum principle}
\label{subsubsect_local_maximum_Euler}

Regarding the local maximum principle previously introduced, the natural system counterpart would be to apply the previous criteria to the Riemann invariants. However, in the non-linear system case, those quantities are not easy to get nor to manipulate. We could have used a linearized version of the Riemann invariants, as in \cite{vilar_lag_14} for instance, but for sake of simplicity we naively apply the local maximum principle to one of the conserved variable. This local maximum principle is then relaxed, by means of the same smooth extrema detector previously introduced, but this time based on the chosen conservative variable. In the numerical results Section~\ref{sect_results_GLMP}, we choose to either work with the density or the total energy, as these physical quantities would be sensitive to any type of wave. Also, as theoretically there is no local maximum for the conserved variables, we add the FV Riemann intermediate state to the local bounds, as

\begin{align}
  \label{NAD_Euler}
\alpha_m^c:=\Min_{S_q^w\,\in \,\mc{N}(S_m^c)}\Big(\ov{v}_q^{\,w,\,n},\,v_{mq}^{\ast,\text{\,\tiny FV}}\Big) \;\leq\;  \ov{v}_m^{\,c,n+1} \;\leq\; \Max_{S_q^w\,\in \,\mc{N}(S_m^c)}\Big(\ov{v}_q^{\,w,\,n},\,v_{mq}^{\ast,\text{\,\tiny FV}}\Big):=\beta_m^c,
\end{align}

where $v\in\,\Big\{\rho,\, q^x,\,q^y,\, E\Big\}$. The blending coefficient $\theta_{mp}$ is then defined as previously \eqref{LMP_SCL}. Similarly to the SCL case, this local maximum principle has to be relaxed to preserve accuracy in the presence of smooth extrema.

\subsection{Numerical results: global and local maximum principles}
\label{sect_results_GLMP}

Similarly to the numerical results section devoted to entropy stability, Section~\ref{sect_results_entropy}, we make use here of several widely addressed and challenging test cases to demonstrate the performance and robustness of this \scheme scheme ensuring Global Maximum Principle and positivity in the system case (GMP), as well as a relaxed Local Maximum Principle (LMP) to reduce the apparition of non-physical oscillations. In all following test cases, if not stated otherwise, the simple case of global Lax-Friedrichs numerical flux will be used for both the DG scheme and subcell first-order FV one. Regarding the cell decomposition into subcells, as it has been observed in \cite{vilar_aplsc_2D}, this does not have a major influence on the quality of the numerical results, especially in the non-linear case. Consequently, the simple case of the quad/tri subdivision, Figure~\ref{fig_tri1}, will be used if not specified differently. Last, for the 2D non-linear problems, we make use of the blending coefficients smoothening procedures, see Section~\ref{smoother}. The less diffusive one, introduced in Section~\ref{subsect_blended_fluxes}, will be used if not stated otherwise. 

\subsubsection{1D linear advection case}

\begin{figure}[!ht]
  \begin{center}
    \includegraphics[height=7.5cm]{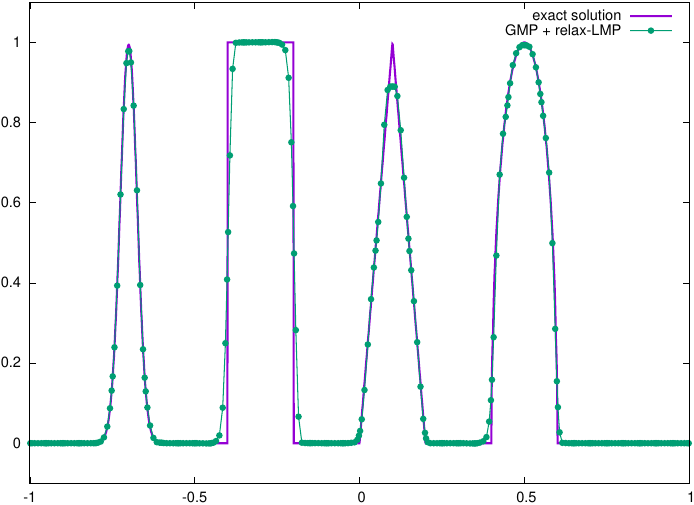}
    \caption{$\poly{6}$-DG/FV with GMP and relaxed-LMP on 40 cells}
    \label{fig_advect_GLMP}
  \end{center}
\end{figure}

As for the entropy stability section, we first consider the very simple case of 1D linear advection of composite signal. In Figure~\ref{fig_advect_GLMP}, the numerical solution obtained by means of the $\poly{6}$ monolithic DG/FV scheme ensuring GMP and a relaxed-LMP one, on a coarse grid made of only 40 cells, is displayed.
One can see in Figure~\ref{fig_advect_GLMP} how the scheme behaves, producing an extremely accurate solution while ensuring the preservation of the global maximum principle and a very low oscillatory profile. We can also observe how the smooth extrema relaxation has permitted to accurately capture the smooth part of the solution. 

\subsubsection{1D Buckley non-convex case}

\begin{figure}[!ht]
  \begin{center}
    \subfigure[Test case n°1]{\includegraphics[height=6.cm]{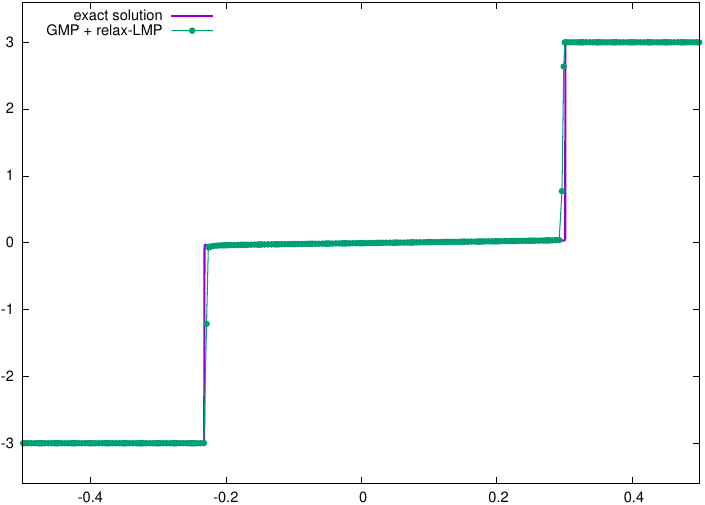}\label{fig_buckley_shu}}\hspace*{5.mm}
    \subfigure[Test case n°2]{\includegraphics[height=6.cm]{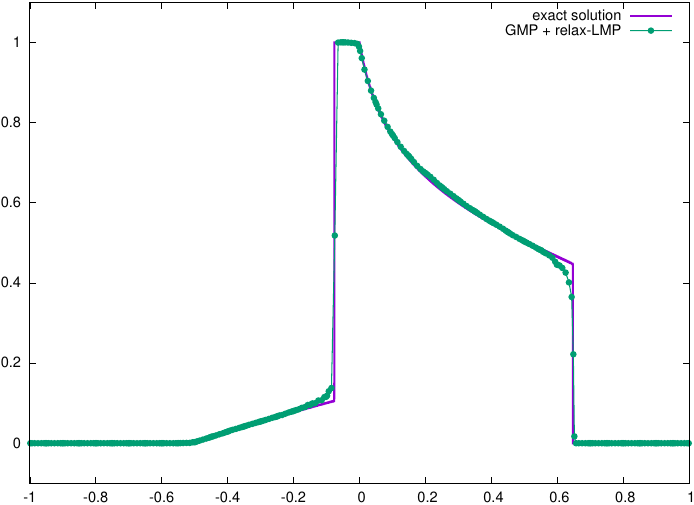}\label{fig_buckley}}
    \caption{$\poly{6}$-DG/FV with GMP and relaxed-LMP on 40 cells}
    \label{fig_buckley_2}
  \end{center}
\end{figure}

Now, by means of the same two test cases used previously in the context of 1D non-convex Buckley SCL, we display in Figure~\ref{fig_buckley_2} the numerical solutions obtained through the GMP and relaxed-LMP $\poly{6}$ monolithic DG/FV using only 40 cells.
While in Section~\ref{sect_results_entropy} we have shown that entropy stability for a given entropy is generally not enough to capture the unique entropic solution, in Figure~\ref{fig_buckley_2} one can see how the two maximum principles imposed here allow the very accurate and robust resolution of the problem under consideration, as even in this very coarse grid context the numerical solutions are extremely close to the exact entropic solutions.

\subsubsection{2D non-linear Burgers problem case}

\begin{figure}[!ht]
  \begin{center}
    \subfigure[Solution submean values]{\includegraphics[height=6.cm]{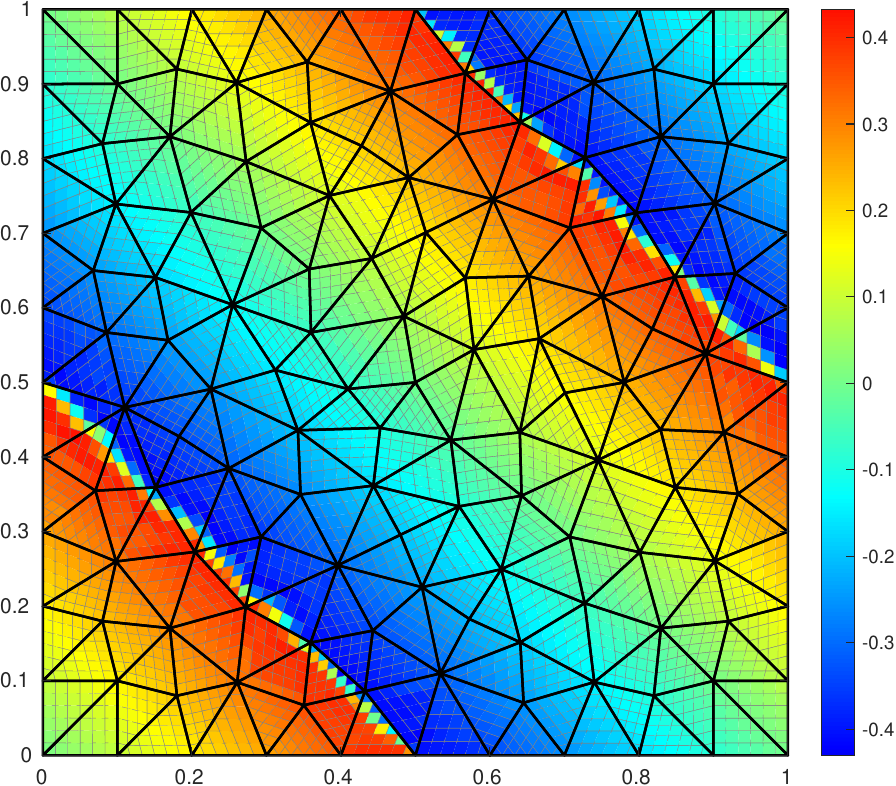}\label{fig_burgers_sol}}\hspace*{5.mm}
    \subfigure[Blending coefficients]{\includegraphics[height=6.cm]{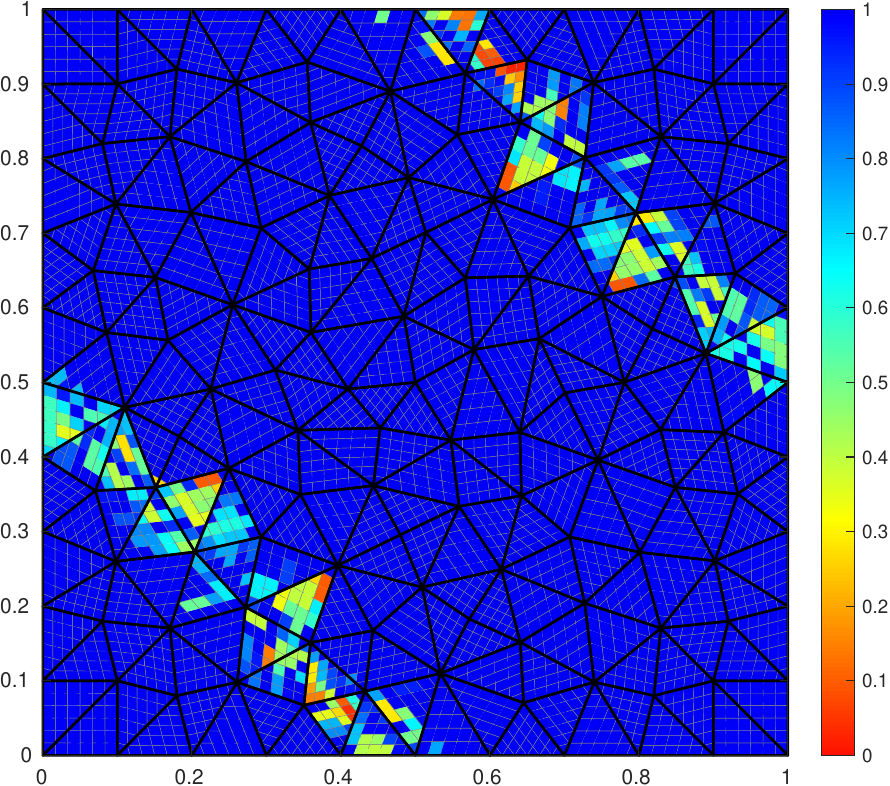}\label{fig_burgers_theta}}
    \caption{$\poly{5}$-DG/FV scheme with GMP and relaxed-LMP on 242 cells}
    \label{fig_burgers}
  \end{center}
\end{figure}

Now, to illustrated how the monolithic fluxes blending operates, let us make use of the Burgers equation, defined through \eqref{lcs1} and flux function $\bs{F}(u)=\demi\(u^2,\,u^2\)\tra$, with the smooth initial solution $u_0(\bs{x})=\sin(2 \pi \,(x+y))$. The domain is chosen as the unit square $[0,\,1]^2$ with periodic boundary condition. Through time, the exact solution will exhibit two stationary shocks along the lines defined by $\(\bs{x}\in[0,\,1]^2,\;x+y=0.5\)$ and $\(\bs{x}\in[0,\,1]^2,\;x+y=1.5\)$. We run this test case until $t=0.5$ with a sixth-order monolithic DG/FV scheme ensuring GMP and a relaxed-LMP, on a very coarse unstructured grid made of 242 cells. In Figure~\ref{fig_burgers_sol}, we display the subcells' mean values while in Figure~\ref{fig_burgers_theta} the subcell blending coefficients, see definition \eqref{subcell_blending_coeff}, are shown.

First, Figure~\ref{fig_burgers_theta} illustrates very well how the monolithic DG/FV scheme works and is able to accurately capture the discontinuities, as only the subcells in a small vicinity of the shocks will be computed trough a convex blending of high-order DG reconstructed fluxes and low-order FV fluxes. Elsewhere the blending coefficients are automatically set to one, which tells us that only the high-order reconstructed fluxes, which gives the equivalency with a pure DG scheme, are used. 

\subsubsection{2D KPP non-convex case}

\begin{figure}[!ht]
  \begin{center}
    \includegraphics[height=7.5cm]{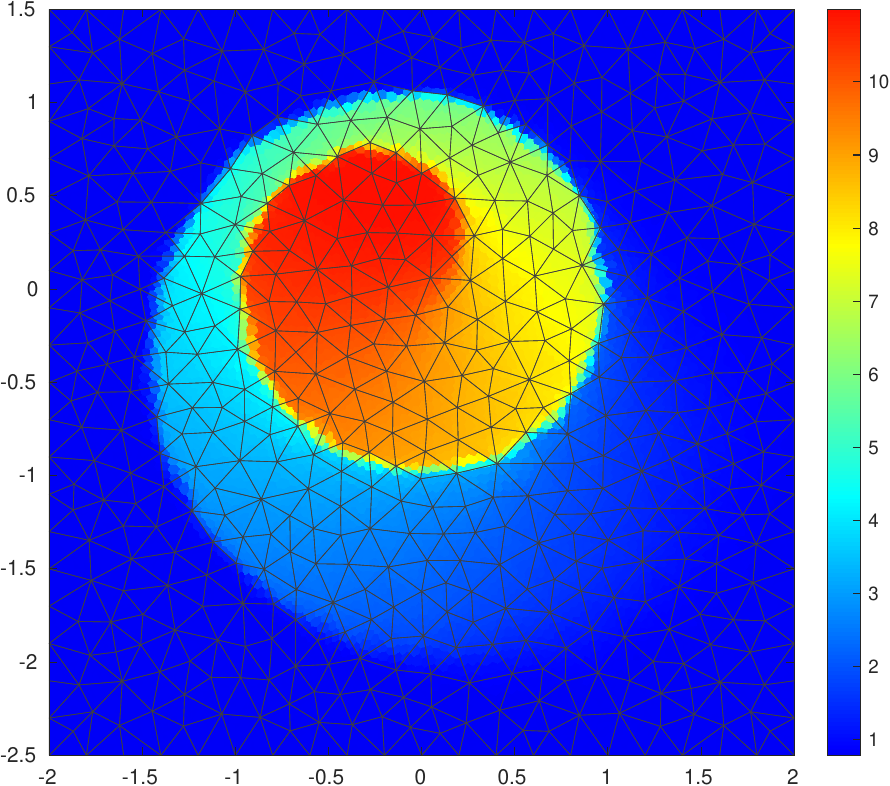}
    \caption{$\poly{3}$-DG/FV scheme with GMP and relaxed-LMP on 1054 cells}
    \label{fig_KPP_GLMP}
  \end{center}
\end{figure}

Now, we consider once more the non-convex flux case of KPP SCL. A similar set up is used here, but whilst in Figure~\ref{fig_KPP_order4_entr} cell entropy stability was enforced, in Figure~\ref{fig_KPP_GLMP} we make use of GLM and relaxed-LMP, combined with the blending coefficient smoother introduced in Section~\ref{smoother}. 
While we have seen that the high-order entropy stable monolithic scheme fails to capture the entropic solution, see Figure~\ref{fig_KPP_order4_entr}, here the two maximum principles, GMP and relaxed-LMP, allows the correct approximation of the unique entropic solution, even in the difficult context of high-order schemes and coarse grids.

\subsubsection{1D modified Sod shock tube problem}

\begin{figure}[!ht]
  \begin{center}
    \subfigure[Density]{\includegraphics[height=5.7cm]{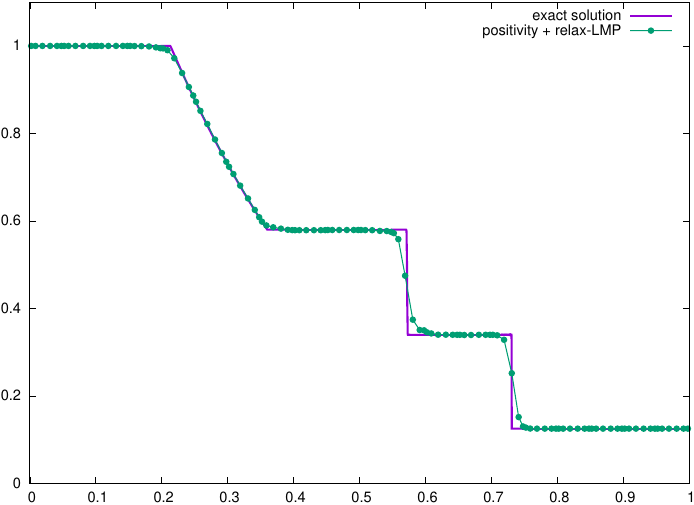}\label{fig_modif_sod_den}}\hspace*{5.mm}
    \subfigure[Velocity]{\includegraphics[height=5.7cm]{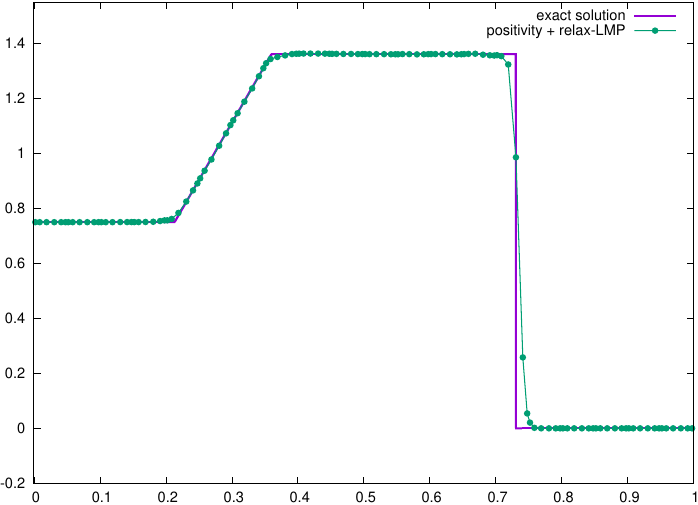}\label{fig_modif_sod_vel}}
    \caption{$\poly{5}$-DG/FV scheme with positivity and relaxed-LMP on 20 cells}
    \label{fig_modif_sod}
  \end{center}
\end{figure}

\begin{figure}[!ht]
  \begin{center}
    \includegraphics[height=7.5cm]{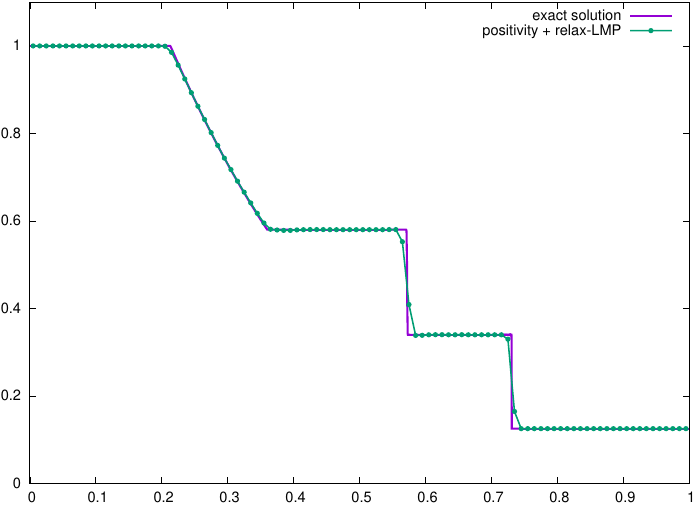}
    \caption{$\poly{2}$-DG/FV scheme with positivity and relaxed-LMP on 100 cells: cells' mean values}
  \label{fig_sod_mdf_order3}
  \end{center}
\end{figure}

Once more, we make use of the modified Sod shock tube problem. In Section~\ref{test_Sod_mdf_entropy}, we have seen how even a pure DG scheme with a positivity-preserving limiter was able to capture the correct solution and that entropy stability increases robustness and reduces spurious oscillations. Now, in a similar configuration, instead of entropy stability, we display the solution obtained through the $\poly{5}$ monolithic DG/FV scheme ensuring positivity and relaxed-LMP, see Figure~\ref{fig_modif_sod}.
Compared to Figure~\ref{fig_modif_sod_entr}, one can observed on Figure~\ref{fig_modif_sod} how those two principles allow us to obtain excellent results, as the numerical solution is non-oscillatory and extremely close to the exact one, even in this extreme coarse grid and very high order context. Let us however emphasize that this \scheme scheme is obviously not limited to the case of very high-order of accuracy on coarse grids. It also preforms very well at second or third order, see Figure~\ref{fig_sod_mdf_order3}.
In the light of Figure~\ref{fig_sod_mdf_order3}, one can see how the third-order monolithic scheme on 100 cells, with positivity and relaxed-LMP conditions, does produce a numerical solution very close to the exact one.

\subsubsection{1D smooth isentropic solution}

\begin{figure}[!ht]
  \begin{center}
    \includegraphics[height=7.5cm]{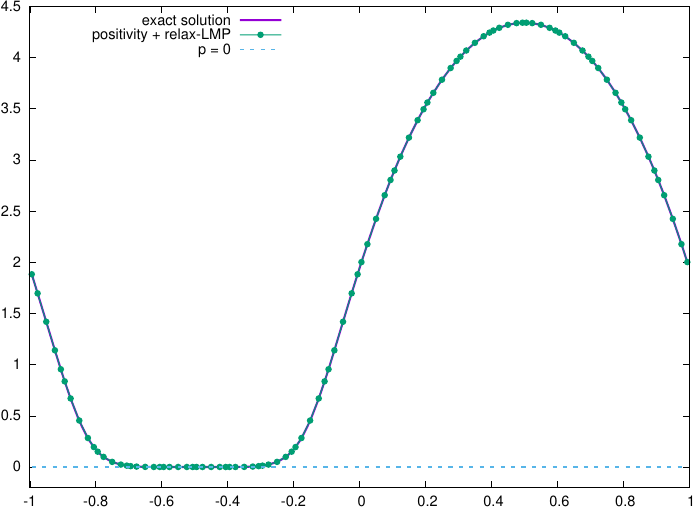}
    \caption{$\poly{4}$-DG/FV scheme with positivity and relaxed-LMP on 20 cells: pressure}
    \label{fig_smooth_Euler}
  \end{center}
\end{figure}

To test the accuracy of the \scheme scheme in the case of system, we make use of a smooth test case initially introduced in \cite{vilar_lag_10}. This example has been derived in the isentropic case, for the perfect gas equation of state with the polytropic index $\gamma=3$. In this special situation, the characteristic curves of the Euler equations become straight lines and the governing equations reduce to two Burgers equations. It is then simple to solve analytically this problem. Here, similarly to \cite{vilar_aplsc_1D}, we modify the initial data to yield a more challenging example, as
\begin{align*}
\rho^0(x)=1+0.9999999 \sin(\pi x),\quad u^0(x)=0,\quad p^0(x)=\rho^0(x)^{\gamma}, \quad x\in [-1,1],
\end{align*}
provided with periodic conditions. This means that initially $\rho^0(-\demi)=1.E-7$ and $p^0(-\demi)=1.E-21$. The density and pressure being so close to zero, any numerical scheme not ensuring a positivity preservation would fail. This is the case of unlimited DG schemes. In Figure~\ref{fig_smooth_Euler}, numerical solutions obtained by means of the $\poly{4}$ monolithic DG/FV scheme with positivity and relaxed-LMP are depicted at time $t=0.1$, using 20 cells.
Figure~\ref{fig_smooth_Euler} demonstrates how robust the monolithic scheme is as no loss of positivity is possible due to the theory, and how accurate the scheme is, the numerical solution being extremely close to the exact one even with only 20 cells. In Table~\ref{table_euler_ordre_5}, we gather the global errors and rates of convergence related to the 5th order scheme, along with the global minimum and the average, in space over the whole domain and in time covering the whole calculation, of the blending coefficients. The results confirm the expected fifth-order rate of convergence, even though the solution has been locally corrected.

\begin{table}[!ht]
  \begin{center}
    \begin{tabular}{|c||c|c||c|c||c|c|}
      \hline & \multicolumn{2}{c||}{$L_1$} & \multicolumn{2}{c||}{$L_2$} & \multicolumn{2}{c|}{$\theta_{mp}$}\\
      \hline $h$ & $E^h_{L_1}$ & $q^h_{L_1}$ & $E^h_{L_2}$ & $q^h_{L_2}$ & min. $\theta_{mp}$ & aver. $\theta_{mp}$\\
      \hline $\frac{1}{10}$ & 9.07E-4 & 5.86 & 1.23E-3 & 5.90 & 2.00E-1 & 0.981 \\
      \hline $\frac{1}{20}$ & 1.56E-5 & 4.03 & 2.05E-5 & 3.83 & 1.92E-1 & 0.997 \\
      \hline $\frac{1}{40}$ & 9.53E-7 & 4.89 & 1.44E-6 & 4.85 & 5.65E-4 & 0.999 \\
      \hline $\frac{1}{80}$ & 3.21E-8 & 4.80 & 5.00E-8 & 4.87 & 3.48E-5 & 0.999\\
      \hline $\frac{1}{160}$ & 1.15E-9 & - & 1.71E-9 & - & 1.00 & 1.00 \\
      \hline 
    \end{tabular}
  \end{center}
  \caption{Convergence rates computed on the pressure for the $\poly{4}$-DG/FV scheme with positivity and relaxed-LMP}
  \label{table_euler_ordre_5}
\end{table}

We can also notice in the light of Table~\ref{table_euler_ordre_5} that refining the mesh, the amount of first-order FV required to stabilize the scheme decreases more and more, as the minimum and averaged blending coefficients tend to one. One could have expected such conclusion since, in this smooth solution context, pure DG is expected to converge to the exact solution. Hence no blending should be required in the end.

\subsubsection{2D Sod shock tube problem}

To close this numerical application section and assess once again the high capability of the \scheme method presented here, the 2D Euler compressible gas dynamics system \eqref{Euler_system} case will be now addressed. First, we consider the extension of the classical Sod shock tube \cite{Sod} to the case of the cylindrical geometry. This problem consists of a cylindrical shock tube of unity radius. The interface is located at $r=0.5$. At the initial time, the states on the left and on the right sides of the interface are constant. The left state is a high pressure fluid characterized by $(\rho_0^L,p_0^L,\bs{v}_0^L)=(1,1,\bs{0})$, the right state is a low pressure fluid defined by $(\rho_0^R,p_0^R,\bs{v}_0^R)= (0.125,0.1,\bs{0})$. The gamma gas law is defined by $\gamma=\frac{7}{5}$. The computational domain is defined in polar coordinates by $(r,\theta)\in [0,1]\times [0,\frac{\pi}{4}]$. We prescribe symmetry boundary conditions at the boundaries $\theta=0$ and $\theta=\frac{\pi}{4}$, and an outflow condition at $r=1$. The exact solution consists of three circular waves, a shock followed by a contact discontinuity and rarefaction wave. The aim of this test case is then to assess the \scheme scheme accuracy while ensuring a non-oscillatory behavior, and its ability to preserve the radial symmetry.
\begin{figure}[!ht]
  \begin{center}
    \subfigure[Density map]{\includegraphics[height=5.45cm]{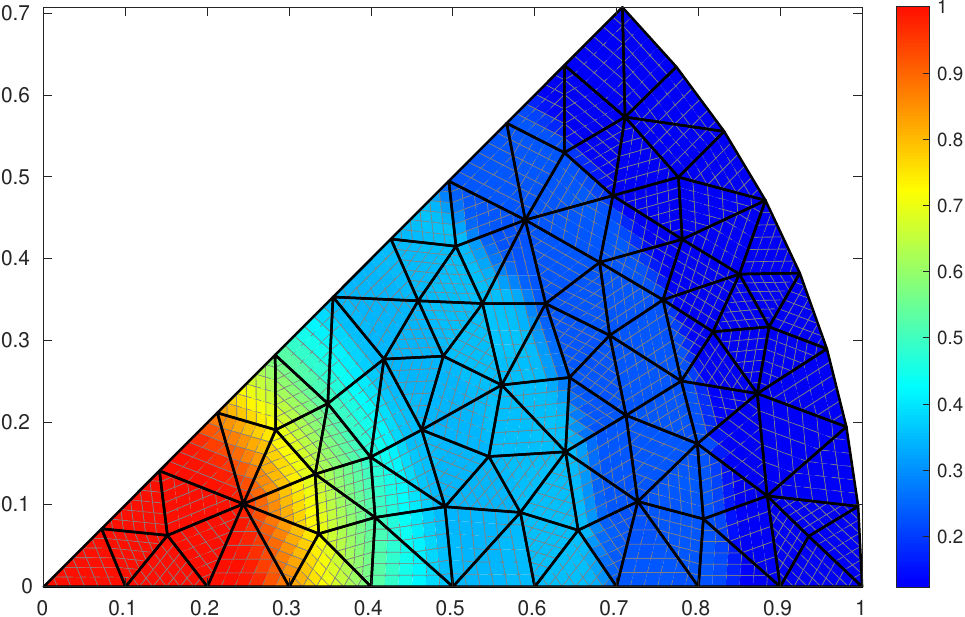}\label{fig_sod_2D_den}}\hspace*{5.mm}
    \subfigure[Density profile]{\includegraphics[height=5.45cm]{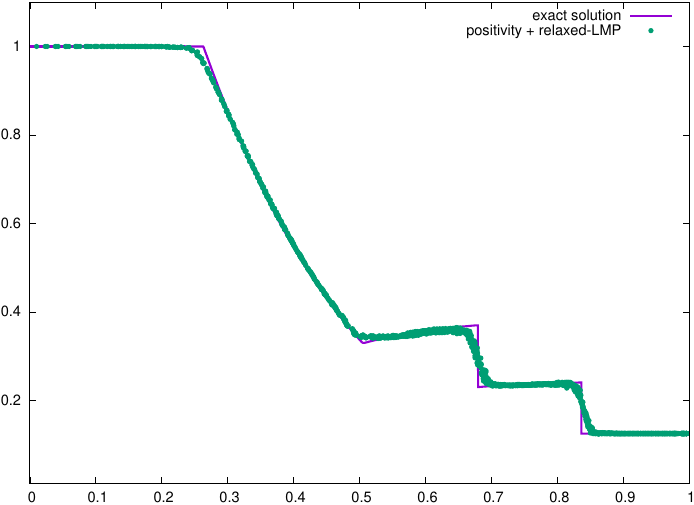}\label{fig_sod_2D_den_profile}}
    \caption{$\poly{5}$-DG/FV scheme with positivity and relaxed-LMP on a $110$ cells mesh}
    \label{fig_sod_2D}
  \end{center}
\end{figure}
In Figure~\ref{fig_sod_2D}, the $\poly{5}$ monolithic DG/FV scheme with positivity and relaxed-LMP has been used on a very coarse anisotropic mesh made of only 110 triangular cells. In the light of Figure~\ref{fig_sod_2D_den}, one can see how the radial wave structure has been accurately capture, even in this extremely coarse mesh case, and how the three types of waves, meaning expansion wave, contact discontinuity and shock wave, travel and go through the large cells. Figure~\ref{fig_sod_2D_den_profile}, where the subcells' mean values versus the subcell centroid radii $\sqrt{x^2+y^2}$ are displayed, confirms this statement as the different points for a given radius do coincide.

\subsubsection{2D Sedov point blast problem}

We consider the Sedov problem for a point-blast in a uniform medium. An exact solution based on self-similarity arguments is available, see for instance \cite{Kamm07}. The initial conditions are characterized by $(\rho_0,p_0,\bs{v}_0)=(1,10^{-14},\bs{0})$, and the polytropic index is equal to $\frac{7}{5}$. We set an initial delta-function energy source at the origin prescribing the pressure in a control volume, yet to be defined, containing the origin as follows, $p_{or}=(\gamma-1)\, \frac{\varepsilon_0}{v_{or}}$, where $v_{or}$ denotes the measure of the chosen control volume and $\varepsilon_0$ the total amount of release energy. By choosing $\varepsilon_0=0.244816$, as suggested in \cite{Kamm07}, the solution consists of a diverging infinite strength shock wave whose front is located at radius $r=1$ at $t=1$, with a peak density reaching 6. The computational domain is defined in polar coordinates by $(r,\theta)\in [0,\, 1.2]\times [0,\frac{\pi}{4}]$. Similarly to the polar Sod shock tube problem, we prescribe symmetry boundary conditions at the boundaries $\theta=0$ and $\theta=\frac{\pi}{4}$, and an outflow condition at $r=1.2$.
Regarding the control volume in which the delta-function energy will be dropped off, generally the cell containing the origin is considered. Here, similarly to \cite{vilar_aplsc_2D} and to make this test case even more challenging, we choose to restrict the energy source only to the one subcell containing the origin. This means that initially, in one grid element the pressure in one subcell will be set to $p_{or}$, while in the remainder of the cell the pressure will be $10^{-14}$. Let us further emphasize that generally in this test case, because one cannot simulate vacuum, the initial pressure is set to $10^{-6}$ over the domain, except at the origin. Here, to make it once again more challenging, we set the initial pressure to $10^{-14}$.
\begin{figure}[!ht]
  \begin{center}
    \subfigure[Energy map]{\includegraphics[height=5.4cm]{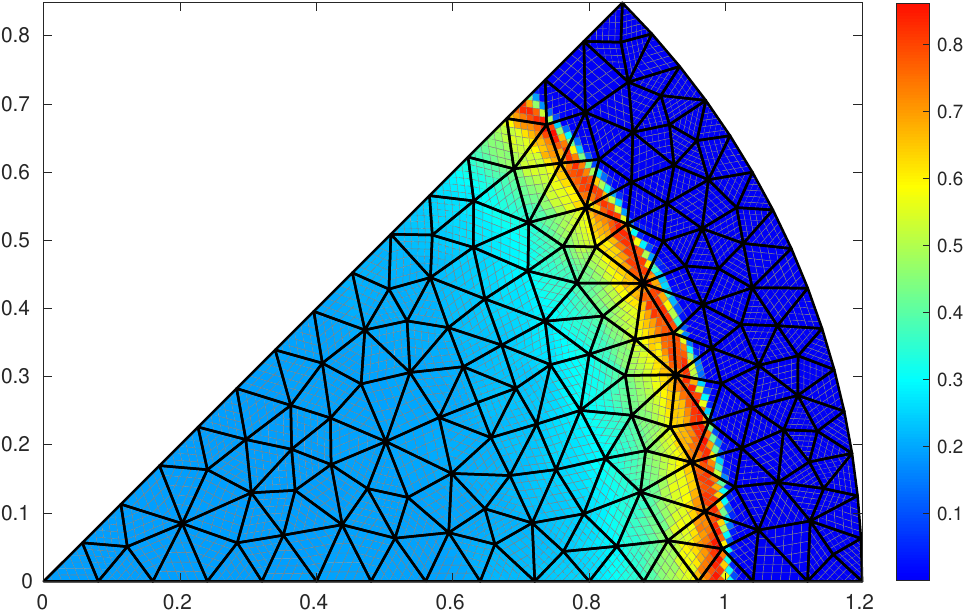}\label{fig_sedov_2D_order6_ene}}\hspace*{5.mm}
    \subfigure[Density profile]{\includegraphics[height=5.4cm]{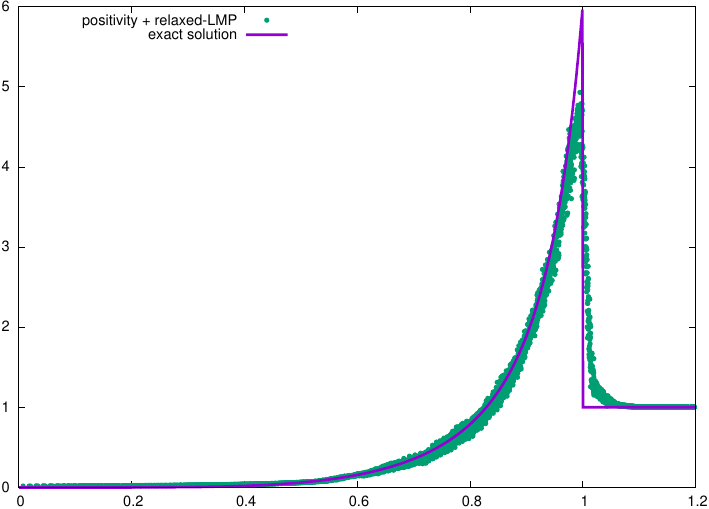}\label{fig_sedov_2D_order6_den_profile}}
    \caption{$\poly{5}$-DG/FV scheme with positivity and relaxed-LMP on a $271$ cells mesh}
    \label{fig_sedov_2D_order6}
  \end{center}
\end{figure}
We run this modified Sedov point blast problem with the $\poly{5}$ monolithic DG/FV scheme, with positivity and relaxed-LMP conditions, on a very coarse grid made of $271$ cells. In this particular case, the amount of total energy contained in the subcell located at the origin reaches $1947.5$, while in the rest of the cell as well as in the remainder of the domain the total energy is set to 2.5E-14. Any scheme lacking positivity-preserving property or a rigorous stabilization technique would fail solving this test problem. In Figure~\ref{fig_sedov_2D_order6}, one can see how the circular aspect of the shock has been accurately captured by the scheme and the shock wave front is correctly located. This latter further goes inside and through different cells, enlightening the very robust and precise subcell resolution of this \scheme method. The numerical solution produced remains quite close to the one-dimensional self-similar exact solution, see Figure~\ref{fig_sedov_2D_order6_den_profile}.\\

\begin{figure}[!ht]
  \begin{center}
    \subfigure[Cells' mean total energy map]{\includegraphics[height=5.5cm]{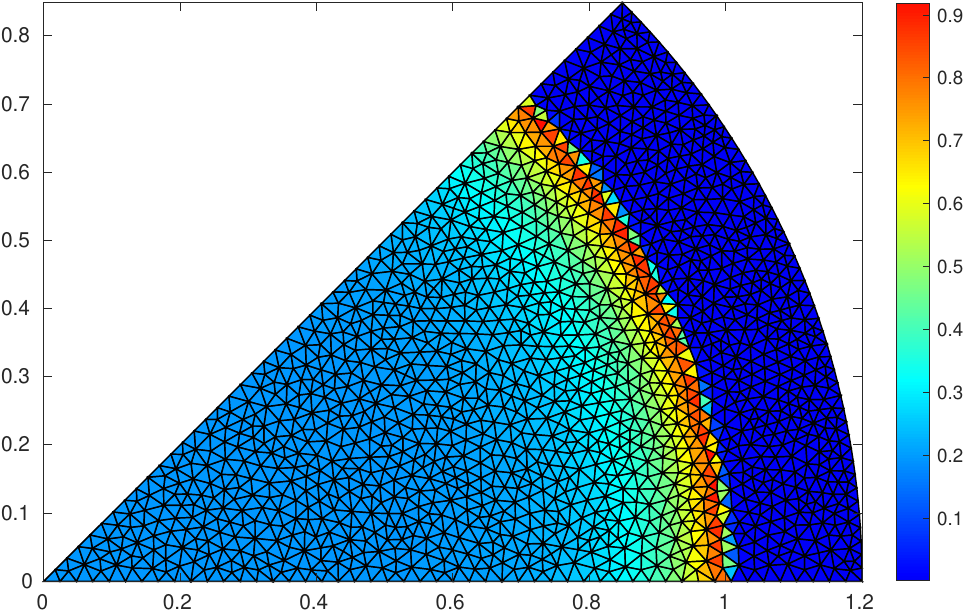}\label{fig_sedov_3rd_ene}}\hspace*{2mm}
    \subfigure[Subcells' mean density versus radii]{\includegraphics[height=5.5cm]{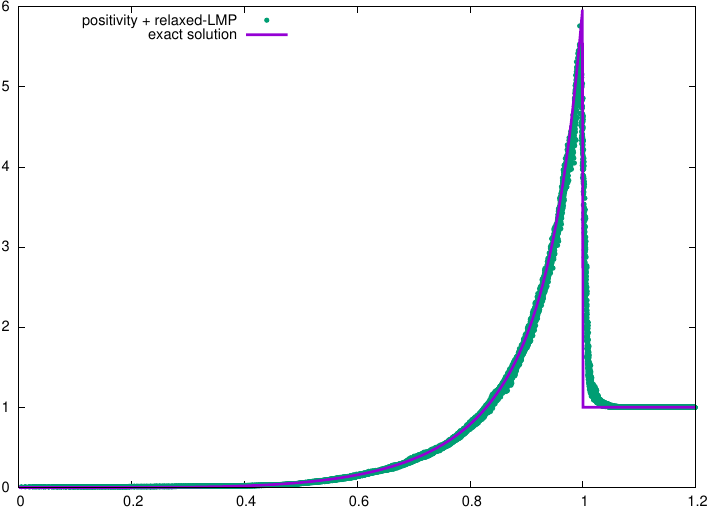}\label{fig_sedov_3rd_den_profile}}
    \caption{$\poly{2}$-DG/FV scheme with positivity and relaxed-LMP on $2894$ cells}
  \label{fig_sedov_3rd-order}
  \end{center}
\end{figure}

Once again, this \scheme scheme performs also very well for lower order methods, as depicted in Figure~\ref{fig_sedov_3rd-order} where a finer grid made of 2894 cells has been used. Indeed, the $\poly{2}$ monolithic scheme solution is very close to the one-dimensional analytical solution. In Figure~\ref{fig_sedov_3rd_ene}, only the cell total energy means values are represented, and not the submean values as we generally do, for a better readability of the results in this finer grid context.

\subsubsection{2D forward-facing step problem}

\begin{figure}[!ht]
  \begin{center}
    \subfigure[$\poly{1}$ on 84108 cells]{\includegraphics[height=5.cm]{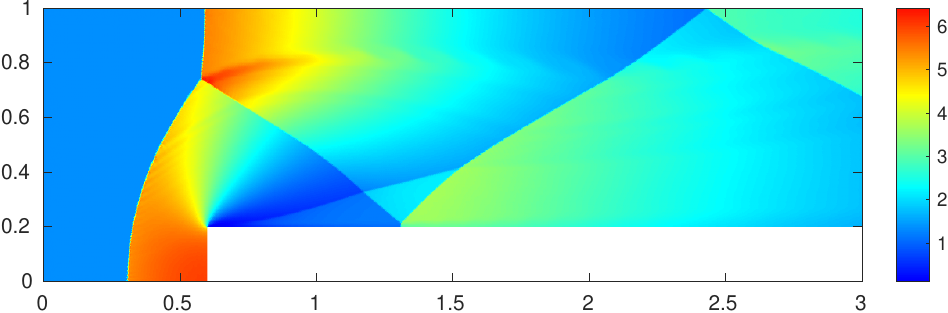}\label{fig_step_P1}}\\[0.5mm]
    \subfigure[$\poly{3}$ on 680 cells]{\includegraphics[height=5.cm]{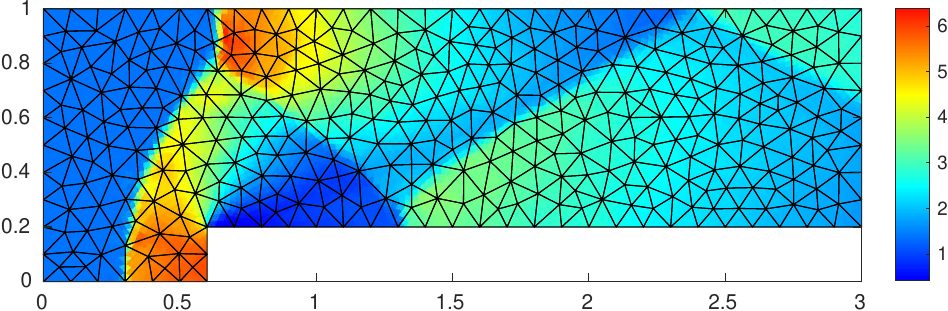}\label{fig_step_P3}}
    \caption{\Scheme scheme: subcells' density mean values}
    \label{fig_step}
  \end{center}
\end{figure}

We now consider the forward facing step problem, which has been initially introduced by A. Emery in \cite{emery_1968}, and further studied by P. Woodward and P. Colella in \cite{woodward_1984}. This challenging test case consists in a Mach 3 flow in a 3 units in length and 1 unit in width wind tunnel. Initially, the tunnel is filled with a gamma gas law with $\gamma=\frac{7}{5}$, which everywhere has density $\rho_0=1.4$, pressure $p_0=1$ and velocity $\bs{v}_0=\(3,\, 0\)\tra$. The 0.2 high step being located at $x=0.6$, the computational domain is then $\([0,\,3]\times[0,\,1]\)\setminus\([0.6,\,3]\times[0.2,\,1]\)$. Gas with this density, pressure and velocity is continually fed in from the left-hand boundary. Let us emphasize that unlike as it is generally done, we did not refine the mesh near the corner, see Figure~\ref{fig_step_P3} for instance, nor modify in any way our monolithic DG/FV scheme.

In Figure~\ref{fig_step}, we compare the numerical solutions obtained by means of our \scheme scheme, ensuring positivity and a relaxed-LMP, respectively with a $\poly{1}$ on a fine grid made of 84108 cells and $\poly{3}$ on coarse grid made of 680 cells. Firstly, let us note that the $\poly{1}$ monolithic DG/FV scheme has produced a quite satisfactory solution close to the expected one and has been able to capture the Kelvin-Helmholtz instabilities in the top of the channel. Secondly, let us point out that due to the complexity of the flow, with multiple shocks waves and walls interacting, the benefit of high-order schemes over low order schemes may be limited. However, as depicted in Figure~\ref{fig_step_P3}, high-order schemes on coarse grids may be a good solution to obtain the main features of the solution under investigation. Indeed, despite the coarseness of the mesh used in Figure~\ref{fig_step_P3}, the shocks and the rarefaction fan created around the corner are quite well resolved, while ensuring a low oscillatory behavior. Now, to capture the finer structures of the solution, a finer mesh has to be used, keeping in mind that high-order monolithic DG/FV scheme will always outclass the lower-order ones, but being obviously more computationally costly.

\subsubsection{2D hypersonic flow over half cylinder problem}
\label{test_Mach_20}

\begin{figure}[!ht]
  \begin{center}
    \subfigure[Grid used]{\includegraphics[height=12.6cm]{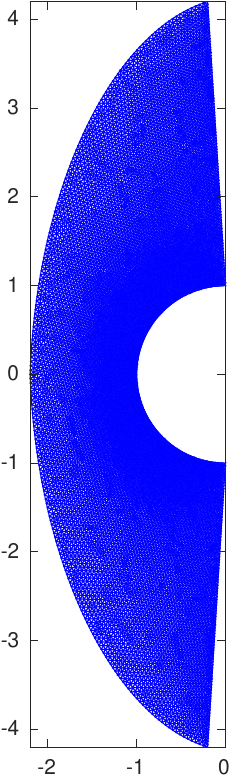}\label{fig_half_cyl_FV_grid}}\hspace*{1.4cm}
    \subfigure[HLL numerical flux]{\includegraphics[height=12.6cm]{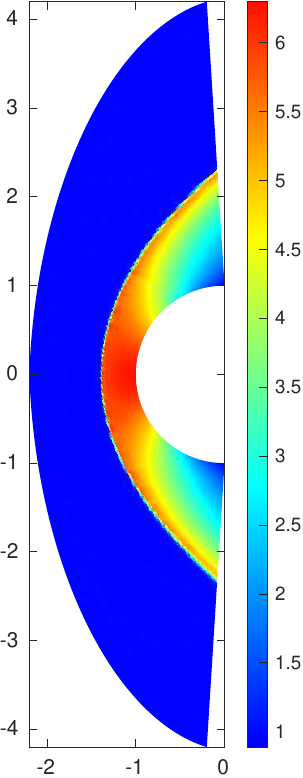}\label{fig_half_cyl_FV_HLL}}\hspace*{0.6cm}
    \subfigure[HLL-C numerical flux]{\includegraphics[height=12.6cm]{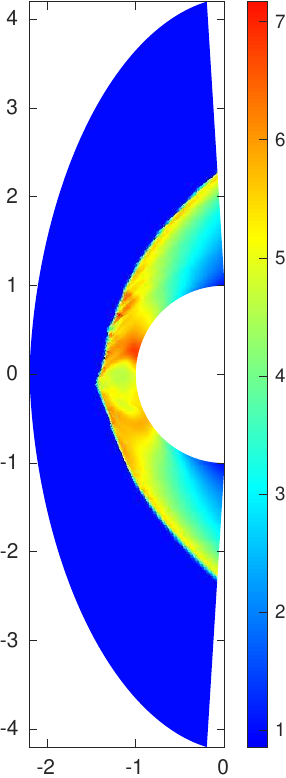}\label{fig_half_cyl_P0_HLLC}}
    \caption{1st-order FV scheme on a grid made of 25266 cells with HLL and HLL-C numerical fluxes}
    \label{fig_half_cyl_FV}
  \end{center}
\end{figure}

The hypersonic flow over a half-cylinder test case is a well-documented test case used to challenge numerical methods. Specifically, some schemes may develop the infamous carbuncle phenomenon, even using classical FV scheme. Instead of producing a smooth bow shock profile upstream of the half-cylinder, the carbuncle issue manifests as a pair of oblique shocks ahead of the stagnation region, compromising the overall flow predictions around the cylinder. Following the approach in \cite{Rodionov_half_cyl}, we simulate an inviscid flow at Mach $M_a=20$ around a half-cylinder blunt body subjected to an incoming hypersonic flow characterized by $\(\rho_i,\,p_i,\,\bs{v}_i\)=\(1,\,1,\,(M_a\,\sqrt{\gamma},\,0)\tra\)$ with $\gamma=\frac{7}{5}$. The steady-state resulting flow is simulated using an explicit time-marching procedure, ending at time $t=2.5$. The computational domain is sufficiently large, containing half of a cylinder centered at the origin with a radius $r=1$ and a left incoming hypersonic flow. At the cylinder surface, a wall-type boundary condition is applied, while the bottom and upper boundary conditions are free outflow and an inflow condition is applied at the left boundary. First, to exhibit the so-called carbuncle effect, we display in Figure~\ref{fig_half_cyl_FV} the numerical solutions obtained using a first-order FV scheme, based respectively on HLL and HLL-C numerical fluxes, on a fine grid made of 25266 cells, see Figure~\ref{fig_half_cyl_FV_grid}.

\begin{figure}[!ht]
  \begin{center}
    \subfigure[HLL-C/HLL-C]{\includegraphics[height=14.3cm]{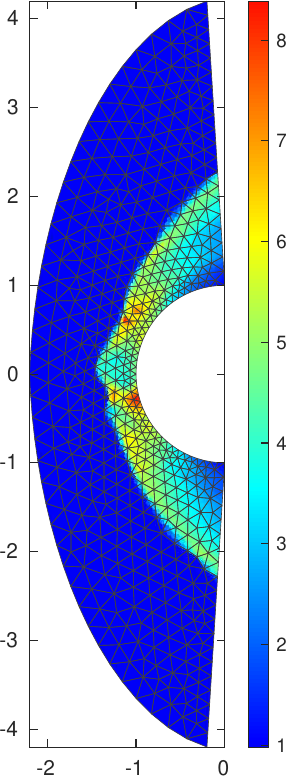}\label{fig_half_cyl_P2_HLLC}}\hspace*{6.mm}
    \subfigure[HLL-C/HLL]{\includegraphics[height=14.3cm]{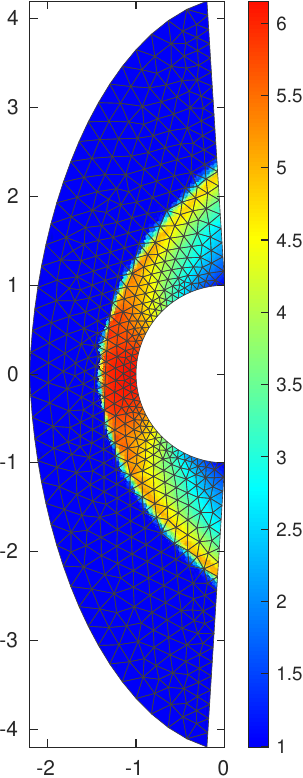}\label{fig_half_cyl_P2_HLL_HLLC}}
    \caption{$\poly{2}$-DG/FV scheme with positivity and relaxed-LMP on $1044$ cells: density}
    \label{fig_half_cyl_P2}
  \end{center}
\end{figure}

As expected, the use of HLL-C numerical flux does trigger the carbuncle effect, see Figure~\ref{fig_half_cyl_P0_HLLC}, while the use of HLL numerical flux does not, see Figure~\ref{fig_half_cyl_FV_HLL}. Now, in addition to assess the robustness of our monolithic scheme in this Mach 20 hypersonic flow context, we want to assess how the carbuncle effect translates going to higher orders of accuracy. To do so, we run our $\poly{2}$ monolithic DG/FV scheme, ensuring positivity and a relaxed-LMP, on a coarse mesh made of 1044 cells, and display the solutions in Figure~\ref{fig_half_cyl_P2}.

\begin{figure}[!ht]
  \begin{center}
    \subfigure[HLL-C/HLL: density]{\includegraphics[height=14.3cm]{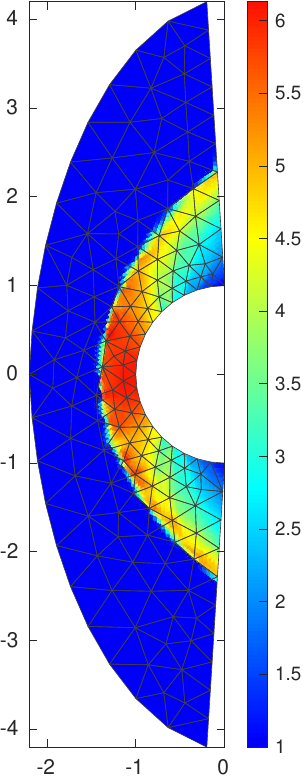}\label{fig_half_cyl_P5_HLL_HLLC}}\hspace*{6.mm}
    \subfigure[HLL-C/HLL: blending coeff.]{\includegraphics[height=14.3cm]{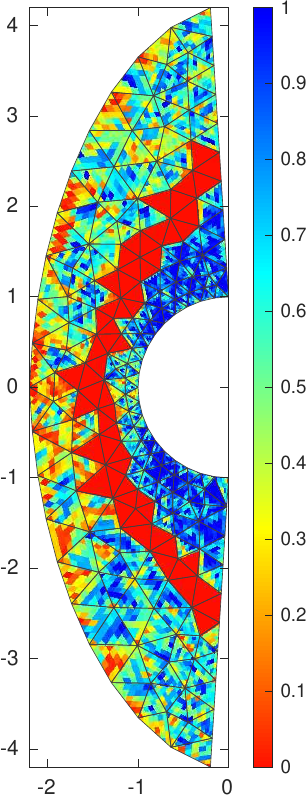}\label{fig_half_cyl_P5_coeff}}
    \caption{$\poly{5}$-DG/FV scheme with positivity and relaxed-LMP on $292$ cells}
    \label{fig_half_cyl_P5}
  \end{center}
\end{figure}

Let us first emphasize that it is absolutely not required to use the same numerical flux to compute the DG residual in \eqref{recons_flux_definition} and for the first-order FV flues in \eqref{blended_flux}. Thus, to see the repercussion of such choices on the potential trigger of the carbuncle effect, we compare in Figure~\ref{fig_half_cyl_P2} the solutions obtained by means of our $\poly{2}$ monolithic DG/FV scheme respectively using HLL-C numerical flux for both DG (hence the reconstructed fluxes) and the first-order FV fluxes in the first case, Figure~\ref{fig_half_cyl_P2_HLLC}, and HLL-C for DG and HLL for the FV fluxes in the second case, Figure~\ref{fig_half_cyl_P2_HLL_HLLC}. In Figure~\ref{fig_half_cyl_P2_HLLC}, one can see that even going to high order of accuracy, the use of HLL-C numerical flux does trigger the carbuncle effect. However, if in the monolithic scheme, the first-order numerical method forming the safe base of this blended scheme uses HLL numerical flux, the approximated solution will not present any carbuncle effect. This highlights the fact that, with appropriate choices for the DG and FV numerical fluxes, it should be possible to ensure some properties on the high-order parts of solution and others on the low-order parts of the solution. Finally, to test once more the high capability and robustness of this monolithic scheme going to very high-order of accuracy and very coarse meshes, even in the simulation of this complex high Mach hypersonic flow, we show in Figure~\ref{fig_half_cyl_P5} the solution obtained by means of the $\poly{5}$ monolithic DG/FV scheme, ensuring positivity and a relaxed-LMP, on an extremely coarse mesh made of only 292 cells.

In the light of Figure~\ref{fig_half_cyl_P5_HLL_HLLC}, one can see once more how robust and accurate the monolithic scheme is, taking into account the extreme coarseness of the grid used, as the scheme has been able to capture the bow structure of the shock. Furthermore, the blending coefficients displayed on Figure~\ref{fig_half_cyl_P5_coeff} allow us to illustrate what has been said in Remark~\ref{remark_NaN}. Indeed, one can see that, in cells containing the bow shock, every subcells are computed through a first-order FV scheme, as the $\theta_m^{\,c}$ are equal to zero everywhere in the cell. This shows that the reconstructed fluxes may have developed non-admissible values, revealing that DG scheme has totally failed in those cells. This generally comes from the computation of square root of negative values in the evaluation of the DG numerical fluxes. But as long as the first-order scheme is robust, the monolithic and its associated simulation code could not crash.

\section{Conclusion}
\label{conclusion}

This article is concerned with the construction of a new type of monolithic scheme, based on general unstructured grids, blending locally at the subcell level DG scheme and first-order FV scheme. This subcell blending procedure relies on the expression of DG methods as a finite volume scheme on a subgrid. By means of this theoretical part, we combine in a convex manner, at the subcell level, the so-called reconstructed fluxes and first-order FV numerical fluxes, through a blending coefficient $\theta_{mp}$. The next step is then to determine those blending coefficients to achieve all the desired properties. In the first part of this paper, we have first focused our attention of the issue of entropy stability. Different types of entropy stabilities were hence introduced, along with the corresponding blending coefficients. In particular, a semi-discrete cell entropy stability, for a given entropy, has been proposed which allows the preservation of the high-order accuracy of the \scheme scheme. However, the numerical results showed that this entropy criterion might be, in some complex case, too relaxed to capture the unique entropic weak solution. Furthermore, the cost of such entropy condition is very high as the first-order FV numerical fluxes have to be specifically modified, which may lead to the loss of other desired properties as positivity for instance. In a second part, we focus on the imposition of different maximum principles, a global one to ensure the numerical solution to remain in a convex admissible set, and a local one to address the issues of spurious oscillations. A wide number of test cases on different problems have been used to depict the very good performance and robustness of the presented monolithic scheme using those principles.\\

In a very near future, we expect to apply those monolithic schemes to the problematic of coastal flows simulation and their coupling with a moving object. We also intend to extend this \scheme framework to multi-dimensional approximated Riemann solvers, and then to moving grid configurations, both in ALE and Lagrangian formalisms. Finally, we plan as well to extend the present scheme to the three-dimensional case, as the theory developed should remain exactly the same.


\newpage
\appendix

\appendixpage
\section{Some properties on FV schemes and E-fluxes}
\label{sect_E-flux}

This appendix aims at recalling some properties yield by a FV scheme
\begin{align}
  \label{subcell_FV}
  \ov{u}_m^{c,\,n+1}=\ov{u}_m^{c,\,n}-\Frac{\Dt}{|S_m^c|}\,\Sum_{S_{p}^v\,\in\mc{V}_m^{\,c}} l_{mp}\, \mc{F}\(\ov{u}_m^{c,\,n},\ov{u}_p^{v,\,n},\bs{n}_{mp}\),
\end{align}
relying on a general numerical flux \eqref{num_flux}, in the simple case of SCL. All the following properties can be easily extended to systems and other types of fluxes as HLL for instance.

\subsection{Discrete maximum principle}
\label{sect_FV_DMP}

First, let us show that following discrete maximum principle holds
\begin{align}
  \label{DMP}
\min\big(\ov{u}_m^{c,\,n},\,\Min_{S_{p}^v\,\in\mc{V}_m^{\,c}}\ov{u}_p^{v,\,n}\big)\leq \ov{u}_m^{c,\,n+1} \leq  \max\big(\ov{u}_m^{c,\,n},\,\Max_{S_{p}^v\,\in\mc{V}_m^{\,c}}\ov{u}_p^{v,\,n}\big).
\end{align}
Following the steps presented in Section~\ref{sect_monolithic}, scheme \eqref{subcell_FV} can be recast into the following Godunov-type form
\begin{align}
  \label{convex_combo_FV}
  \ov{u}_m^{c,n+1}= \Big(1-\Frac{\Dt}{|S_m^c|}\,\Sum_{S_{p}^v\,\in\mc{V}_m^{\,c}}l_{mp}\,\gamma_{mp}\Big)\,\ov{u}_m^{c,n}+\Frac{\Dt}{|S_m^c|}\,\Sum_{S_{p}^v\,\in\mc{V}_m^{\,c}}l_{mp}\,\gamma_{mp}\,u_{mp}^{\ast,\text{\,\tiny FV}},
\end{align}
where $u_{mp}^{\ast,\text{\,\tiny FV}}:=u^\ast(\ov{u}_m^{c,n},\ov{u}_p^{v,n},\bs{n}_{mp})$, the FV Riemann intermediate state, writes as
\begin{align}
  \label{intermediate_FV}
  u^\ast(u_L,u_R,\bs{n})=\Frac{u_L+u_R}{2}-\Frac{\big(\bs{F}(u_R)-\bs{F}(u_L)\big)\pds\bs{n}}{2\,\gamma(u_L,u_R,\bs{n})}.
\end{align}
A sufficient condition to ensure \eqref{DMP} is then to show that $u^\ast(u_L,u_R,\bs{n})$ lies in $I(u_L,u_R)$. This condition is evident, as
\begin{align*}
  u^\ast(u_L,u_R,\bs{n})=u_L\,\(\Inv{2}+\Frac{\beta}{2\,\gamma}\)+u_R\,\(\Inv{2}-\Frac{\beta}{2\,\gamma}\),
\end{align*}
with $\gamma:=\gamma(u_L,u_R,\bs{n})$ and $\beta=\frac{(\bs{F}(u_R)-\bs{F}(u_L))\,\pds\,\bs{n}}{u_R-u_L}$. Since $\gamma\geq\max_{w\,\in\, I(u_L,u_R)}(|\bs{F}'(w)\pds\bs{n}|)$, it directly follows that $|\beta|\leq \gamma$. The intermediate state $u^\ast$ hence writes as a convex combination of $u_L$ and $u_R$.

\subsection{Discrete entropy stability}
\label{sect_FV_entropy}

\begin{figure}[!ht]
  \begin{center}
    \includegraphics[width=6cm]{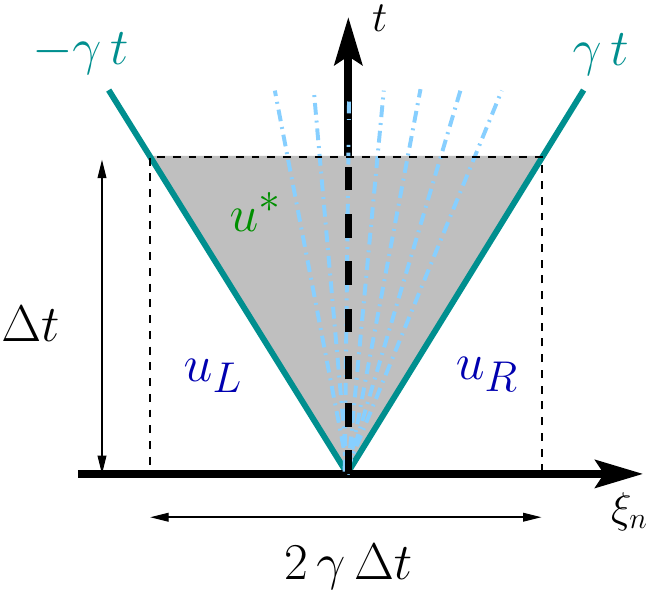}
    \caption{Riemann fan}
  \label{Fig_Riemann_fan}
  \end{center}
\end{figure}

Now, let us recall that the discrete scheme \eqref{subcell_FV} is entropy stable for any entropy, as it implies the following inequality
\begin{align}
  \label{subcell_FV_entropy}
  \eta\big(\ov{u}_m^{c,\,n+1}\big)\leq \eta\big(\ov{u}_m^{c,\,n}\big)-\Frac{\Dt}{|S_m^c|}\,\Sum_{S_{p}^v\,\in\mc{V}_m^{\,c}} l_{mp}\, \phi^\ast\(\ov{u}_m^{c,\,n},\ov{u}_p^{v,\,n},\bs{n}_{mp}\),
\end{align}
with function $\phi^\ast$, a consistent numerical entropy flux, such that $\phi^\ast(u,u,\bs{n})=\bs{\phi}(u)\pds\bs{n}$. To be coherent with the numerical flux $\mc{F}$ under consideration \eqref{num_flux}, we define the numerical entropy flux as
\begin{align}
  \label{num_entr_flux}
  \phi^\ast(u_L,u_R,\bs{n})=\Frac{\big(\bs{\phi}(u_L)+\bs{\phi}(u_R)\big)\pds\bs{n}}{2}-\Frac{\gamma(u_L,u_R,\bs{n})}{2}\,\big(\eta(u_R)-\eta(u_L)\big).
\end{align}
To demonstrate \eqref{subcell_FV_entropy}, let us first see that, by convexity of the entropy, the convex relation \eqref{convex_combo_FV} leads to
\begin{align}
  \label{convex_combo_FV_entropy}
  \eta\(\ov{u}_m^{c,n+1}\)\leq \eta\(\ov{u}_m^{c,n}\)-\Frac{\Dt}{|S_m^c|}\,\Sum_{S_{p}^v\,\in\mc{V}_m^{\,c}}l_{mp}\,\gamma_{mp}\,\(\eta\(\ov{u}_m^{c,n}\)-\eta\(u_{mp}^{\ast,\text{\,\tiny FV}}\)\).
\end{align}
Thus, it directly follows that a sufficient condition to obtain \eqref{subcell_FV_entropy} is
\begin{align*}
  -\,\gamma_{mp}\,\Big(\eta\big(\ov{u}_m^{c,n}\big)-\eta\big(u_{mp}^{\ast,\text{\,\tiny FV}}\big)\Big)\leq -\Big(\phi^\ast\big(\ov{u}_m^{c,\,n},\ov{u}_p^{v,\,n},\bs{n}_{mp}\big)-\bs{\phi}\big(\ov{u}_m^{c,n}\big)\pds\bs{n}_{mp}\Big).
\end{align*}
By means of the numerical entropy flux \eqref{num_entr_flux}, this sufficient condition can be reformulated as follows
\begin{align}
  \label{CS_FV_entropy}
  \eta\big(u^\ast(u_L,u_R,\bs{n})\big)\leq \Frac{\eta(u_L)+\eta(u_R)}{2}-\Frac{\big(\bs{\phi}(u_R)-\bs{\phi}(u_L)\big)\pds\bs{n}}{2\,\gamma(u_L,u_R,\bs{n})}.
\end{align}
To ensure that condition \eqref{CS_FV_entropy} is indeed guaranteed, let us introduce the following Riemann problem
\begin{subequations}
\label{Riemann_pb}
\begin{alignat}{2}[left = \empheqlbrace\,]
&\vdt{\,u}(\bs{x},t)+\divx{\bs{F}\(u(\bs{x},t)\)}=0, \qqquad& (\bs{x},t)\in\,\R^d\times\R^+, \label{Riemann1}\\
  &u(\bs{x},0)=\left\{\begin{array}{ll}u_L \quad &\text{if }\; (\bs{x}\pds\bs{n})<0,\\u_R \quad &\text{if }\; (\bs{x}\pds\bs{n})>0.
  \end{array}\right. \label{Riemann2}
\end{alignat}
\end{subequations}
Introducing rotated space variables $\xi_n=(\bs{x}\pds\bs{n})$ and $\xi_\tau=(\bs{x}\pds\bs{\tau})$, with $\bs{n}$ a given unit normal and $\bs{\tau}=\bs{n}^\perp$, and due to the invariance by rotation of SCL, equation \eqref{Riemann1} rewrites as
\begin{align*}
  \vdt{\,u}+\partial_{\xi_n}\(\bs{F}(u)\pds\bs{n}\)+\partial_{\xi_\tau}\(\bs{F}(u)\pds\bs{\tau}\)=0.
\end{align*}
Finally, because $u_0$ does only depends on $\xi_n$, the solution $u$ does as well, and the PDE reduces to the following 1D problem
\begin{align}
\label{Riemann_pb_normal}
  \vdt{\,u}+\partial_{\xi_n}F_n(u)=0,
\end{align}
where $F_n(u)=\(\bs{F}(u)\pds\bs{n}\)$. Now, let us show that $u^\ast(u_L,u_R,\bs{n})$, defined in \eqref{intermediate_FV}, is nothing but the average value of $\mc{W}\big(\frac{\xi_n}{t};\,u_L,\,u_R\big)$, the unique entropic weak solution of the considered Riemann problem, as long as $\gamma:=\gamma(u_L,u_R,\bs{n})\geq\max_{w\,\in\, I(u_L,u_R)}(|\bs{F}'(w)\pds\bs{n}|)$. In this case, the waves produced by the initial discontinuity will remain left and right bounded by respectively $\xi_n=\pm\,\gamma\,t$, as depicted by Figure~\ref{Fig_Riemann_fan}.
First, integrating equation \eqref{Riemann_pb_normal} onto the time-space box $[-\gamma\,t,\gamma\,t]\times[0,\Dt]$, displayed in Figure~\ref{Fig_Riemann_fan}, one gets
\begin{align}
  \label{intermediate_FV_entropy1}
  \Inv{2\,\gamma\,\Dt}\Int_{-\gamma\,t}^{\gamma\,t} \mc{W}\Big(\frac{\xi_n}{\Dt};\,u_L,\,u_R\Big)\,\dd\, \xi_n= \Frac{u_L+u_R}{2}-\Frac{\big(\bs{F}(u_R)-\bs{F}(u_L)\big)\pds\bs{n}}{2\,\gamma}:=u^\ast.
\end{align}
Because $\mc{W}\big(\frac{\xi_n}{t};\,u_L,\,u_R\big)$, the unique entropic weak solution, lies in $I(u_L,u_R)$, it furthers demonstrates that $u^\ast\in I(u_L,u_R)$ as well. Now, let us recall that the unique solution $u(\bs{x},t):=\mc{W}\big(\frac{\xi_n}{t};\,u_L,\,u_R\big)$ ensures, in a weak sens, the following entropic inequalities
\begin{align}
\label{Riemann_pb_normal_entropy}
  \vdt{\,\eta(u)}+\partial_{\xi_n}\phi_n(u)\leq 0,
\end{align}
for any couple entropy - entropy flux $(\eta,\bs{\phi})$, where $\phi_n(u)=\(\bs{\phi}(u)\pds\bs{n}\)$. Similarly as before, integrating \eqref{Riemann_pb_normal_entropy} onto $[-\gamma\,t,\gamma\,t]\times[0,\Dt]$, it follows that
\begin{align}
  \label{intermediate_FV_entropy2}
  \Inv{2\,\gamma\,\Dt}\Int_{-\gamma\,t}^{\gamma\,t} \eta\(\mc{W}\Big(\frac{\xi_n}{\Dt};\,u_L,\,u_R\Big)\)\,\dd\, \xi_n\leq \Frac{\eta(u_L)+\eta(u_R)}{2}-\Frac{\big(\bs{\phi}(u_R)-\bs{\phi}(u_L)\big)\pds\bs{n}}{2\,\gamma}.
\end{align}
Finally, due to entropy convexity and using Jensen's inequality
\begin{align*}
  \eta\(u^\ast\):=\eta\(\Inv{2\,\gamma\,\Dt}\Int_{-\gamma\,t}^{\gamma\,t} \mc{W}\Big(\frac{\xi_n}{\Dt};\,u_L,\,u_R\Big)\,\dd\, \xi_n\)\leq\Inv{2\,\gamma\,\Dt}\Int_{-\gamma\,t}^{\gamma\,t} \eta\(\mc{W}\Big(\frac{\xi_n}{\Dt};\,u_L,\,u_R\Big)\)\,\dd\, \xi_n,
\end{align*}
relation \eqref{intermediate_FV_entropy2} reduces to the desired condition \eqref{CS_FV_entropy}.

\subsection{Two-point Tadmor relation}
\label{sect_FV_Tadmor}

Here, we show that the use of a numerical flux \eqref{num_flux} with $\gamma(u_L,u_R,\bs{n})\geq\max_{w\,\in\, I(u_L,u_R)}(|\bs{F}'(w)\pds\bs{n}|)$, for which we have just displayed how it ensures a discrete entropy inequality for any entropy, also guarantees the Tadmor inequality \eqref{entropy_FV_Tadmor}. Such inequality ensures the semi-discrete FV scheme to be entropy stable for a given entropy, \cite{Tadmor_entropy,Tadmor_entropy_acta}. To this end, we will exhibit how a numerical flux \eqref{num_flux} can be put into the following form, with $D\geq0 $
\begin{align}
  \label{num_flux_entropy}
  \mc{F}\(u_L,u_R,\bs{n}\)=\Frac{\bs{\Psi}(v_R)-\bs{\Psi}(v_L)}{v_R-v_L}\pds\bs{n}-\Frac{D}{2}\,(v_R-v_L),
\end{align}
where $v_{L/R}=v(u_{L/R})$. If the entropy viscosity dissipation coefficient $D=0$, the semi-discrete FV scheme will be entropy conservative, while being entropy dissipative for $D>0$. The combination of definitions \eqref{num_flux} and \eqref{num_flux_entropy} states that the entropy dissipation coefficient is given by
\begin{align*}
  D&=\gamma\,\(\Frac{u_R-u_L}{v_R-v_L}\)+\Frac{2}{(v_R-v_L)^2}\(\bs{\psi}_R-\bs{\psi}_L-\Frac{(\bs{F}_L+\bs{F}_R)}{2}(v_R-v_L)\)\pds\bs{n},\\[3mm]
  &= \gamma\,\(\Frac{u_R-u_L}{v_R-v_L}\)+\Frac{2}{(v_R-v_L)^2}\(\Int_{v_L}^{v_R} \bs{\Psi}'(v)\;\dd\, v-\Frac{(\bs{F}_R+\bs{F}_L)}{2}\Int_{u_L}^{u_R} v'(u)\;\dd\, u\)\pds\bs{n},
\end{align*}
where $\gamma:=\gamma(u_L,u_R,\bs{n})$, while $\bs{\psi}_{L/R}=\bs{\psi}(u_{L/R})=\bs{\Psi}(v_{L/R})$ and $\bs{F}_{L/R}=\bs{F}(u_{L/R})$. By definition of the entropy potential flux, we have that $\bs{\Psi}'(v(u))=\bs{F}(v(u))$. By a change of variable in the first integral, it follows that
\begin{align*}
  D&= \gamma\,\(\Frac{u_R-u_L}{v_R-v_L}\)+\Frac{2}{(v_R-v_L)^2}\(\Int_{u_L}^{u_R} \bs{F}(u)\,v'(u)\;\dd\, u-\Frac{(\bs{F}_R+\bs{F}_L)}{2}\Int_{u_L}^{u_R} v'(u)\;\dd\, u\)\pds\bs{n},\\[3mm]
  &= \gamma\,\(\Frac{u_R-u_L}{v_R-v_L}\)-\Inv{(v_R-v_L)^2}\,\Int_{u_L}^{u_R} \bs{n}\pds\big(\bs{F}_R+\bs{F}_L-2\,\bs{F}(u)\big)\,v'(u)\;\dd\, u,\\[3mm]
  &= \(\Frac{u_R-u_L}{v_R-v_L}\)^2\,\Inv{u_R-u_L}\,\Int_{u_L}^{u_R} \underbrace{\(\gamma-\Frac{\big(\bs{F}_R+\bs{F}_L-2\,\bs{F}(u)\big)\pds\bs{n}}{u_R-u_L}\)}_{\Gamma(u)}\,v'(u)\;\dd\, u.
\end{align*}
Finally, thanks to the entropy convexity, $v'(u)>0$, if $\forall\, u\in I(u_L,u_R), \Gamma(u)\geq 0$ then $D\geq 0$. This is actually the case because
\begin{align*}
  \Gamma(u)&\geq \gamma-\left|\Frac{\big(\bs{F}_R+\bs{F}_L-2\,\bs{F}(u)\big)\pds\bs{n}}{u_R-u_L}\right|\geq \gamma-\Frac{\big|\big(\bs{F}_R-\bs{F}(u)\big)\pds\bs{n}\big|+\big|\big(\bs{F}_L-\bs{F}(u)\big)\pds\bs{n}\big|}{|u_R-u_L|},\\[3mm]
 &= \gamma-\(\left|\Frac{u-u_L}{u_R-u_L}\right|\,\left|\Frac{\big(\bs{F}(u)-\bs{F}_L\big)\pds\bs{n}}{u-u_L}\right|+\left|\Frac{u_R-u}{u_R-u_L}\right|\,\left|\Frac{\big(\bs{F}_R-\bs{F}(u)\big)\pds\bs{n}}{u_R-u}\right|\),\\[3mm]
 &\geq \gamma-\(\Frac{\left|u-u_L\right|+\left|u_R-u\right|}{\left|u_R-u_L\right|}\)\max_{w\,\in\, I(u_L,u_R)}(|\bs{F}'(w)\pds\bs{n}|),\\[3mm]
 &= \gamma-\max_{w\,\in\, I(u_L,u_R)}(|\bs{F}'(w)\pds\bs{n}|).
\end{align*}
Under the condition that $\gamma\geq\max_{w\,\in\, I(u_L,u_R)}(|\bs{F}'(w)\pds\bs{n}|)$, the numerical flux \eqref{num_flux} can indeed be expressed as in \eqref{num_flux_entropy} with an entropy dissipation coefficient $D\geq 0$.



\begin{thebibliography}{}
  
  \bibitem{batten}
P.~Batten, N.~Clarke, C.~Lambert, and Causon.
\newblock {On the choice of wavespeeds for the HLLC Riemann solver}.
\newblock {\em SIAM J. Sci. Comput.}, 18:1553--1570, 1997.

\bibitem{bisw}
R.~Biswas, K.~Devine, and J.E. Flaherty.
\newblock Parallel adaptive finite element methods for conservation laws.
\newblock {\em Applied Numerical Mathematics}, 14:255--284, 1994.

\bibitem{Boris_FCT}
J.-P. Boris and D.-L. Book.
\newblock {Flux-corrected transport. I. SHASTA, a fluid transport algorithm
  that works}.
\newblock {\em J. Comp. Phys.}, 11(1):38--69, 1973.

\bibitem{burb01}
A.~Burbeau, P.~Sagaut, and C.-H. Bruneau.
\newblock {A problem-independent limiter for high-order Runge Kutta
  discontinuous Galerkin methods}.
\newblock {\em J. Comp. Phys.}, 169:111--150, 2001.

\bibitem{IDP_renac}
V.~Carlier and F.~Renac.
\newblock Invariant domain preserving high-order spectral discontinuous
  approximations of hyperbolic systems.
\newblock {\em SIAM J. Sci. Comput.}, 45(3):A1385--A1412, 2023.

\bibitem{entropy_carpenter}
M.~H. Carpenter, T.~Fisher, E.~Nielsen, and S.~Frankel.
\newblock {Entropy Stable Spectral Collocation Schemes for the Navier--Stokes
  Equations: Discontinuous Interfaces}.
\newblock {\em SIAM J. Sci. Comput.}, 36:B835--B867, 2014.

\bibitem{DG_SBP_chan}
J.~Chan.
\newblock On discretely entropy conservative and entropy stable discontinuous
  galerkin methods.
\newblock {\em J. Comp. Phys.}, 362:346--374, 2018.

\bibitem{SBP_shu}
T.~Chen and C.-W. Shu.
\newblock {Entropy stable high order discontinuous Galerkin methods with
  suitable quadrature rules for hyperbolic conservation laws}.
\newblock {\em J. Comp. Phys.}, 345:427--461, 2017.

\bibitem{mood1}
S.~Clain, S.~Diot, and R.~Loub\`{e}re.
\newblock {A high-order finite volume method for hyperbolic systems:
  Multi-dimensional Optimal Order Detection (MOOD)}.
\newblock {\em J. Comp. Phys.}, 230:4028--4050, 2011.

\bibitem{Cockburn_lcs4}
B.~Cockburn, S.~Hou, and C.-W. Shu.
\newblock {The Runge-Kutta local projection discontinuous Galerkin finite
  element method for conservation laws IV: The multidimensional case}.
\newblock {\em Math. Comp.}, 54:545--581, 1990.

\bibitem{Cockburn_lcs5}
B.~Cockburn, S.~Hou, and C.-W. Shu.
\newblock {The Runge-Kutta Discontinuous Galerkin Method for Conservation Laws
  V: Multidimensional Systems}.
\newblock {\em J. Comp. Phys.}, 141:199--224, 1998.

\bibitem{Cockburn_lcs2}
B.~Cockburn and C.-W. Shu.
\newblock {TVB Runge-Kutta local projection discontinuous Galerkin finite
  element method for conservation laws II: General framework}.
\newblock {\em Math. Comp.}, 52:411--435, 1989.

\bibitem{dumbser_subFV_lim_tri}
M.~Dumbser and R.~Loub\`{e}re.
\newblock {A simple robust and accurate a posteriorisub-cell finite volume
  limiter for the discontinuousGalerkin method on unstructured meshes}.
\newblock {\em J. Comp. Phys.}, 319:163--199, 2016.

\bibitem{emery_1968}
A.~Emery.
\newblock {An evaluation of several differencing methods for inviscid flow
  problems}.
\newblock {\em J. Comp. Phys.}, 2(3):306--331, 1968.

\bibitem{DG_SBP_renac}
{F. Renac}.
\newblock Entropy stable, robust and high-order dgsem for the compressible
  multicomponent euler equations.
\newblock {\em J. Comp. Phys.}, 445:110584, 2021.

\bibitem{SBP_carpenter}
T.~C. Fisher and M.~H. Carpenter.
\newblock {High-order entropy stable finite difference schemes for nonlinear
  conservation laws: Finite domains}.
\newblock {\em J. Comp. Phys.}, 252:518--557, 2013.

\bibitem{DG_SBP_gassner}
{G. Gassner}.
\newblock {A Skew-Symmetric Discontinuous Galerkin Spectral Element
  Discretization and Its Relation to SBP-SAT Finite Difference Methods}.
\newblock {\em SIAM J. Sci. Comput.}, 35:1233--1253, 2013.

\bibitem{SBP_gassner}
G.~J. Gassner, A.~R.Winters, and D.~A. Kopriva.
\newblock {Split form nodal discontinuous Galerkin schemes with
  summation-by-parts property for the compressible Euler equations}.
\newblock {\em J. Comp. Phys.}, 327:39--66, 2016.

\bibitem{IDP_guermond_2018}
J.-L. Guermond, M.~Nazarov, B.~Popov, and I.~Tomas.
\newblock Second-order invariant domain preserving approximation of the euler
  equations using convex limiting.
\newblock {\em SIAM Journal on Scientific Computing}, 40(5):A3211--A3239, 2018.

\bibitem{IDP_guermond_2016}
J.-L. Guermond and B.~Popov.
\newblock Invariant domains and first-order continuous finite element
  approximation for hyperbolic systems.
\newblock {\em SIAM J. Numer. Anal.}, 54(4):2466--2489, 2016.

\bibitem{IDP_guermond_2019}
J.-L. Guermond, B.~Popov, and I.~Tomas.
\newblock Invariant domain preserving discretization-independent schemes and
  convex limiting for hyperbolic systems.
\newblock {\em Comput. Methods Appl. Mech. and Engrg.}, 347:143--175, 2019.

\bibitem{vilar_aplsc_NSW_1D}
A.~Haidar, F.~Marche, and F.~Vilar.
\newblock A posteriori finite-volume local subcell correction of high-order
  discontinuous {G}alerkin schemes for the nonlinear shallow-water equations.
\newblock {\em J. Comp. Phys.}, 452:110902, 2022.

\bibitem{vilar_aplsc_NSW_1D_fixed}
A.~Haidar, F.~Marche, and F.~Vilar.
\newblock Free-boundary problems for wave structure interactions in
  shallow-water: Dg-ale description and local subcell correction.
\newblock {\em J. Sci. Comput.}, 98(2), 2024.

\bibitem{monolithic_hajduk_2021}
H.~Hajduk.
\newblock Monolithic convex limiting in discontinuous galerkin discretizations
  of hyperbolic conservation laws.
\newblock {\em Computers and Mathematics with Applications}, 87:120--138, 2021.

\bibitem{Hesthaven_nodal_DG}
J.S. Hesthaven and T.~Warburton.
\newblock {\em Nodal Discontinuous Galerkin Methods: Algorithms, Analysis, and
  Applications}.
\newblock Springer Publishing Company, Incorporated, 2007.

\bibitem{unstructured_SBP_16}
J.E. Hicken, D.C. Del~Rey Fern\'{a}ndez, and D.W. Zingg.
\newblock Multidimensional summation-by-parts operators: General theory and
  application to simplex elements.
\newblock {\em SIAM J. Sci. Comput.}, 38(4):A1935--A1958, 2016.

\bibitem{hou}
S.~Hou and X.-D. Liu.
\newblock {Solutions of Multi-dimensional Hyperbolic Systems of Conservation
  Laws by Square Entropy Condition Satisfying Discontinuous Galerkin Method}.
\newblock {\em J. Sci. Comput.}, 31:127--151, 2007.

\bibitem{peraire_2012}
A.~Huerta, E.~Casoni, and J.~Peraire.
\newblock {A simple shock-capturing technique for high-order discontinuous
  Galerkin methods}.
\newblock {\em Int. J. Numer. Meth. Fluids}, 69:1614--1632, 2012.

\bibitem{MLP_unstruct}
{J. S. Park and S.-H. Yoon and C. Kim}.
\newblock {Multi-dimensional limiting process for hyperbolic conservation laws
  on unstructured grids}.
\newblock {\em J. Comp. Phys.}, 229:788--812, 2010.

\bibitem{jiang}
G.-S. Jiang and C.-W. Shu.
\newblock On cell entropy inequality for discontinuous galerkin method for a
  scalar hyperbolic equation.
\newblock {\em Mathematics of Computation}, 62:531--538, 1994.

\bibitem{jiang_weno}
G.-S. Jiang and C.-W. Shu.
\newblock Efficient implementation of weighted eno schemes.
\newblock {\em J. Comp. Phys.}, 126:202--228, 1996.

\bibitem{Kamm07}
J.R. Kamm and F.X. Timmes.
\newblock {On efficient generation of numerically robust Sedov solutions}.
\newblock Technical Report LA-UR-07-2849, Los Alamos National Laboratory, 2007.

\bibitem{kriv07}
L.~Krivodonova.
\newblock {Limiters for high-order discontinuous Galerkin methods}.
\newblock {\em J. Comp. Phys.}, 226:879--896, 2007.

\bibitem{Popov}
A.~Kurganov, G.~Petrova, and B.~Popov.
\newblock {Adaptive semi-discrete central-upwind schemes for non convex
  hyperbolic conservation laws}.
\newblock {\em SIAM J. Sci. Comput.}, 29:2381--2401, 2007.

\bibitem{Kuzmin09}
D.~Kuzmin.
\newblock {A vertex-based hierarchical slope limiter for p-adaptative
  discontinuous Galerkin methods}.
\newblock {\em J. Comp. Appl. Math.}, 233:3077--3085, 2009.

\bibitem{Kuzmin_monolithic_2020}
D.~Kuzmin.
\newblock Monolithic convex limiting for continuous finite element
  discretizations of hyperbolic conservation laws.
\newblock {\em Comput. Methods Appl. Mech. and Engrg.}, 361:112804, 2020.

\bibitem{Kuzmin_subcell_flux_limiting}
D.~Kuzmin and M.Q. {de Luna}.
\newblock Subcell flux limiting for high-order bernstein finite element
  discretizations of scalar hyperbolic conservation laws.
\newblock {\em J. Comp. Phys.}, 411:109411, 2020.

\bibitem{Kuzmin_FCT_limiting_2022}
D.~Kuzmin, M.Q. de~Luna, D.I. Ketcheson, and J.~Gr\"ull.
\newblock {Bound-preserving flux limiting for high-order explicit Runge Kutta
  time discretizations of hyperbolic conservation laws}.
\newblock {\em J. Sci. Comput.}, 91, 2022.

\bibitem{Kuzmin_smooth}
D.~Kuzmin and F.~Schieweck.
\newblock A parameter-free smoothness indicator for high-resolution finite
  element schemes.
\newblock {\em Centr. Eur. J. Math.}, 11:1478--1488, 2013.

\bibitem{LeVeque2}
R.~J. LeVeque.
\newblock {High-resolution conservative algorithms for advection in
  compressible flow}.
\newblock {\em SIAM J. Numer. Anal.}, 33:627--665, 1996.

\bibitem{Li_vertex_based_lim}
L.~Li and Q.~Zhang.
\newblock {A new vertex-based limiting approach for nodal discontinuous
  Galerkin methods on arbitrary unstructured meshes}.
\newblock {\em Computers and Fluids}, 159:316--326, 2017.

\bibitem{Chan_knapsack}
Y.~Lin and J.~Chan.
\newblock High order entropy stable discontinuous galerkin spectral element
  methods through subcell limiting.
\newblock {\em J. Comp. Phys.}, 498:112677, 2024.

\bibitem{Kuzmin_FCT_limiting_2017}
C.~Lohmann, D.~Kuzmin, J.N. Shadid, and S.~Mabuza.
\newblock {Flux-corrected transport algorithms for continuous Galerkin methods
  based on high order Bernstein finite elements}.
\newblock {\em J. Comp. Phys.}, 344:151--186, 2017.

\bibitem{Osher_Eflux}
S.~Osher.
\newblock {Riemann solvers, the entropy condition and difference
  approximations}.
\newblock {\em SIAM J. Numer. Anal.}, 21:217--235, 1984.

\bibitem{Pazner_2021}
W.~Pazner.
\newblock Sparse invariant domain preserving discontinuous galerkin methods
  with subcell convex limiting.
\newblock {\em Comput. Methods Appl. Mech. and Engrg.}, 382:113876, 2021.

\bibitem{Reed}
W.~H. Reed and T.~R. Hill.
\newblock {Triangular Mesh Methods for the Neutron Transport Equation}.
\newblock Technical Report LA-UR-73-479, Los Alamos National Laboratory, 1973.

\bibitem{Rodionov_half_cyl}
A.V. Rodionov.
\newblock {Artificial viscosity Godunov-type schemes to cure the carbuncle
  phenomenon}.
\newblock {\em J. Comp. Phys.}, 345:308--329, 2017.

\bibitem{Kuzmin_Gassner_monolithic_2024}
A.~Rueda-Ram\`irez, B.~Bolm, D.~Kuzmin, and G.~Gassner.
\newblock Monolithic convex limiting for legendre-gauss-lobatto discontinuous
  galerkin spectral-element methods.
\newblock {\em Commun. Appl. Math. Comput.}, 2024.

\bibitem{subcell_limiting_rueda}
A.M. Rueda-Ram\`irez, W.~Pazner, and G.J. Gassner.
\newblock Subcell limiting strategies for discontinuous galerkin spectral
  element methods.
\newblock {\em Computers and Fluids}, 247:105627, 2022.

\bibitem{subcell_shock_rueda}
{S. Hennemann and A.M. Rueda-Ram\`irez and F.J. Hindenlang and G.J. Gassner}.
\newblock A provably entropy stable subcell shock capturing approach for high
  order split form dg for the compressible euler equations.
\newblock {\em J. Comp. Phys.}, 426:109935, 2021.

\bibitem{Osher}
C.-W. Shu and S.~Osher.
\newblock {Efficient implementation of essentially non-oscillatory
  shock-capturing schemes}.
\newblock {\em J. Comp. Phys.}, 77:439--471, 1988.

\bibitem{Sod}
G.~A. Sod.
\newblock A survey of several finite difference methods for systems of
  non-linear hyperbolic conservation laws.
\newblock {\em J. Comp. Phys.}, 27:1--31, 1978.

\bibitem{munz_subFV_2014}
M.~Sonntag and C.~D. Munz.
\newblock {Shock capturing for discontinuous Galerkin methods using finite
  volume subcells}.
\newblock In {\em Finite Volumes for Complex Applications VII}, pages 945--953.
  Springer, 2014.

\bibitem{Tadmor_Escheme}
E.~Tadmor.
\newblock Numerical viscosity and the entropy condition for conservative
  difference schemes.
\newblock {\em Mathematics of Computation}, 168(43):369--381, 1984.

\bibitem{Tadmor_entropy}
E.~Tadmor.
\newblock The numerical viscosity of entropy stable schemes for systemes of
  conservation laws. i.
\newblock {\em Mathematics of Computation}, 179(49):91--103, 1987.

\bibitem{Tadmor_entropy_acta}
E.~Tadmor.
\newblock Entropy stability theory for difference approximations of nonlinear
  conservation laws and related time-dependent problems.
\newblock {\em Acta Numerica}, 12:451--512, 2003.

\bibitem{toro_book}
E.~F. Toro.
\newblock {\em {Riemann solvers and numerical methods for fluid dynamics}}.
\newblock Springer-Verlag, 1999.

\bibitem{vilar_aplsc_1D}
F.~Vilar.
\newblock {A posteriori correction of high-order discontinuous Galerkin scheme
  through subcell finite volume formulation and flux reconstruction}.
\newblock {\em J. Comp. Phys.}, 387:245--279, 2018.

\bibitem{vilar_aplsc_2D}
F.~Vilar and R.~Abgrall.
\newblock {A Posteriori Local Subcell Correction of High-Order Discontinuous
  Galerkin Scheme for Conservation Laws on Two-Dimensional Unstructured Grids}.
\newblock {\em SIAM J. Sci. Comput.}, 46(2):A851--A883, 2024.

\bibitem{vilar_lag_10}
F.~Vilar, P.-H. Maire, and R.~Abgrall.
\newblock {Cell-centered discontinuous Galerkin discretizations for
  two-dimensional scalar conservation laws on unstructured grids and for
  one-dimensional Lagrangian hydrodynamics}.
\newblock {\em Computers and Fluids}, 46(1):498--604, 2011.

\bibitem{vilar_lag_14}
F.~Vilar, P.-H. Maire, and R.~Abgrall.
\newblock {A discontinuous Galerkin discretization for solving the
  two-dimensional gas dynamics equations written under total Lagrangian
  formulation on general unstructured grids}.
\newblock {\em J. Comp. Phys.}, 276:188--234, 2014.

\bibitem{woodward_1984}
P.~Woodward and P.~Collela.
\newblock {The numerical-simulation of two-dimensional fluid-flow with strong
  shocks}.
\newblock {\em J. Comp. Phys.}, 54(1):115--173, 1984.

\bibitem{bound_preserving_shu}
K.~Wu and C.-W. Shu.
\newblock Geometric quasilinearization framework for analysis and design of
  bound-preserving schemes.
\newblock {\em SIAM Review}, 65(4):1031--1073, 2023.

\bibitem{MYang}
M.~Yang and Z.J. Wang.
\newblock {A parameter-free generalized moment limiter for high-order methods
  on unstrucured grids}.
\newblock {\em Adv. Appl. Math. Mech.}, 4:451--480, 2009.

\bibitem{Zalesak_FCT}
S.T. Zalesak.
\newblock {Fully multidimensional flux-corrected transport algorithms for
  fluids}.
\newblock {\em J. Comp. Phys.}, 31(3):335--362, 1979.

\bibitem{zshu1}
X.~Zhang and C.-W. Shu.
\newblock {On maximum-principle-satisfying high order schemes for scalar
  conservation laws}.
\newblock {\em J. Comp. Phys.}, 229:3091--3120, 2010.

\bibitem{zshu5}
X.~Zhang, Y.~Xia, and C.-W. Shu.
\newblock Maximum-principle-satisfying and positivity-preserving high order
  discontinuous galerkin schemes for conservation laws on triangular meshes.
\newblock {\em J. Sci. Comput.}, 50:29--62, 2012.

\end{thebibliography}
\end{document}